\documentclass[11pt]{article}

\usepackage{amsmath, amsthm, amssymb}
\usepackage{enumitem}
\usepackage{pdflscape}
\usepackage{caption}
\usepackage{bm}
\usepackage{ifpdf}
	\ifpdf
\usepackage[pdftex]{graphicx}
	\else
\usepackage[dvips]{graphicx}
	\fi
\usepackage[all]{xy}
\usepackage{tocvsec2}
\usepackage{bbm}
	\input xy
	\xyoption{all}
\usepackage[pdftex,plainpages=false,hypertexnames=false,pdfpagelabels]{hyperref}
	\newcommand{\arxiv}[1]{\href{http://arxiv.org/abs/#1}{\tt arXiv:\nolinkurl{#1}}}
	\newcommand{\arXiv}[1]{\href{http://arxiv.org/abs/#1}{\tt arXiv:\nolinkurl{#1}}}
	\newcommand{\doi}[1]{\href{http://dx.doi.org/#1}{{\tt DOI:#1}}}

	\newcommand{\googlebooks}[1]{(preview at \href{http://books.google.com/books?id=#1}{google books})}
	
\usepackage{xcolor}
	\definecolor{dark-red}{rgb}{0.7,0.25,0.25}
	\definecolor{dark-blue}{rgb}{0.15,0.15,0.55}
	\definecolor{medium-blue}{rgb}{0,0,.8}
	\definecolor{DarkGreen}{RGB}{0,150,0}
	\definecolor{rho}{named}{red}
	\hypersetup{
	   colorlinks, linkcolor={purple},
	   citecolor={medium-blue}, urlcolor={medium-blue}
	}
\usepackage{longtable}
\usepackage{fullpage}

\usepackage{mathbbol}

\theoremstyle{plain}
\newtheorem{thm}{Theorem}[section]
\newtheorem*{thm*}{Theorem}
\newtheorem{thmalpha}{Theorem}

\newtheorem{cor}[thm]{Corollary}
\newtheorem{coralpha}[thmalpha]{Corollary}
\newtheorem*{cor*}{Corollary}

\newtheorem*{conj*}{Conjecture}
\newtheorem{lem}[thm]{Lemma}
\newtheorem{prop}[thm]{Proposition}

\newtheorem*{quest*}{Question}
\newtheorem*{claim*}{Claim}

\theoremstyle{definition}
\newtheorem{defn}[thm]{Definition}

\newtheorem{construction}[thm]{Construction}

\newtheorem{ex}[thm]{Example}
\newtheorem*{ex*}{Example}
\newtheorem{sub-ex}[thm]{Sub-Example}
\newtheorem{counter-ex}[thm]{Counter-Example}

\newtheorem*{rem*}{Remark}
\newtheorem{remark}[thm]{Remark}

\DeclareMathOperator{\Ad}{Ad}

\DeclareMathOperator{\coev}{coev}

\DeclareMathOperator{\End}{End}
\DeclareMathOperator{\ev}{ev}
\DeclareMathOperator{\Hom}{Hom}

\DeclareMathOperator{\op}{op}

\DeclareMathOperator{\supp}{supp}
\DeclareMathOperator{\id}{id}
\DeclareMathOperator{\Isom}{Isom}

\DeclareMathOperator{\Irr}{Irr}

\DeclareMathOperator{\Tr}{Tr}
\DeclareMathOperator{\tr}{tr}

\newcommand{\comment}[1]{}
\newcommand{\set}[2]{\left\{#1 \middle| #2\right\}}

\newcommand{\bbOne}{\mathbbm{1}}

\newcommand{\noshow}[1]{}
\newcommand{\MR}[1]{}

\newcommand{\Rep}{{\sf Rep}}

\newcommand{\Mod}{{\sf Mod}}
\newcommand{\Proj}{{\sf Proj}}

\newcommand{\Bim}{{\sf Bim}}

\newcommand{\fgpBim}{{\sf Bim_{fgp}}}
\newcommand{\fgpMod}{{\sf Mod_{fgp}}}
\newcommand{\fgpBimtr}{{\sf Bim_{fgp}^{tr}}}
\newcommand{\spbfBim}{{\sf Bim_{bf}^{sp}}}
\renewcommand{\Vec}{{\sf Vec}}

\newcommand{\Hilb}{{\sf Hilb}}
\newcommand{\fdHilb}{{\sf Hilb_{fd}}}
\newcommand{\rCorr}{{\mathsf{C^{*}Alg}}}

\newcommand{\cCCAlgs}{{\sf C^*Alg}(\cC)}
\newcommand{\CDisc}{{\sf C^*Disc}}

\newcommand{\PQR}{{\sf PQR}}
\newcommand{\PQN}{{\sf PQN}}
\newcommand{\<}{\langle}
\renewcommand{\>}{\rangle}

\def\semicolon{;}
\def\applytolist#1{
    \expandafter\def\csname multi#1\endcsname##1{
        \def\multiack{##1}\ifx\multiack\semicolon
            \def\next{\relax}
        \else
            \csname #1\endcsname{##1}
            \def\next{\csname multi#1\endcsname}
        \fi
        \next}
    \csname multi#1\endcsname}

\def\calc#1{\expandafter\def\csname c#1\endcsname{{\mathcal #1}}}
\applytolist{calc}QWERTYUIOPLKJHGFDSAZXCVBNM;
\def\bbc#1{\expandafter\def\csname bb#1\endcsname{{\mathbb #1}}}
\applytolist{bbc}QWERTYUIOPLKJHGFDSAZXCVBNM;
\def\bfc#1{\expandafter\def\csname bf#1\endcsname{{\mathbf #1}}}
\applytolist{bfc}QWERTYUIOPLKJHGFDSAZXCVBNM;
\def\sfc#1{\expandafter\def\csname s#1\endcsname{{\sf #1}}}
\applytolist{sfc}QWERTYUIOPLKJHGFDSAZXCVBNM;
\def\fc#1{\expandafter\def\csname f#1\endcsname{{\mathfrak #1}}}
\applytolist{fc}QWERTYUIOPLKJHGFDSAZXCVBNM;

\usepackage{tikz}
\usepackage{tikz-cd}
\usetikzlibrary{arrows,backgrounds,patterns.meta}
\usetikzlibrary{positioning,shadings,cd}
\usetikzlibrary{shapes}
\usetikzlibrary{backgrounds}
\usetikzlibrary{decorations,decorations.pathreplacing,decorations.markings}
\usetikzlibrary{fit,calc,through}
\usetikzlibrary{external}
\usetikzlibrary{arrows}
\tikzset{vertex/.style = {shape=circle,draw,fill=black,inner sep=0pt,minimum size=5pt}}
\tikzset{edge/.style = {->,> = latex', bend right}}
\tikzset{
	super thick/.style={line width=3pt}
}
\tikzset{
    quadruple/.style args={[#1] in [#2] in [#3] in [#4]}{
        #1,preaction={preaction={preaction={draw,#4},draw,#3}, draw,#2}
    }
}
\tikzstyle{shaded}=[fill=red!10!blue!20!gray!30!white]
\tikzstyle{unshaded}=[fill=white]
\tikzstyle{empty box}=[circle, draw, thick, fill=white, opaque, inner sep=2mm]
\tikzstyle{annular}=[scale=.7, inner sep=1mm, baseline]
\tikzstyle{rectangular}=[scale=.75, inner sep=1mm, baseline=-.1cm]
\tikzstyle{mid>}=[decoration={markings, mark=at position 0.5 with {\arrow{>}}}, postaction={decorate}]
\tikzstyle{mid<}=[decoration={markings, mark=at position 0.5 with {\arrow{<}}}, postaction={decorate}]
\tikzstyle{over}=[double, draw=white, super thick, double=]

\tikzdeclarepattern{
  name=primeddots,
  type=uncolored,
  bounding box={(-.6pt,-.6pt) and (.6pt,.6pt)},
  tile size={(5pt,5pt)},
  tile transformation={rotate=60},
  parameters={none}, 
  code={
    \fill(0pt,0pt) circle (.35pt);
  }
}

\tikzdeclarepattern{
  name=primeddots2,
  type=uncolored,
  bounding box={(-.3pt,-.3pt) and (.3pt,.3pt)},
  tile size={(2.5pt,2.5pt)},
  tile transformation={rotate=60},
  parameters={none}, 
  code={
    \fill(0pt,0pt) circle (.35pt);
  }
}

\tikzstyle{primedregion}[none]=[
	preaction={fill=#1},
	pattern=primeddots,
  draw=#1,
]

\tikzstyle{primedregion2}[none]=[
	preaction={fill=#1},
	pattern=primeddots2,
  draw=#1,
]

\newcommand{\roundNbox}[6]{
	\draw[rounded corners=5pt, very thick, #1] ($#2+(-#3,-#3)+(-#4,0)$) rectangle ($#2+(#3,#3)+(#5,0)$);
	\coordinate (ZZa) at ($#2+(-#4,0)$);
	\coordinate (ZZb) at ($#2+(#5,0)$);
	\node at ($1/2*(ZZa)+1/2*(ZZb)$) {#6};
}

\newcommand{\tikzmath}[2][]
     {\vcenter{\hbox{\begin{tikzpicture}[#1]#2
                     \end{tikzpicture}}}
     }

\newcommand{\BColor}{gray!55}

\tikzcdset{scale cd/.style={every label/.append style={scale=#1},
    cells={nodes={scale=#1}}}}

\let\OLDthebibliography\thebibliography
\renewcommand\thebibliography[1]{
  \OLDthebibliography{#1}
  \setlength{\parskip}{0pt}
  \setlength{\itemsep}{0pt plus 0.3ex}
}

\begin{document}

\title{\textbf{Discrete Inclusions of C*-algebras}}
\author{Roberto Hern\'{a}ndez Palomares\footnote{Department of Mathematics, University of Waterloo \hfill \url{robertohp.math@gmail.com}} and Brent Nelson\footnote{Department of Mathematics, Michigan State University \hfill \url{brent@math.msu.edu}}}
\date{}
\maketitle
\begin{abstract}
\noindent We introduce the category of C*-discrete inclusions of C*-algebras $A\subset B$ with a faithful conditional expectation $E:B\twoheadrightarrow A$. 
This class includes many examples such as finite Watatani index inclusions, and also abundant infinite index inclusions like crossed products by outer actions of discrete (quantum) groups and unitary tensor categories. We prove irreducible $(A'\cap B= \bbC1)$ C*-discrete inclusions are precisely crossed products by outer actions of unitary tensor categories and certain C*-algebra objects. 
\end{abstract}

\section*{Introduction}

V. F. R. Jones started the modern theory of subfactors with his extraordinary discovery of the Index Rigidity Theorem for $\rm{II}_1$-subfactors \cite{MR696688}: the range of indices is exactly $\{4\cos^2\left(\pi/n\right)\}_{n\geq3}\cup[4,\infty].$ As the theory matured it became clear that subfactors have \emph{quantum symmetries}---vastly generalizing discrete groups and their representations---and that a huge class, the \emph{discrete subfactors}, could be understood in terms of actions by unitary tensor categories (UTCs). In the setting of C*-algebras, subfactors are replaced by an inclusion of C*-algebras $A\subset B$ admitting a faithful conditional expectation $E\colon B\twoheadrightarrow A$, which we denote by $A\overset{E}{\subset} B$. The Index Rigidity Theorem was adapted to the setting of C*-algebras by Watatani in \cite{MR996807}, along with many other tools from subfactor theory such as \cite{ MR1228532, MR1615672,  MR1604162, MR1642530,  MR1604162, MR1624182, MR1742862, MR1862184,  MR1900138, MR1990630, MR2085108, MR2125398, MR2590623, MR3145747, MR4010423, MR4083877, AranoRokhlin, 2021arXiv210511899R, MR4419534, 2022arXiv220711854C}. We follow in this tradition by studying the C*-analogue of discrete subfactors, taking inspiration from the recent work of C. Jones and Penneys \cite{MR3948170}.

Discreteness was first introduced for semifinite subfactors in \cite{MR1055223}, and was later refined to be suitable for general subfactors in \cite[Definition 3.7]{MR1622812}. This latter definition states that a subfactor $N\subset M$ admits a faithful normal conditional expectation $E\colon M\twoheadrightarrow N$ whose dual operator valued weight from the basic construction $\<M,e_N\>$ to $M$ is semifinite on the relative commutant $N'\cap \<M,e_N\>$. Adapting this definition to C*-algebras presents many challenges, since several of the concepts involved either lack analogues in the C*-setting or behave erratically in the infinite index setting \cite{MR1642530, MR2085108}. Instead, we appeal to a more intrinsic characterization of discreteness offered by \cite[Proposition 3.22]{MR3948170}, which states that for certain irreducible subfactors discreteness is equivalent to the \emph{quasi-normalizer} of the inclusion being weakly dense---a condition that was studied extensively in \cite{MR1729488, MR3801484}. Using this characterization, \cite{MR3948170} introduced a tensor category description of discrete subfactors, focusing  on actions of UTCs on $\rm{II}_1$-factors and highlighting the importance of so called \emph{W*-algebra objects}, which are infinite dimensional analogues of Q-systems \cite{JP17}.  Roughly, they showed that underlying any irreducible discrete subfactor  $N\subset M$ there is an outer action of some UTC on $N$, denoted $\cC\overset{F}{\curvearrowright}N,$ and a W*-algebra object $\bbM$. Formally, $F:\cC\to\fgpBim(N)$ is a fully-faithful unitary tensor functor (i.e. bijective on morphisms) and $\bbM:\cC^{\op}\to\Vec$ is a linear functor, which we also refer to as a $\cC$-graded W*-algebra. Thus, we refer to the image $F[\cC]$ as \emph{the support of the inclusion}, or $\cC_{N\subset M}$ for clarity.  This data allows one to parameterize discrete subfactors over $N$ in terms of the dynamical data $\cC\overset{F}{\curvearrowright}N$ and the algebraic data $(\cC,\bbM)$. These ideas can be translated to the setting of C*-algebras and underpin our main results and the definition of C*-discreteness, which we now endeavor to explain.

For a unital inclusion of C*-algebras with expectation $A\overset{E}{\subset} B$, let $\cB$ denote the $B$-$A$ right C*-correspondence generated by $B$ equipped with the $A$-valued inner product $\<b_1 | b_2 \>_A = E(b_1^* b_2)$, and let $\fgpBim(A)$ denote the collection of finitely generated projective (fgp) Hilbert $A$-$A$-bimodules (see Section~\ref{sec:modules_correspondences_bimodules}). We denote by $B\Omega$ the dense copy of $B$ in $\cB$, and notice the inclusion has finite Pimsner--Popa index if and only if $\cB=B\Omega$ \cite[Theorem 1]{MR1642530}. Due to the fact that we work with Hilbert C*-modules rather than Hilbert spaces, the quasi-normalizer in the subfactor setting is replaced by the following: the \emph{projective quasi-normalizer} of $A\overset{E}{\subset} B$ is defined as the intermediate $*$-algebra $\PQN(A\subset B)$ consisting of the $b\in B$ such that there exists a $K\in \fgpBim(A)$ satisfying $b\in K\subset B\Omega \subset \cB$. Therefore, the $*$-algebra structure of $\PQN(A\subset B)$ is determined---and in fact dominated---by the fusion algebra of some UTC. We then define a \emph{C*-discrete} inclusion $A\overset{E}{\subset}B$ to be one where $\PQN(A\subset B)$ is norm-dense in $B$ (see Definition~\ref{defn:CDisc}).
C*-discrete inclusions assemble into a category with expectation preserving ucp morphisms, denoted $\CDisc$ (see Definition~\ref{defn:CStarDisc}).

In our context, $\fgpBim(A)$ forms a UTC consisting of dualizable right $A$-$A$ correspondences (see Section~\ref{subsec:category_of_C*-alg}), and so an outer action of a UTC on $A$ is a fully-faithful unitary tensor functor $F:\cC\to \fgpBim(A)$. 
Given $A\subset B,$ we denote by $\cC_{A\subset B}$ the UTC generated by the dualizable $A$-$A$ sub-bimodules of $\cB$.
Additionally, under the circumstances that we will consider, \emph{C*-algebra objects} (a categorification of C*-algebras) simply coincide with the W*-algebra objects (see Section \ref{sec:AlgObjs}).
The following is then the state-free C*-analogue of \cite[Theorem A]{MR3948170}:
\begin{thmalpha}[{Theorem~\ref{thm:DiscCharacterization}}]\label{thmalpha:Main}
    Let $A$ be a unital C*-algebra with trivial center, and $F:\cC\to \fgpBim(A)$ be an outer action by a UTC.  We then have a correspondence: 
     \begin{align*}
    &\left\{
        \begin{aligned}
            &\left(A\overset{E}{\subset}B\right)\in\CDisc\ \\
        \end{aligned}
         \middle|\
        \begin{aligned}
              & A'\cap B =\bbC1\\
              &\ \  \cC_{A\subset B}\subseteq F[\cC] 
        \end{aligned}
    \right\}\ \ \
     \cong&&
    \left\{
        \begin{aligned}
            & \text{C*-algebra object } \bbB\text{ in } \cC\\
        \end{aligned}
         \middle|\
          \begin{aligned}
              & \bbB(1_\cC)\cong \bbC\\
        \end{aligned}
    \right\}\!.
    \end{align*}
\end{thmalpha}

In the above correspondence, the map from left to right is essentially taking the projective quasi-normalizer of the inclusion, while the map from right to left is taking the \emph{reduced crossed product} by a connected $\cC$-graded C*-algebra $\bbB$ through the action $\cC\overset{F}{\curvearrowright}A$ (Construction \ref{const:realizedredC*alg}). For the latter, we adapt the \emph{realization}/crossed product and \emph{underlying algebra} functors of C. Jones and Penneys to our setting. The heart of this proof involves showing that $\PQN(A\subset B)$ is a *-algebra (see Lemma~\ref{lem:deltaproperties}), which demands heavy use of the abstract machinery of \emph{operator algebras in UTCs} developed in \cite{JP17, MR3948170}.

It follows as a direct corollary to Theorem~\ref{thmalpha:Main} that every irreducible C*-discrete inclusion arises as a crossed product by an outer action of a UTC.

\begin{coralpha}[{Corollary~\ref{cor:reconstruction}}]\label{coralpha:StdInvs}
Let $\left( A\overset{E}{\subset} B\right)\in\CDisc$ be unital and irreducible. Then 
    $$\left(A\overset{E}{\subset} B\right) \cong \left( A\overset{E'}{\subset} A\underset{F,r}{\rtimes} \langle B\rangle\right).$$
    Here, $E'$ is the canonical faithful conditional expectation mapping elements in the crossed product onto the coefficient algebra $A$ and $F:\cC_{A\subset B}\hookrightarrow \fgpBim(A)$ is the inclusion. 
\end{coralpha}

Despite \cite{MR3948170} offering a clear road map to Theorem~\ref{thmalpha:Main}, there are nonetheless many crucial differences between the von Neumann and C*-algebraic settings to overcome. For example, the standard form $L^2(M,\tau_M)$ offers a canonical right action, and this is stratum for the \emph{Frobenius Reciprocity Theorem} \cite[Theorem 3.13]{MR3948170} used to unveil the structure of the quasi-normalizer in terms of $\cC_{N\subset M}$-graded W*-algebras. This is in stark contrast to $\cB$ which admits \emph{no} bounded right $B$-action. Nevertheless, passing back and forth between the $B$-$A$ correspondence $\cB$ and $B$ with its trivial $B$-$B$ structure, one can completely bypass this analytic obstruction and arrive at a \emph{C*-Frobenius Reciprocity Theorem} (see Theorem~\ref{thm:FR}). This can then be used to describe the projective quasi-normalizer in terms of $\cC_{A\subset B}$-graded C*-algebras; namely,
$$
\PQN(A\subset B)\cong \bigoplus_{K\in \Irr(\cC_{A\subset B})} K\otimes_\bbC \Hom_{A-A}(K\to \cB)^\diamondsuit,
$$ 
where $\Irr(\cC)$ is any complete list of representatives of irreducible fgp $A$-$A$ bimodules (see Proposition \ref{prop:PQNAlgRealization}). 
The symbol $\diamondsuit$ means we only take those adjointable $A$-$A$ intertwiners valued in $B\Omega\subset \cB,$ and we refer to these as the ``diamond spaces,'' which are categorical analogues of \emph{$A$-central vectors}.
Equipped with the C*-Frobenius Reciprocity Theorem, we can link the diamond spaces with the \emph{cyclic $\cC$-module C*-category} generated by $B\in\fgpBim(B)$, which by a unitary version of \emph{Ostrik Theorem} gives a $\cC$-graded C*-algebra object \cite[Theorem 2]{JP17}.

A key tool to understand a subfactor is its standard invariant, formally  given by the lattice of higher relative commutants. Depending on the context,  it has been axiomatized in different ways such as  Ocneanu's \emph{paragroups} for finite-depth subfactors \cite{MR996454, MR1642584}, more generally as Popa's $\lambda$-\emph{lattices} \cite{MR1334479}, Jones' \emph{Planar algebras} \cite{1999math......9027J}, and for type $\rm{III}$-factors  Longo's Q-systems \cite{MR1257245}, later given a categorical description by M\"uger in \cite{MR1966524}. The planar algebraic approach was notoriously effective in the classification of \emph{finite-index} $\rm{II}_1$-subfactors \cite{MR3166042,MR3345186,MR4565376}. The algebraic data $(\cC,\bbM)$ appearing in \cite{MR3948170} offers yet another axiomatization of the standard invariant, and it was shown \cite[Corollary C]{MR3948170} that any $(\cC, \bbM)$ always arises from some discrete subfactor as a crossed product by the diagrammatic GJS-action \cite{GJS10,JSW10}. This should be compared to Popa's results for subfactors \cite{MR1055708, MR1334479,MR1339767}. Motivated by this and Theorem \ref{thmalpha:Main}, we therefore define the \emph{standard invariant} of $A\subset B$ to be:
\begin{itemize}
    \item the support of the inclusion $\cC_{A\subset B}$ and the canonical embedding $F:\cC_{A\subset B}\to \fgpBim(A);$ 
    \item $\bbB,$ a connected $\cC_{A\subset B}$-graded C*-algebra  (i.e. $\bbB(1_{\cC_{A\subset B}})\cong\bbC).$
\end{itemize}
The following is the C*-analogue of \cite[Corollary C]{MR3948170}. It relies on the GJS-construction as expressed by Hartglass and the first-named author \cite{MR4139893} to obtain a unital simple C*-algebra $A=A(\cC)$ admitting an outer action $\cC\curvearrowright A$.

\begin{coralpha}[{Corollary~\ref{cor:RedReconstruction}}]\label{coralpha:RedReconstruction}
Let $(\cC,\bbB)$ be an abstract standard invariant. Then, there exists $(A\subset B)\in\CDisc$ whose standard invariant is isomorphic to $(\cC,\bbB).$ Moreover, $A$ can be chosen to be separable, unital, simple, stable rank 1, monotracial, and exact.
\end{coralpha}

Examples of (irreducible) C*-discrete inclusions are ubiquitous and can also be produced from actions of UTCs and their internal structure. 
On the one hand, besides outer actions by discrete groups, there are numerous \emph{genuine} examples of UTC-actions $F$ on various C*-algebras in the literature which are not \emph{group-automorphic}  \cite{MR1211581,MR1228532,MR1604162,MR4139893, MR4328058,2022arXiv220711854C, MR4566007}. 
On the other hand, producing $\cC$-graded C*-algebras is entirely a problem internal to $\cC$, independent of the dynamical ingredient $F$.  Notoriously, any irreducible discrete subfactor yields examples of connected $\cC$-graded C*-algebras, since locally finite dimensional $\cC$-graded W*-algebras are automatically C*. Very importantly, any compact quantum group $\bbG$ and its \emph{Fiber functors} $\Rep_f(\bbG\to\fdHilb)$ from Tannaka-Krein duality also furnish examples of connected $\cC$-graded C*-algebras \cite{JP17}. Other fundamental families of C*-discrete inclusions stem from Cuntz algebra cores by their canonical \emph{gauge symmetries}, as well as from $A$-valued semicircular systems, and we will elaborate on these in forthcoming articles  \cite{HPNIII} and \cite{HPNII}. 

\subsection*{Acknowledgements}
The authors are grateful to David Penneys, who introduced them to the themes of this paper and provided many helpful comments along its development, and also to Corey Jones for offering useful suggestions, viewpoints and valuable advice. We are indebted to Matthew Lorentz who provided tremendous support in the initial stages of this project.
We thank also Michael Brannan, Michael Frank, Sergio Gir\'on Pacheco, Lucas Hataishi, Matthew Kennedy, David Larson, Bojan Magajna, Dimitri Shlyakhtenko, Roger Smith, Karen Strung, Reiji Tomatsu, Moritz Weber, Stuart White, and Makoto Yamashita for helpful conversations. 
The first-named author was partially supported by the AMS-Simons Travel Grant 2022, NSF grant DMS-2001163, NSF grant DMS-2000331 and NSF grant DMS-1654159. The second-named author was supported by NSF grant DMS-1856683.

\tableofcontents

\section{Preliminaries}

\subsection{Hilbert C*-(bi)modules and C*-correspondences}\label{sec:modules_correspondences_bimodules}

We recall the definitions of right (and left) Hilbert C*-modules, right (and left) C*-correspondences, Hilbert C*-bimodules, and related notions. Throughout, $A$ will denote a unital C*-algebra.

\subsubsection{Hilbert C*-modules}
A \textbf{right Hilbert $A$-module} is a vector space $X$ endowed with a right action of $A$, denoted by $\xi\lhd a$ for $\xi\in X$ and $a\in A$, and an $A$-valued inner product $\langle\ \cdot\mid -\rangle_A: X\times X\to A$ such that for each $\eta,\xi,\zeta\in X,$ and $a\in A$ we have: 
\begin{itemize}
    \item 
        $\langle \eta \mid \xi\lhd a +\zeta \rangle_A = \langle\eta\mid  \xi\rangle_A\cdot a + \langle\eta\mid  \zeta\rangle_A$, 
    \item 
        $\langle\eta \mid \xi\rangle_A = (\langle\xi\mid  \eta\rangle_A)^*,$
    \item 
        $\langle\xi \mid \xi\rangle_A\geq 0,$ and $\langle\xi\mid  \xi\rangle_A = 0$ if and only if $\xi=0.$
\end{itemize}
By the non-degeneracy of the inner product, we consider the canonical norm $||\xi||^2_A:=||\langle \xi \mid \xi\rangle_A||$ on $X$, which by assumption makes $(X,\ ||\cdot||_A)$ a Banach space. Some standard references on Hilbert C*-modules include \cite{wegge1993k,  MR1325694, MR2125398, 2006math.....11349F}.

For right Hilbert C* A-modules $X$ and $Y$, we consider right $A$-linear \textbf{adjointable operators} $T:X\to Y$, which together with some right $A$-linear map $T^*:Y\to X$ satisfy:  
    $$
        \langle T^*\eta\mid \xi \rangle^X_A = \langle \eta\mid T\xi\rangle^Y_A \qquad \xi\in X,\ \eta\in Y.
    $$
The space of all adjointable right $A$-linear maps is denoted by $\Hom^\dag_{-A}(X\to Y)$ or $\Hom^\dag(X_A\to Y_A),$ and $\End^\dag(X_A):=\Hom_A(X\to X).$
We remark that adjointable operators are necessarily bounded, and that $\End^\dag(X_A)$ is a unital C*-algebra.

There is an analogous notion of \textbf{left-$A$ Hilbert C*-modules} ${}_AX$, whose inner product is denoted by ${}_A\langle-,\ \cdot \rangle$, which contrastingly is $A$-linear on the left. Similarly, we denote by ${^\dag}\End({_A}X)$ the C*-algebra of all bounded endomorphisms of $X$ which are left $A$-modular and have an adjoint with respect to the left $A$-valued inner product.

Given a C*-algebra $A$, and a right Hilbert $A$-module $X_A,$ we say $X_A$ is \textbf{right finitely generated projective} (\textbf{fgp}) if there exists an $n\in\bbN$ and a right Hilbert $A$-module $Y_A$ such that $X\oplus Y\cong A^n$ as right Hilbert $A$-modules. 
Equivalently, for some projection $p\in M_n(A),$ we get $X_A\cong p[A^{\oplus n}]_A.$ 
    There is an analogous notion of left finitely generated left Hilbert $A$-module

A crucial characterization of right fgp Hilbert C*-modules is the existence of finite \textbf{right Pimsner--Popa basis} \cite{MR860811}.
This is, if $X_A\cong p[A^{\oplus n}]_A$ as above, by letting $u_i$ be the image of the columns of $p$ under the given isomorphism for $i=1,...,n,$ we obtain the following identity on $X$:
$$\xi=\sum_{i=1}^n u_i\lhd\langle u_i\mid  \xi \rangle_A.$$ 
For left fgp Hilbert C*-modules, there is an obvious analogous notion of a finite left Pimsner--Popa basis in terms of the left $A$-valued inner product:  
$$\xi =\overset{m}{\underset{j=1}{\sum}}{}_A\langle\xi, v_j\rangle\rhd v_j,$$ 
where $\{v_j\}_{j=1}^m\subset X$ is a finite \textbf{left Pimsner--Popa basis} obtained similarly.

\subsubsection{C*-correspondences}
Let $B$ be another unital C*-algebra. A \textbf{right $A$-$B$ C*-correspondence} ${}_AX_B$ 
consists of a right Hilbert $B$-module together with a unital *-homomorphism $A\to \End^\dag(X_B)$.
In practice, we refer to this unital *-homomorphism as a left $A$-action, denoted $a\rhd \xi$ for $a\in A$ and $\xi\in X,$ where the adjointability condition explicitly means that for every $\xi,\eta\in X$ and every $a\in A$ we have 
$$\langle\eta\mid  a\rhd\xi\rangle_B= \langle a^*\rhd\eta\mid  \xi\rangle_B.$$
Analogously, a \textbf{left $A$-$B$ C*-correspondence} is a left Hilbert $A$-module $Y$ with a unital $*$-homomorphism $B\to \End^\dag({}_A Y)$, which we refer to as a right $B$-action. Many of the following properties and constructions have analogous definitions for Hilbert C*-modules and left C*-correspondences, but given the focus of this paper we state them here for right C*-correspondences.

Given C*-algebras $A,B$ and a right C* correspondence ${}_A X_B$, we say $\{\xi_i\in X\colon i\in I\}$ \textbf{algebraically generates} $X$ if
    \[
        X=A\rhd\{\xi_i\}_{i\in I}\lhd B := \left\{\sum_{i\in F} a_i\rhd\xi_i \lhd b_i\bigg|\ F\subset I \text{ finite, }a_i\in A, b_i\in B\right\}.
    \]
For example, if $X$ is fgp as a right Hilbert $B$-module then any right finite Pimsner--Popa basis will generate $X$ algebraically.  
Furthermore, we say $\{\xi_i\}_{i\in I}\subset X$  \textbf{topologically generates} $X$ provided that $X=\overline{A\rhd\{\xi_i\}_i \lhd B}^{\|\cdot\|_B}$.

There are two types of tensor products one can define for right C*-correspondences (\emph{interior} and \emph{exterior}) and we will have occassion to use both.
Given right C*-correspondences ${}_AX_B$ and ${}_BY_C$, let $N$ denote the subspace of the algebraic tensor product $X\odot Y$ spanned by elements of the form $(\xi\lhd b)\boxtimes \eta - \xi \boxtimes (b\rhd \eta)$ for $\xi\in X$, $\eta\in Y$, and $b\in B$. The \textbf{interior tensor product} $X\boxtimes_B Y$ is the completion of $(X\odot Y)/N$ with respect to the $C$-valued inner product
    \[
        \<\xi_1\boxtimes \eta_1\mid \xi_2\boxtimes \eta_2\>_C:= \<\eta_1 \mid \<\xi_1 \mid\xi_2\>_B\rhd \eta_2\>_C,
    \]
which inherits the structure of a right $A$-$C$ C*-correspondence from ${}_AX$, $Y_C$ in the obvious way. When interior algebra $B$ is clear from context, we will suppress it in the subscript: $X\boxtimes Y$. (Note that this notation is not standard in the literature, but we adopt it here in order to avoid visual confusion with our next inner product.)

Given right C*-correspondences ${}_AX_B$ and ${}_CY_D$, the following defines a $(B\otimes_{\text{min}} D)$-valued inner product on $X\odot Y$:
    \[
        \<\xi_1\otimes \eta_1\mid \xi_2\otimes \eta_2 \>_{B\otimes_{\text{min}} D}:= \<\xi_1 \mid \xi_2\>_B \otimes \<\eta_1\mid \eta_2\>_D.
    \]
The \textbf{exterior tensor product} $X\otimes Y$ is the completion of $X\odot Y$ with respect to this inner product, which inherits the structure of a right $(A\otimes_{\text{min}} C)$-$(B\otimes_{\text{min}} D)$ C*-correspondence from ${}_AX_B$, ${}_CY_D$ in the obvious way. Note that no separation is required since this inner product is already positive definite (see \cite[Chapter 6]{MR1325694}, specifically pages 62-63).

\begin{remark}\label{rem:ext_tens_with_Hilb}
If we regard a Hilbert space $\cH$ as a right $\bbC$-$\bbC$ C*-correspondence, then for any right C*-correspondence ${}_AX_B$ the exterior tensor product $X\otimes \cH$ is a right $A$-$B$ C*-correspondence since $A\otimes_{\text{min}} \bbC\cong A$ and $B\otimes_{\text{min}} \bbC \cong B$. Explicitly, the left and right actions are given by
    \[
        a\rhd (\xi\otimes \eta)\lhd b = (a\rhd \xi \lhd b)\otimes \eta,
    \]
for $a\in A$, $b\in B$, $\xi\in X$, and $\eta\in \cH$, and the $B$-valued inner product is given by
    \[
        \< \xi_1\otimes \eta_1 \mid \xi_2\otimes \eta_2\>_A = \<\xi\mid \xi_2\>_A (\eta_1 \mid \eta_2)_{\bbC},
    \]
for $\xi_1,\xi_2\in X$ and $\eta_1,\eta_2\in \cH$.
\end{remark}

Given a family $\{ {}_AX^{(i)}_B\colon i\in I\}$ of right C*-correspondences, their \textbf{direct sum} is denoted
    \[
        \overline{\bigoplus_{i\in I}}\ X^{(i)}:= \left\{ (\xi_i)_{i\in I} \in \prod_{i\in I} X^{(i)}\colon \sum_{i\in I} \<\xi_i\mid \xi_i\>_B \text{ converges in B}\right\}.
    \]
To be precise, above we mean that the net of finite partial sums converges in norm. For such tuples
    \[
        \< (\xi_i)_{i\in I}\mid (\eta_i)_{i\in I}\>_B:= \sum_{i\in I} \<\xi_i\mid \eta_i\>_B
    \]
gives a right $B$-valued inner product and the actions
    \[
        a\rhd (\xi_i)_{i\in I}\lhd b:= (a\rhd \xi_i \lhd b)_{i\in I} \qquad a\in A,\ b\in B,
    \]
make the direct sum into a right $A$-$B$ C*-correspondence.
It is routine to check the left action of $A$ is by adjointable operators on the right $B$-module $(\overline{\bigoplus_i}X^{(i)})_B.$ If $X^{(i)}=X$ for all $i\in I$, we denote the direct sum by $\ell^2(I,X)$ and note that it is isomorphic to the exterior tensor product $X\otimes \ell^2(I)$.  

We will denote the \textbf{algebraic direct sum} by
    \[
        \bigoplus_{i\in I} X^{(i)}  := \left\{ (\xi_i)_{i\in I} \in \prod_{i\in I} X^{(i)}\colon \xi_i=0 \text{ for all but finitely many }i\in I\right\}.
    \]
Note that for finite $I$, the direct sum and algebraic direct sum coincide. In this case, if $X^{(i)}=X$ for all $i\in I$, we denote the direct sum by $X^{\oplus |I|}$ and note that it is isomorphic to the exterior tensor product $X\otimes \bbC^{|I|}$. For further detail see \cite[Example 1.3.5]{MR2125398}, \cite[Definition 8.1.9]{BlLM04}, and \cite[Proposition 1.2]{MR2085108} (for Hilbert C*-modules see \cite[Chapter 15]{wegge1993k}).

\subsubsection{Hilbert C*-bimodules}
We say a right $A$-$B$ C*-correspondence ${}_AX_B$ is a \textbf{Hilbert $A$-$B$-bimodule} if there is a left $A$-valued inner product ${}_A\!\langle\cdot,- \rangle$ on $X$, such that ${}_AX$ is a left $A$-$B$ C*-correspondence and such that the left $A$-norm and right $B$-norm generate equivalent topologies. We say a Hilbert $A$-$B$-bimodule is \textbf{finitely generated projective} (\textbf{fgp}) if ${}_A X$ and $X_B$ are left and right fgp, respectively. Equivalently, if it admits both finite left and right Pimsner--Popa bases.

Generalizing the notion of index for $\rm{II}_1$-factors, we can define the \textbf{right} (resp. \textbf{left}) \textbf{Watatani-index} for a Hilbert $A$-$B$-bimodule ${}_A X_B$ \cite{MR996807, MR1624182}, provided it has a finite right (resp. left) Pimsner--Popa basis $\{u_i\}_{i=1}^n\subset X$ ($\{v_j\}_{j=1}^m\subset X$): 
\begin{align}\label{eqn:WI}
    &r-\mathsf{Ind}(X):=\sum_{i=1}^n {}_A\langle u_i, u_i\rangle, && \ell\mathsf{-ind}(X):=\sum_{j=1}^m \langle v_j\mid  v_j\rangle_A.
\end{align}
In case both left and right Watatani indices exist, we say that $X$ has \textbf{finite Watatani index} given by 
$$\mathsf{Ind}_W(X):=r-\mathsf{Ind}(X)\cdot \ell-\mathsf{Ind}(X),$$ 
which can be shown to be a positive element in the center of $A$.

\subsection{Unitary tensor categories and C* 2-categories}
In this section we shall succinctly introduce the necessary definitions and relevant tools, while we also establish the notation. 
A more involved description of these concepts can be found in \cite{MR4419534, 2022arXiv220711854C}, 
and for general references on  2-categories see \cite{MR4261588}. 
\begin{defn}\label{defn:UTC}
    A \textbf{C*-category} $\cC$ is a $\bbC$-linear category satisfying the following conditions:
    \begin{enumerate}[label=(C*\arabic*)]
        \item For each pair of objects $a,b\in \cC,$ there is a conjugate-linear involution $*:\cC(a\to b)\to\cC(b\to a)$ such that for each pair of composable morphisms $f$ and $g$ we have $(f\circ g)^*=g^*\circ f^*.$
        \item There is a C*-norm on $\cC(a\to b);$ this is, the identities $||f\circ f^*||=||f^*\circ f||= ||f||^2$ hold on $\cC(a\to b)$.
        \item for each $f\in\cC(a\to b),$ there exists $g\in\cC(a\to a)$ such that $f^*\circ f = g^*\circ g.$
    \end{enumerate}
    A $*$-functor $F:\cC\to \cD$ between C*-categories satisfies $F(f^*)=F(f)^*$ for every morphism $f$ in $\cC.$
\end{defn}

\begin{defn}\label{defn:C*2cat}
    A \textbf{C* 2-category} is a 2-category $\cC$ where every $\hom$ 1-category is a C*-category, and the horizontal composition $-\boxtimes-$ of 2-morphisms is compatible with the $*$-structure.
    For objects $a,b\in\cC$ and 1-morphisms ${}_aX_b, {}_aZ_b\in\cC(a\to b)$ and ${}_bY_c\in\cC(b\to c),$ we write 1-composition $-\boxtimes-$ from left to right, e.g. ${}_aX\boxtimes_b Y_c\in\cC(a\to c)$. (We shall drop the leg notation whenever it is convenient or no confusion arises.)
    
    For 1-morphisms $X,Z,W\in\cC(a\to b)$ and 2-morphisms $f\in\cC(X\Rightarrow W)$ and $g\in\cC(W\Rightarrow Z)$ we write the vertical 2-composition $-\circ-$ from right to left. e.g.  $g\circ f\in\cC(X\Rightarrow Z)$. 
    
    We will assume that all hom 1-categories in our C* 2-categories are \emph{Cauchy complete}.
    This means that for each $a,b\in\cC$, the C* category $\cC(a\to b)$ admits all orthogonal direct sums, and also that all idempotents split via an isometry.
    For further details see \cite[Assumption 2.7]{MR4419534}.
\end{defn}

We say a 2-category $\cC$ is \textbf{rigid/has duals} if for every 1-morphism $X\in\cC(a\to b)$ there exists some  $\overline{X}\in\cC(b\to a)$ together with \emph{evaluation and coevaluation maps} $\ev_X\in\cC(\overline{X}\boxtimes_a X\Rightarrow 1_a)$ and $\coev_X\in\cC(1_b\Rightarrow X\boxtimes_b \overline{X})$ 
satisfying the \emph{Zig-Zag equations}: $\id_X=(\id_X\boxtimes \ev_X)\circ(\coev_X\boxtimes \id_X)$ and $\id_{\overline X} = (\ev_X\boxtimes \id_{\overline X})\circ(\id_{\overline X}\boxtimes\coev_X).$ 
Additionally, we assume that each 1-morphism $X\in\cC(a\to b)$ has a predual object $X_\vee\in\cC(b\to a)$ with $\overline{X_\vee}\cong X$ in $\cC(a\to b).$ 
Throughout this manuscript, whenever $\cC$ is a rigid C* 2-category, we will assume we have fixed a particular  \emph{unitary dual functor} $\overline{\ \cdot\ }$ on $\cC$ so that $(\overline{f})^*=\overline{f^*}$ for all 2-morphisms in $\cC.$ 
Furthermore, the choice of unitary dual functor can be made \emph{balanced} whenever the corresponding Hom C*-category has a simple unit. The term balanced means that the left and right categorical traces induced by the chosen dual match. 
We refer the interested reader for more details on unitary dual functors on C* 2-categories and UTCs to  \cite[Definition 2.9]{MR4419534} and \cite{MR4133163}. 

Recall that a \textbf{$\rm C^*$-tensor category} can be viewed a C* 2-category over a single object \cite[Remark 2.10]{MR4419534}. In practice, we forget about the object itself, but retain the data of its Hom 1-category and treat it as such. Now we introduce the categories we are most interested in in this manuscript:  
\begin{defn}
A \textbf{unitary tensor category} (\textbf{UTC}) $\cC$ is a semisimple rigid $\rm C^*$-tensor category with simple unit object $1_{\cC}$. 
Rigidity means that every object has a unitary dual, and simplicity of the tensor unit is the requirement that $\End_\cC(1_\cC)\cong\bbC$. 
If the isomorphism classes of simple objects in a unitary tensor category form a finite set, we say the category is a \textbf{unitary fusion category}. 
\end{defn}

\noindent Equivalently, one can realize a unitary tensor category as a $\rm C^*$ 2-category with a unique object, and such that every 1-morphism has a unitary dual.

We now describe the relevant functors between unitary tensor categories: 
A {\bf unitary tensor functor} is a triple $(F,\ F^1,\ F^2)$ consisting of a $*$-functor $F:\cC\to \cD,$ a chosen unitary isomorphism $F^1\in\cD(1_\cD\to F(1_\cC)),$ and a unitary natural isomorphism $F^2=\{F^2_{a,b}: F(a)\otimes_\cD F(b)\to F(a\otimes_\cC b)\}_{a, b\in\cC},$ sometimes called the \emph{tensorator} of $F.$

\subsection{The C* 2-category of unital C*-algebras}\label{subsec:category_of_C*-alg}
We now describe the main C* 2-category of interest, alongside with some important subcategories.

\begin{defn}
We denote by $\rCorr$ the C* 2-category whose: objects are unital C*-algebras denoted by $A, B, C\in\rCorr$; 1-morphisms ${}_AX_B, {}_AY_B\in\rCorr(A\to B)$ are right C*-correspondences; and 2-morphisms $f, g: {}_AX_B\Rightarrow {}_AY_B \in \rCorr_{A-B}(X\to Y)$ are adjointable $A$-$B$ bimodular maps. When the objects are clear from context, we omit the subscripts. The 1-composition $-\boxtimes-$ is given by the interior tensor product, and the 2-composition $-\circ-$ is the usual map composite. 
\end{defn}

Within this category we want to further specialize and focus attention on the dualizable correspondences: 
\begin{defn}
We denote the C* 2-subcategory of $\rCorr(A\to B)$ consisting of Hilbert $A$-$B$ bimodules by  $\Bim(A\to B)$, and $\Mod(A)$ the $\rm C^*$-category $\rCorr(\bbC\rightarrow A)$. We define $\fgpBim(A)$ to be the full unitary tensor subcategory of dualizable objects in $\rCorr(A\rightarrow A)$; that is, the fgp Hilbert $A$-$A$ bimodules. We denote by $\fgpMod(A)$ the C*-category of right fgp Hilbert $A$-modules.
\end{defn}

\begin{remark}\label{remark:conjugateBim}
    By \cite[Theorem 4.13]{MR2085108} we know that any dualizable C*-correspondence $X\in\rCorr(A\to A)$ can be canonically equipped with a left $A$-valued inner product, making it a finite (left and right) Watatani index bimodule. In turn, by \cite[Example 2.31]{MR2085108} dualizabilty of ${}_AX_A$ makes $X_A$ and ${}_AX$ finitely generated projective.
    
    In fact we can concretely realize the unitary dual of $K\in\fgpBim(A)$ as \textbf{the conjugate bimodule} ${_A}\overline{K}_A.$ 
    Here, $\overline K=\left\{\overline \xi\right\}_{\xi\in K}$ as a vector space, and whose $A$-$A$ bimodule structure is given by $a_1\rhd\overline{\xi}\lhd a_2:= \overline{a_2^*\rhd\xi\lhd a_1^*}$ for $\overline\xi\in \overline K$ and $a_1,a_2\in A.$ 
    The left and right $A$-valued inner products ${_A}\langle \overline{\eta},\ \overline{\xi}\rangle :=\langle \eta\mid  \xi\rangle_A$ and $\langle \overline\eta\mid  \overline\xi\rangle_A:= {_A}\langle \eta,\ \xi \rangle. $
    The corresponding evaluation and coevaluation can be expressed with respect to a finite right-$A$ Pimsner--Popa basis $\{u_i\}_1^n\subset K$ as:
    \begin{align*}
        \ev_K:\overline{K}\boxtimes_A K&\to {_A}A_A \qquad \text{and }\qquad &&\coev_K:{_A}A_A\to K\boxtimes_A\overline K\\
        \overline{\eta}\boxtimes\xi&\mapsto \langle \eta\mid  \xi\rangle _A \qquad && a\mapsto a\rhd\sum_{i=1}^n u_i\boxtimes\overline{u_i}. 
    \end{align*}
    Further discussion on this structure are found in \cite[Example 2.15]{MR4139893}. 
\end{remark}

Let A be a unital C*-algebra with trivial center ($A'\cap A=\bbC1$). 
We say that \textbf{$K\in\fgpBim(A)$ is simple} if $\End({}_AK_A) = \rCorr_{A-A}(K\to K)= \id_K\bbC,$ and call \textbf{$K$ irreducible} if its only fgp $A$-$A$ subbimodules are $K$ and $\{0\}.$
Since $\End({}_AK_A)$ is a finite dimensional C*-algebra as $K$ is fgp, it is easy to see that $K$ is simple if and only if it is irreducible. 
One can see that $\End(K_A)$ (and $\End({_A}K)$) is a corner of a matrix algebra over $A,$ and so every endomorphism is adjointable and an $A$-linear combination of the rank-one operators $|\eta\rangle\langle \xi|$ for $\eta, \xi\in K,$ ie \emph{finite rank operators}, making all right (or left) $A$-linear endomorphisms of $K$ are automatically adjointable.
The following lemma highlights additional features of $\fgpBim(A)$ that will be particularly useful to us in this paper. 

\begin{lem}\label{lem:fgp_properties}
Let $A$ be a unital C*-algebra with trivial center and let ${}_A X_A$ be a right C*-correspondence. For $K\in \fgpBim(A)$ and a bounded $A$-$A$ bimodular map $f\colon K\to X$, the following statements hold.
    \begin{enumerate}[label=(\arabic*)]
    \item $f\in \rCorr_{A-A}(K\to X)$; that is, $f$ is adjointable.
    \item\label{part:fgp_closed_comp_range} $f[K]$ is closed, orthogonally complemented, and fgp.
    \item\label{part:fgp_closed_comp_kernel} $\mathsf{Ker}(f)\subset K$ is closed, orthogonally complemented and fgp.
    \item If $K$ is irreducible and $f\neq 0$, then ${}_ A K_A \cong {}_A f[K]_A$ unitarily.
    \end{enumerate}
\end{lem}
\begin{proof}
We first remark that $K_1\oplus \cdots \oplus K_n\in \fgpBim(A)$ for any $n\in \bbN$ and $K_1,\ldots, K_n\in \fgpBim(A)$. Indeed, simply the take the unions of the left and right Pimsner--Popa bases for the $K_j$ to obtain a basis for the direct sum. This observation will be used implicitly in the proof below. Also, for $f=0$ all statements are immediate so we suppose throughout that $f\neq 0$.\\

\noindent (1): Notice that $K_A$ is a \emph{self-dual} Hilbert C*-module; \emph{i.e.} $K\cong K',$ the Banach C*-module of all anti-linear maps $K\to A$ (c.f. \cite[\S 2.5]{MR2125398}), since it is fgp and $A$ is unital. 
Adjointability then follows from \cite[Prop 2.5.2]{MR2125398}, stating that every bounded $A$-modular map out of any self-dual Hilbert C*-module into any arbitrary Hilbert C*-module is necessarily adjointable.\\

\noindent (2): First suppose $K$ is irreducible. Then $f^*f= \alpha \id_K$ for some non-zero $\alpha>0$, and so rescaling if necessary we obtain $f^*f = \id_K$. That is, $f$ is an isometry and consequently has closed range. Since it also adjointable by the previous part, $f[K]$ is complemented by\cite[Proposition 3.6]{MR1325694}. Restricting the co-domain of $f$ to $f[K]$ and applying \cite[Theorem 3.5]{MR1325694} shows ${}_ A K_A \cong {}_A f[K]_A$ unitarily. It follows that $f[K]$ right fgp where the right Pimsner--Popa basis is the image of the corresponding basis for $K$. If we then define a left $A$-valued inner product on $f[K]$ by ${}_A \<f(\xi)\mid f(\eta)\>:= {}_A \<\xi\mid \eta\>$ for $\xi,\eta\in K$, then the left Pimsner--Popa basis for $K$ is mapped to one of $f[K]$. Hence $f[K]\in \fgpBim(A)$.

Now, for arbitrary $K\in \fgpBim(A)$, we first decompose it into a direct sum $K=\bigoplus_{j=1}^n K_j$ of irreducible $K_j\in \fgpBim(A)$. Then by the above, $f[K] = \bigoplus_{j=1}^n f[K_j]$ is a closed, orthogonally complemented, and fgp.\\

\noindent (3): This follows from \cite[Theorem 15.3.8]{wegge1993k}, and the identity $\mathsf{Ker}(f) = (f^*f)[K]^\perp\subseteq K$.\\

\noindent (4): This was established in the proof of (2).
\end{proof}

In the previous proof, we saw that for irreducible $K\in \fgpBim(A)$, each $f\in \rCorr_{A-A}(K\to X)$ is either trivial or a multiple of an isometry. This allows us to give $\rCorr_{A-A}(K\to X)$ a Hilbert space structure, which will be essential to many of our arguments.

\begin{lem}\label{lem:Hilbert_space_structure}
Let $A$ be a unital C*-algebra with trivial center, let ${}_A X_A$ be a right C*-correspondence, and let $K\in \fgpBim(A)$ be irreducible. Then
    \[
        (f \mid g)_{\bbC} := f^* g\in \bbC \qquad f,g\in \rCorr_{A-A}(K\to X),
    \]
defines an inner product on  $\rCorr_{A-A}(K\to X)$ satisfying $(f\mid f)_\bbC^{1/2} = \|f\|$.
\end{lem}
\begin{proof}
For $\xi\in K$ one has
    \[
        \< f(\xi)\mid f(\xi)\>_A = \< f^*f(\xi)\mid \xi\>_A = (f\mid f)_\bbC \<\xi\mid \xi\>_A.
    \]
Taking norms in $A$ and a supremum over $\|\xi\|_A= 1$ yields $\|f\|^2 = (f\mid f)$. So the inner product is positive definite and the other properties follow easily.
\end{proof}

\begin{remark}\label{rmk:MapsAform}
For a given unital C*-algebra, and any given an arbitrary $K\in \fgpBim(A)$ and a finite right Pimsner--Popa basis $\{u_i\}_1^n\subset K$, we can explicitly describe it as ${}_AK_A \cong {}_ApA^n_A,$ where $p\in \Proj(M_n(A))$ is given by $p:= \left(\langle u_i\mid  u_j\rangle\right)_{i,j=1}^n.$ 
The right $A$-action is diagonal, and the left action on $pA^n$ is given by $a\rhd- = \left(\langle u_i\mid  a\rhd u_j\rangle\right)_{i,j=1}^n,$ which can directly be seen to commute with $p$.
This gives a unital $*$-monomorphism $-\rhd:A\to \End(pA^n_A) \cong pM_n(A)p.$ Furthermore, recall that (c.f. \cite[Prop 5.11]{FalguieresThesis})
    \begin{align*}
    \rCorr_{A-A}(K \to X) &\cong \left\{\vec{\xi}\in ((\bbC^n)^*\otimes X)p \Big|\  \vec{\xi} \lhd a = a\rhd \vec{\xi}\right\}=:X^K\subset X^{\oplus n},\\
\text{via } \quad f\mapsto \vec{\xi_f}:=(f(u_i))_i&, \quad\text{ and } \quad
\vec{\xi}\mapsto\left(f_{\vec{\xi}}:k\mapsto \vec{\xi}(\langle u_i\mid  k \rangle)_i\right).
    \end{align*} 
Here, $(\bbC^n)^*$ stands for row vectors and the right $A$-action on a \textbf{central vector} $\vec\xi \in X^K$ is given by the matrix multiplication $\vec{\xi} \lhd a = \left(\xi_i\right)_i\cdot\left(\langle u_i\mid  a\rhd u_j\rangle\right)_{i,j} = \left(\sum_k \xi_k\langle u_k\mid  a\rhd u_j\rangle\right)_j$. If $K$ is irreducible, notice that by Lemma~\ref{lem:Hilbert_space_structure}
    \[
        \langle \xi_f\mid \xi_g\rangle^{X^n} = \sum_{i=1}^n\langle f(u_i)\mid g(u_i)\rangle^{X} = (f\mid  g)_{\bbC}\cdot \sum_{i=1}^n  \langle u_i\mid  u_i\rangle^K.
    \]
Notice moreover the geometric significance of this inner product, where  $f[K]\perp g[K]$ if and only if $ \xi_f\perp\xi_g.$
(See also \cite[Thm 2.7]{MR1624182} for a similar characterization of central vectors.)
\end{remark}


\subsection{W*/C*-Algebra objects in \texorpdfstring{$\Vec({\cC})$}{Vec({\cC})}}\label{sec:AlgObjs}
For a given unitary tensor category $\cC$, the tensor category $\Vec(\cC)$ of $\cC$-graded vector spaces is a certain type of completion of $\cC$ that no longer consists of dualizable objects, but allows for infinite direct sums.
The \emph{Q-systems} arising from finite-index subfactors are particular objects in $\Vec(\cC),$ which actually lie inside of $\cC,$ and so this enlargement is necessary to capture infinite subfactors.

We shall now give a brief account about the involutive tensor category $\Vec(\cC)$ of $\cC$-graded vector spaces. (A more detailed discussion can be found at \cite[\S2,\S3]{JP17}.)
For once, we fix a complete set of representatives for the irreducible objects of $\cC,$ and we denote it by $\Irr(\cC).$ 

\begin{defn}\label{defn:VeccC}
    Given a UTC $\cC$ we define the involutive tensor category $\Vec(\cC)$ whose class of objects is $\{\bbV:\cC^{\op}\to \Vec\ |\ \bbV \text{ is a linear functor}\},$  and whose morphisms are natural transformations of functors. Here, $\cC^{op}$ denotes the \emph{opposite} category with the same objects as $\cC$ but with arrows and compositions in the hom spaces are reversed. 
    For $a\in \cC,$ $\bbV(a)$ are the \emph{fibers} of $\bbV,$ and notice they can be arbitrary complex vector spaces and not just finite dimensional. 
    Notice that we have the following identity for any $a\in\cC:$
    $$\bbV(a)\cong\bigoplus_{c\in\Irr(\cC)}\cC(a\to c)\otimes\bbV(c).$$
    The tensor product of $\bbV,\bbW\in\Vec(\cC)$ is given by
    $$(\bbV\otimes \bbW)(\cdot) = \bigoplus_{a,b\in\Irr(\cC)} \bbV(a)\otimes\cC(\cdot\to a\otimes b)\otimes \bbW(b).$$
    The involutive structure on $\Vec(\cC)$ maps $\bbV\in\Vec(\cC)$ to $\overline{\bbV}\in\Vec(\cC),$ where $\overline{\bbV}(c):=\overline{\bbV(\overline{c})},$ and for $\psi\in\cC(a\to b),$  $\overline{\bbV}(\psi):=\overline{\bbV(\overline \psi)}.$ (The bar denotes the conjugate vector space/mapping.)
    \noindent
    We will omit the associator, pivotal structure, and related compatibility structures and axioms. 
\end{defn}

We can extrapolate the known graphical calculus for $\cC$ into $\Vec(\cC)$ via the \textbf{Yoneda embedding}:

\begin{remark}\label{remark:GraphicalCalculus}
    There is a fully-faithful monoidal functor $\cC \to \Vec(\cC)$ given by $\cC\ni a\mapsto \mathbb{a}:=\cC(\cdot\to a)\in\Vec(\cC)$ called the \emph{Yoneda embedding}. 
    Elements in $\Vec(\cC)$ equivalent to those of the form $\mathbb{a},\mathbb{b},\mathbb{c},...$ are called compact objects in $\Vec(\cC),$ and these are precisely the objects in $\Vec(\cC)$ which are dualizable. 
    Notice that the tensor unit is given by $\mathbb{1} = \cC(\cdot\to 1_\cC).$
    
    The Yoneda Lemma gives $\bbV(a)\cong \Vec(\cC)(\mathbb{a}\to\bbV),$ so we can represent the vectors in $\bbV(a)$ using the morphisms in $\Vec(\cC)$ and its graphical calculus: 
    $$\bbV(a)\ni\xi\longleftrightarrow
    \tikzmath{
        \node at(0,1) {$\scriptstyle \bbV$};
        \draw (0,.3)--(0,.8);
        \roundNbox{fill=white}{(0,0)}{.3}{.1}{.1}{$\xi$}
        \node at(0,-1) {$\scriptstyle\mathbb{a}$};
        \draw (0,-.3) --(0,-.8);
    }\in\Vec(\cC)(\mathbb{a}\to \bbV).
    $$
    Furthermore, if $\psi\in\cC(b\to a),$ $\Theta\in\Vec(\cC)(\bbV\to\bbW)$, and $\xi\in\bbV(a)$, we represent $\bbV(\psi):\bbV(a)\to\bbV(b)$ and $\theta_a(\xi)$ graphically as:
    \begin{align*}
    &\tikzmath{
        \node at(-.3,1.2) {$\scriptstyle \bbV$};
        \draw (-.3,.3)--(-.3,1);
        \roundNbox{fill=white}{(0,0)}{.3}{.9}{.1}{$\bbV(\psi)(\xi)$}
        \node at(-.3,-1.2) {$\scriptstyle\mathbb{b}$};
        \draw (-.3,-.3) --(-.3,-1);
    } 
    = 
    \tikzmath{
        \node at(0,1.2) {$\scriptstyle \bbV$};
        \draw (0,.7)--(0,1);
        \roundNbox{fill=white}{(0,.45)}{.25}{.1}{.1}{$\xi$}
        \draw (0,0.15) --(0,-0.15);
        \node at(0.15,0) {$\scriptstyle\mathbb{a}$};
        \roundNbox{fill=white}{(0,-.45)}{.25}{.1}{.1}{$\psi$}
        \node at(0,-1.2) {$\scriptstyle\mathbb{b}$};
        \draw (0,-.7) --(0,-1);
    },
    \qquad\qquad\quad\quad    \text{ and }
    &&\tikzmath{
        \node at(-.3,1.2) {$\scriptstyle \bbW$};
        \draw (-.3,.3)--(-.3,1);
        \roundNbox{fill=white}{(0,0)}{.3}{.9}{.1}{$\theta_a(\xi)$}
        \node at(-.3,-1.2) {$\scriptstyle\mathbb{a}$};
        \draw (-.3,-.3) --(-.3,-1);
    }
    =
    \tikzmath{
        \node at(0,1.2) {$\scriptstyle \bbW$};
        \draw (0,.7)--(0,1);
        \roundNbox{fill=white}{(0,.45)}{.25}{.1}{.1}{$\theta$}
        \draw (0,0.15) --(0,-0.15);
        \node at(0.15,0) {$\scriptstyle\bbV$};
        \roundNbox{fill=white}{(0,-.45)}{.25}{.1}{.1}{$\xi$}
        \node at(0,-1.2) {$\scriptstyle\mathbb{a}$};
        \draw (0,-.7) --(0,-1);
    }.
    \end{align*}
    Further details about the graphical calculus for $\Vec(\cC)$ in the Yoneda embedding can be found \cite[\S 2.5]{JP17} including the diagrammatic representation for the tensor product.
\end{remark}

\begin{defn}\label{defn:AlgObjs}
    Let $\cT$ be a tensor category. An \textbf{algebra object} $(A,m,i)$ in $\cT$ consists of an object $A\in\cT$ equipped with a multiplication morphism $m:A\otimes A\to A$ and a unit morphism $i:1_\cT\to A$ satisfying the axioms of a unital associative algebra.
\end{defn}
Algebra objects in $\Vec(\cC)$ precisely correspond to \emph{Lax tensor functors} $\cC^{\op}\to\Vec$ \cite[Proposition 1]{JP17}. (For a lax tensor functor $F$, one only assumes its tensorator $F^2$ to consist of morphisms as opposed to isomorphisms.) 
Here, the multiplication map yields a tensorator and \emph{vice versa}.

\begin{defn}\label{defn:*strucConnected}
    A \textbf{$*$-structure} on an algebra object $\bbA\in\Vec(\cC)$ is a conjugate-linear (cf \cite[Definition 21]{JP17}) natural transformation $j:\bbA\to \bbA$ given by $\{j_c:\bbA(c)\to\bbA(\overline{c})\}_{c\in\cC}$ which is involutive: $j\circ j =\id_{\bbA}$, unital: $j_{1_c} = \bbA(r^{-1})$ and monoidal: $j_{a\otimes b}(m_{a,b}(f\otimes g)) = m_{\bar b, \bar a}(j_b(g), j_a(f))$ for all $a,b\in\cC.$ 
    (Here, the isomorphism $r\in\cC(1_\cC\to \overline{1_\cC})$ is the \emph{real structure} from the involutive structure on $\cC$). 
    Notice we are omitting the coherence structures on $\cC$ from our definition. 

    A \textbf{$*$-algebra} object in $\Vec(\cC)$ is an algebra object with a $*$-structure. 
    If $(\bbA,j^\bbA), (\bbB,j^\bbB)\in\Vec(\cC)$ are $*$-algebra objects, a $*$-natural transformation
    $\theta:\bbA\Rightarrow\bbB$ satisfies 
    $\theta_{\overline a}\left(j_a^\bbA(f)\right)= j^\bbB_a(\theta_a(f))$ for all $a\in \cC$ and $f\in\bbA(a).$ 
    \cite[Definition 22]{JP17}
    We say $\bbA$ is \textbf{connected} if $\bbA(1_\cC)\cong \bbC.$
\end{defn}

By an application of a unitary version of a renowned theorem by Ostrik \cite{MR1976459}, for a fixed C* tensor category $\cC,$ there is a correspondence between the $*$-algebra objects in $\cC$ and cyclic (i.e. singly generated) $\cC$-module dagger categories. \cite[Theorem 2]{JP17} 
(This correspondence holds up to \emph{Morita equivalence}, however, by remembering the pointing/generator of a cyclic $\cC$-module dagger category, one can recover the actual algebra from its Morita class.)  
Using this correspondence, C. Jones and Penneys defined a \textbf{C*/W*-algebra object $\bbA$ in $\Vec(\cC)$} as a $*$-algebra object whose \textbf{corresponding pointed $\cC$-module dagger category $\cM_\bbA$ is C*/W*}  \cite[Definition 25]{JP17}. 
We now briefly describe $\cM_\bbA$ via \cite[Construction 1]{JP17}:
\begin{align*}
&\mathsf{Obj}(\cM_\bbA) = \{\mathbb{a}\otimes \bbA\}_{a\in\cC}, \text{ and morphism spaces given by natural equivalences:}\\
&\qquad \qquad\quad\cM_\bbA(\mathbb{c}\otimes \bbA\to \mathbb{b}\otimes\bbA)\cong\Vec(\cC)(\mathbb{a}\to\mathbb{b}\otimes\bbA)\cong\bbA(\overline{b}\otimes a)\in\Vec.
\end{align*}

Of course, if we fix $\cC=\fdHilb,$ then C*/W*-algebra objects in $\Vec(\cC)\cong \Vec$ correspond to the usual notions of a C*/W*-algebra. 
We now provide several different examples of C*/W* algebra objects arising from different contexts by constructing the explicit cyclic module C*-categories:
\begin{ex}\label{ex:GroupAlg}
Let $\Gamma$ be a countable discrete group, and consider the UTC $\fdHilb(\Gamma)$ from Example \ref{ex:HilbgammOuterAction}. 
\emph{The group algebra object} $\bbC[\Gamma]$ represented by the functor $\fdHilb(\mathsf{For}(-)\otimes\bbC\to \bbC)$ is a C*-algebra object in $\Vec(\fdHilb(\Gamma))$. 
Here, $\mathsf{For}$ is the tensor functor that forgets the $\Gamma$-grading. 
The corresponding cyclic $\fdHilb(\Gamma)$-module C*-category generated by $\bbC$ is given by $\cM_{\bbC[\Gamma]},$ whose objects are of the form $\mathbb{a}\otimes\bbC[\Gamma]$ for $a\in\fdHilb(\Gamma)$, and whose morphisms are given under the natural equivalence by $\cM_{\bbC[\Gamma]}(\mathbb{a}\otimes\bbC[\Gamma]\to \mathbb{b}\otimes\bbC[\Gamma])\cong \Vec(\cC)(\mathbb{a}\to\mathbb{b}\otimes\bbC[\Gamma])\cong \fdHilb(\mathsf{For}(\overline{b}\otimes a)\otimes\bbC\to\bbC).$ 
But then $\cM_{\bbC[\Gamma]}$ is clearly C*. 
Further related examples are described in \cite[\S 6]{MR3948170}.
\end{ex}

\begin{ex}\label{ex:TKAlgObj}\cite[Example 5.35]{JP17}
Tannaka-Krein duality tells us that a compact quantum group $\bbG$ can be reconstructed from an associated UTC denoted by $\Rep(\bbG)$ together with a \emph{fiber functor} 
$\Rep(\bbG)\to \fdHilb.$ 
Each of these fiber functors equips $\Hilb$ with the structure of a $\Rep(\bbG)$-module category and under the aforementioned correspondence yields a W*-algebra object in $\Vec(\Rep(\bbG)).$
When ${\mathsf SU}_q(2)$, all fiber functors were classified in \cite{MR3420332}, in terms of \emph{fair and balanced graphs}.
\end{ex}

The following example provides us with generic means of constructing concrete C*-algebra objects from cyclic $\fgpBim(A)$-module C*-categories:
\begin{ex}\label{ex:GenericCStarAO}
Given a unital inclusion of C*-algebras $A\subset B,$ where $A'\cap A=\bbC1$, and focusing on the UTC $\fgpBim(A)$ we can explicitly produce families of C*-algebra objects in $\Vec(\fgpBim(A))$ as follows: 
Consider the $\fgpBim(A)$-module C*-category $\rCorr(A\to B),$ where for any $K, L\in\fgpBim(A),$ any ${}_AX_B, {}_AY_B,\in\rCorr(A\to B),$ maps  $\psi \in\fgpBim(K\to L)$ and $f\in\rCorr_{A-B}(X\to Y),$ the action is given by 
\begin{align*}
    &K\rhd {}_AX_B:= K\boxtimes_A X_B,\qquad\qquad  \text{ and } && \psi\rhd f:= \psi\boxtimes_A f.
\end{align*}
By picking the basepoint ${}_AB_B\in \rCorr(A\to B)$ there is a corresponding C*-algebra object $\bbB\in\Vec(\cC).$ 
Explicitly,  $\bbB\in \Vec(\cC)$ is given by
\begin{align*}
    \bbB: \fgpBim(A)^{\op}&\to \Vec\\
    K&\mapsto \rCorr_{A-B}(K\boxtimes_A B\to B),\\
    \left(b\xrightarrow[]{\psi}a\right) &\mapsto \bbB(\psi):\bbB(a)\to \bbB(b)\\
    &\hspace{2.4 cm} f\mapsto f\circ (\psi\boxtimes \id_B).
\end{align*}
And the algebra structure for $\bbB$ is: 
\begin{align*}
    \bbB^2_{K,L}:\bbB(K)\otimes_\bbC \bbB(L) &\to \bbB(K\boxtimes_A L)\\
    f\boxtimes g &\mapsto f\circ(\id_{K}\otimes g),
\end{align*}
with unit
\begin{align*}
    \bbB^1_K: \bbOne(K)&\to\bbB(K)\\
    \fgpBim(A)(K\to A)\owns f&\mapsto m_B\circ[f \boxtimes \id_B]\in\rCorr_{A-B}(K\boxtimes_A B\to B).
\end{align*}
Here, $m_B:B\otimes B\to B$ is the multiplication map.
Furthermore, since
$$\bbB({}_AA_A) =\rCorr_{A-B}(A\boxtimes_A B\to B) \cong \rCorr_{A-B}(B\to B)\cong A'\cap B,$$ 
if we demand that $\dim_\bbC(A'\cap B)<\infty$ then $\bbB({}_AA_A)$ is finite dimensional. Therefore, by \cite[Prop 2.8]{MR3948170} it follows that $\dim_\bbC(\bbB(K))\leq\mathsf{Ind}_W(K)<\infty$ for every $K\in\fgpBim(A).$
We postpone describing the $*$-structure on $\bbB$ pulled from its corresponding C* module-category $\cM_\bbB$  until Lemma \ref{lem:StarStruc}. 
\end{ex}

\section{Inclusions of C*-algebras}\label{sec:Inclusions}

Let $A\subset B$ be a unital inclusion of C*-algebras equipped with a faithful conditional expectation  $E: B\to A.$
Let $\cB$ be the completion of $B$ as a right Hilbert C* $A$-module with inner product given by $\langle b_1\mid  b_2\rangle_A := E(b_1^* b_2).$ 
We have an embedding of $B$ inside of $\cB,$ 
whose image will be denoted by $B\Omega\subset\cB.$ 
The norm induced by the $A$-valued inner product on $B\Omega$ is denoted on elements by  $\|b\Omega\|_A:=\|\langle b\Omega\mid b\Omega\rangle_A\|^{1/2}.$ 

Given an inclusion $A\overset{E}{\subset} B$, a sufficient condition so that $\cB\cong B\Omega$ as right C* $A$-modules, ie $\|\cdot\|_B\cong\|\cdot\|_A$ on $B$, so it is $\|\cdot\|_A$-complete, is that there is some constant $c\geq 1$ such that
$\|\cdot\|_B\leq c\|\cdot\|_A,$ since $\|\cdot\|_A\leq \|\cdot\|_B$ is automatic by the positivity of $E$. 
We say that $E$ has finite \textbf{Pimsner--Popa index}, if such $c$ giving the norm comparison exists. 
In this case, we define the value of the Pimsner--Popa index as 
$\textsf{Ind}_p(E)= \textsf{inf}\{c\geq 1|\ (cE-\id)\geq 0\}.$
If no such $c\geq 1$ exists, we say the Pimsner--Popa index for the inclusion is infinite.  \cite[Definition 2]{MR1642530}, \cite{MR860811}. 
We record these conclusions in the following proposition:
\begin{prop}\label{prop:finite_PP_index}
	$E$ as above has finite Pimsner--Popa index if and only if $\cB = B\Omega;$ 
	i.e. $B$ is already complete in $A$-norm. (\cite[Proof of Prop 2.1.5]{MR996807} and \cite[Theorem 1]{MR1642530})
\end{prop}

In \cite[Proposition 2.1.5]{MR996807}, Watatani establishes that finiteness of the Pimsner--Popa index together with unitality of the right $A$-compact operators $\cK(B_A)$ (i.e. the \textbf{reduced C*-basic construction} for $A\subset B$) are equivalent to the finiteness of the Watatani index.
Finiteness of the Pimsner--Popa index for an inclusion $A\overset{E}{\subset}B$ alone does not imply the existence of a finite Pimsner--Popa basis for $E$.  
Indeed, Watatani and Frank-Kirchberg \cite{MR1642530} exhibited examples of this behaviour using unital commutative C*-algebras.
However, as established by Izumi in \cite[Corollary 3.4]{MR1900138}, this disparate behavior disappears, once one restricts attention to unital inclusions $A\overset{E}{\subset} B,$ where $A$ is simple. 
Namely, in this case one gets $\textsf{Ind}_p(E)< \infty$ if and only if $\textsf{Ind}_W(E)< \infty$ if and only if ${}_A\cB_A\in\fgpBim(A),$ ie $B\Omega\cong \cB$ is finitely generated projective as a right and left $A$ Hilbert C*-module.

Recall from Remark~\ref{remark:conjugateBim} that the conjugate bimodule $\overline{K}$ for $K\in \fgpBim(A)$ was defined as a formal conjugate. When $K\subset B\Omega$, this conjugate is implemented by the adjoint on $B$:

\begin{lem}\label{lem:ConcreteConjugateBim}
    Let $A$ be a unital C*-algebra with trivial center. 
    For $\fgpBim(A)\owns K\subset B\Omega,$ the $*$-structure on $B$ gives an $A$-$A$ isomorphism 
    \begin{align*}
        K^* \to\overline{K} \quad \text{  given by  }\quad b^*\Omega \mapsto \overline{b\Omega}. 
    \end{align*}
\end{lem}
\begin{proof}
    By direct computation.
\end{proof}

\subsection{Projective-quasi-normalizers and the Diamond Spaces}
Popa introduced the notion of the \emph{quasi-normalizer} for a unital inclusion of von Neumann algebras in \cite{MR1729488} and used it to compare the properties of their symmetric enveloping inclusions with those given by crossed products of factors by discrete groups and to intrinsically characterize the discrete inclusions studied in \cite{MR1622812}. 
By definition, a subfactor $N\subset M$ is \emph{quasi-regular} if its quasi-normalizer ${\sf QN}(N\subset M)$ is a dense unital $*$-subalgebra of $M.$
C. Jones and Penneys used this machinery together with a generalized crossed product construction called \emph{realization} to produce new examples of infinite-index subfactors and described a correspondence between intermediate subfactors and W*-algebra objects in a unitary tensor category.\cite{MR3948170}

For a unital inclusion of C*-algebras $A\overset{E}{\subset} B$, the \emph{quasi-normalizer} is
\begin{align*}\label{defn:QN}
    {\sf QN}\left(A\subset B\right) 
    := \set{b\in B}{\ \exists\ \{x_i\}_{i=1}^n, \{y_j\}_{j=1}^m\subset B,\ b A \subseteq \sum_{i=1}^n A x_i,\ A\cdot b \subseteq \sum_{j=1}^my_j \cdot A}.
\end{align*}
Clearly, $A\subset {\sf QN}(A\subset B) \subset B$ is an intermediate $*$-algebra.
However, for the purpose of studying infinite-index inclusions of C*-algebras, it is not clear that this definition is strong enough.
So now we introduce a modified quasi-normalizer, which is visible from the viewpoint of unitary tensor categories:

\begin{defn}\label{defn:PQR}
For a unital inclusion of C*-algebras $A\overset{E}{\subset} B$, we define the \textbf{projective quasi-normalizer} as 
	\begin{align*}
	    {B^{\diamondsuit}}:=\PQN\left(A\overset{E}{\subset} B\right)&:=\set{b\in B}{\ \exists\  \fgpBim(A)\owns K\subset B\Omega,\ \overline{A\rhd b\Omega \lhd A}^{||\cdot||_A}\subset K}\\
	    &\ =\set{b\in B}{\ \exists\ \fgpBim(A)\owns K\subset B\Omega,\ b\Omega\in K}.
	\end{align*}
We say the inclusion is \textbf{projective-quasi-regular} (\textbf{PQR}) if $\overline{B^{\diamondsuit}\Omega}^{\|\cdot\|_A}=\cB.$
We will customarily drop the conditional expectation from the notation  unless it is not clearly specified contextually and write simply $\PQN(A\subset B)$ or $B^{\diamondsuit}$. 
\end{defn}

It is clear that $A\subset B^{\diamondsuit}$, and that $B^{\diamondsuit}$ is closed under taking $*$ and under sums. 
Indeed, for $i=1,2,$ if $b_i\Omega\in K_i$ are as above, then clearly $(b_1+b_2)\Omega\in K_1+K_2  \subset B\Omega,$ and $b_i^*\Omega\in K_i^*\subset B\Omega.$ The latter is clearly fgp with the conjugate structure, and $K_1+K_2$ is fgp by Lemma~\ref{lem:fgp_properties}.\ref{part:fgp_closed_comp_range} as the image of $K_1\oplus K_2\in \fgpBim(A)$ under the map $(\xi,\eta)\mapsto \xi+\eta$. However, proving that $B^{\diamondsuit}$ actually a $*$-subalgebra will require more sophisticated categorical technology and is established in Lemma~\ref{lem:deltaproperties}. The first step in this direction is Proposition~\ref{prop:PQNAlgRealization}, where we show $B^{\diamondsuit}$ is a vector space graded over a fusion algebra. We first require some additional notation, which  allow us to work a bit more abstractly and leverage the categorical structure.

\begin{defn}\label{defn:B-valuedmaps}
    Let $A$ be a unital C*-algebra with trivial center. 
    For each $K\in \fgpBim(A),$ we define the \textbf{diamond space}  
    $$\rCorr_{A-A}(K\to \cB)^\diamondsuit := \set{f\in \rCorr_{A-A}(K\to \cB)}{\  f[K]\subset B\Omega}.$$ 
For $f\in \rCorr_{A-A}(K\to \cB)^\diamondsuit$,  define $\check{f}\colon K\to B$ by letting $\check{f}(\xi)\in B$ for $\xi\in K$ be the unique element satisfying $\check{f}(\xi)\Omega = f(\xi)$.
\end{defn}

The key property of $f\in\rCorr_{A-A}(K\to \cB)^\diamondsuit$ is that $\check{f}[K]\subset B^{\diamondsuit}=\PQN(A\subset B)$: for any $\xi\in K$ one has $\check{f}(\xi)\Omega= f(\xi)\in K$. The diamond space in general need not be a closed subspace of $\rCorr_{A-A}(K\to \cB)$. However, when $A\subset B$ is an irreducible inclusion we will see that the diamond spaces are always finite dimensional and hence closed (see Theorem~\ref{thm:fgpInsideB}). This relies on C*-Frobenius reciprocity (see Theorem~\ref{thm:FR}), but can be seen directly for $K={}_A A_A\in \fgpBim(A)$: each $f\in \rCorr_{A-A}(A\to \cB)^\diamondsuit$ is determined by $\check{f}(1)\in A'\cap B=\bbC$.

The following proposition is analogous to \cite[Lemma 2.5.(5)]{MR3801484}, 
\cite[Theorem 2.8.10]{MR2125398}, and \cite[Remark 3.18]{MR3948170}. Recall from Lemma~\ref{lem:Hilbert_space_structure} that when $A$ has trivial center then $\rCorr_{A-A}(K\to \cB)$ is a Hilbert space for irreducible $K\in \fgpBim(A)$. Consequently, the exterior tensor product $K\otimes \rCorr_{A-A}(K\to \cB)$ can be regarded as a right $A$-$A$ C*-correspondence (see Remark~\ref{rem:ext_tens_with_Hilb}).

\begin{prop}\label{prop:PQNAlgRealization}
Let $A\overset{E}{\subset} B$ be a unital inclusion of C*-algebras where $A$ has trivial center, let $\cB$ be the $B$-$A$ C*-correspondence obtained from $E$, and $\cK$ be a complete set of representatives of irreducible finitely generated projective $A$-$A$ subbimodules of $\cB$. Then the map
    \[
        \bigoplus_{K\in\cK} K \odot \rCorr_{A-A}(K\to \cB) \ni \sum_{j=1}^d \xi_j \otimes f_j \mapsto \sum_{j=1}^d f_j(\xi_j)\in \cB,
    \]
extends to an isometry
    \[
        \delta_{\cB} \in \rCorr_{A-A}\left( \left[\overline{\bigoplus_{K\in \cK}} K\otimes \rCorr_{A-A}(K\to \cB)\right] \to \cB \right). 
    \]
Moreover, the image of the subspace
    \[
        \bigoplus_{K\in\cK} K\odot \rCorr_{A-A}(K\to \cB)^\diamondsuit
    \]
under this map is ${\sf PQN}(A\subset B)\Omega$, and therefore one obtains a map
    \begin{align*}
        \delta_B^{\diamondsuit} \colon  \bigoplus_{K\in\cK} K\odot \rCorr_{A-A}(K\to \cB)^\diamondsuit &\to \PQN(A\subset B)\\
        \sum_{j=1}^d \xi_j \otimes f_j &\mapsto \sum_{j=1}^d \check{f_j}(\xi_j).
    \end{align*}
\end{prop}
\begin{proof}
First note that the irreducibility of elements in $\cK$ implies that for distinct $K,K'\in\cK$ one has
    \[
        \< \delta_{\cB}(x), \delta_{\cB}(y)\>_A=0
    \]
for any $x\in K\odot \rCorr_{A-A}(K\to \cB)$ and $y\in K'\odot \rCorr_{A-A}(K'\to \cB)$. Thus is suffices to show $\delta_{\mathcal{B}}$ is an isometry on each summand $K\odot \rCorr_{A-A}(K\to \cB)$. Using Lemma~\ref{lem:Hilbert_space_structure} we have
    \begin{align*}
        \left\< \delta_{\cB}\left( \sum_{j=1}^d \xi_j\otimes f_j \right)\  \Bigg|\ \delta_{\cB}\left( \sum_{k=1}^e \eta_k\otimes g_k \right) \right\>_A &= \sum_{j=1}^d \sum_{k=1}^e \< f_j(\xi_j) \mid g_k(\eta_k) \>_A \\
        &= \sum_{j=1}^d \sum_{k=1}^e \<\xi_j \mid \eta_k\>_A (f_j\mid g_k)_{\bbC}\\
        &= \sum_{j=1}^d \sum_{k=1}^e \<\xi_j\otimes f_j \mid \eta_k\otimes g_k\>_{A}\\
        &= \left\<  \sum_{j=1}^d \xi_j\otimes f_j\  \Bigg|\ \sum_{k=1}^e \eta_k\otimes f_k \right\>_A.
    \end{align*}
Also, for $a,b\in A$ one has
    \[
        a \rhd \delta_{\cB}\left( \sum_{j=1}^d \xi_j\otimes f_j \right) \lhd b = \sum_{j=1}^d f_j( a\rhd \xi_j \lhd b) = \delta_{\cB}\left( a\rhd \left(\sum_{j=1}^d \xi_j\otimes f_j\right) \lhd b\right).
    \]
Thus $\delta_{\cB}$ is $A$-$A$ bimodular.

Next, note for $\xi\in K$ and $f\in \rCorr_{A-A}(K\to \cB)^\diamondsuit$ one has $\delta_{\cB}(\xi\otimes f) = f(\xi) = \check{f}(\xi)\Omega\in \PQN(A\subset B)\Omega$ (see the discussion following Definition~\ref{defn:B-valuedmaps}). Conversely, given $b\in \PQN(A\subset B)$ let $\fgpBim(A)\ni K\subset B\Omega$  be such that $b\Omega\in K$. Note the inclusion map $\iota\colon K\to \cB$ satisfies $\iota\in \rCorr_{A-A}(K\to \cB)^\diamondsuit$, and so if
    \[
        K= f_1[K_1] + \cdots + f_d[K_d]
    \]
for $K_1,\ldots, K_d\in \cK$ then $f_1,\ldots, f_d$ lie in their corresponding diamond spaces. Hence $b\Omega\in \delta_{\cB}\left( \bigoplus_{j=1}^d K_j\odot \rCorr_{A-A}(K_j\to \cB)^\diamondsuit\right)$.
\end{proof}

Later on in Section \ref{subsec:CrossProd}, we will show how to endow $\PQN(A\subset B)$ with a natural $*$-algebra structure. Ultimately, the map $\delta^{\diamondsuit}_B$ as defined above will become a component of a natural isomorphism that will allow us to precisely characterize those inclusions of C*-algebras which are \emph{discrete}. (Compare with Definition \ref{defn:delta}.) It will also be essential in the proof of Lemma~\ref{lem:deltaproperties}, where we finally show $\PQN(A\subset B)$ is closed under multiplication.

At this moment, it is worth mentioning that there are vast families of examples of PQR inclusions. For now, we will delay providing with such examples to Section \ref{subsec:CrossProd} after we have introduced a (potentially stronger) version of projective-quasi-regularity for $A\subset B$ named \emph{C*-discrete} in Definition \ref{defn:CDisc}, which we will be able to characterize abstractly in Section \ref{sec:Characterization}. 
(See discussion after Theorem \ref{thm:CrossProdPQR}.) 
Instead, in the next section we shall explore the structure of PQR inclusions abstractly to establish their main features and properties.

\subsection{C*-Frobenius Reciprocity}
In this section we establish strong and useful properties for PQR inclusions of C*-algebras, and will begin to make clear how some of these features arise from an \emph{underlying quantum symmetry} governed by an action of a UTC. 
Our statements and techniques are inspired from discrete subfactors, and should be contrasted with  \cite[Lemma 2.5]{MR3801484}, and \cite[Prop. 3.11]{MR3948170}.

We now record the C*-algebraic analog of Frobenius reciprocity, inspired by \cite[Theorem 3.13]{MR3948170}, which provides us with a direct link to C*-algebra objects in unitary tensor categories \cite{JP17} as will be evident later. 
Remarkably, in the context of discrete subfactors there is an analogous statement of Frobenius Reciprocity, which relies heavily on the existence of traces or choices of normal faithful states, together with their associated Modular Tomita-Takesaki theory used to obtain a bounded right $B$-action on the GNS space $L^2(B).$ 
So a potent feature of our result is that it not only holds at the Hilbert C*-module/C*-correspondence level, but also that it makes no reference to states at all. 
Furthermore, the result holds without additional assumption on simplicity of the C*-algebras or trivial-center hypotheses.

\begin{thm}[{C*-Frobenius Reciprocity}]\label{thm:FR}
    Let $A\subset B$ be a unital inclusion of C*-algebras admitting a faithful conditional expectation  $E: B\twoheadrightarrow A$ 
    and let $K\in \fgpBim(A).$ 
    Then, the maps
    \begin{align*}
        \Psi_K:\rCorr_{A-A}(K\to \cB)^\diamondsuit &\longrightarrow \rCorr_{A-B}\left(K\underset{A}{\boxtimes}B\to B\right) \\
        \Psi_K(g)(k\boxtimes b)&:= \check{g}(k)b,\\
        \text { and }\quad \Phi_K: \rCorr_{A-B}\left(K\underset{A}{\boxtimes}B\to B\right) &\longrightarrow \rCorr_{A-A}(K\to \cB)^\diamondsuit\\
        \Phi_K(f)(k)&:= f(k\boxtimes 1_A)\Omega
    \end{align*}
    assemble into mutually inverse natural isomorphisms, i.e. morphisms in $\Vec(\fgpBim(A))$:
    \begin{align*}
        \Psi:\rCorr_{A-A}(-\to \cB)^\diamondsuit \Rightarrow \rCorr_{A-B}(-\boxtimes_A B\to B), \\    
        \Phi:\rCorr_{A-B}(-\boxtimes_A B\to B)\Rightarrow \rCorr_{A-A}(-\to \cB)^\diamondsuit.   
    \end{align*}
    \noindent Here, $\check g (k)\Omega=g(k)$ is uniquely determined as an element in $B$ as $E$ is faithful, and we are viewing $B$ as a right $A$-$B$ C*-correspondence.
    Diagrammatically, these maps become: 
    \begin{align*}
    \Psi_K:\ 
    \tikzmath
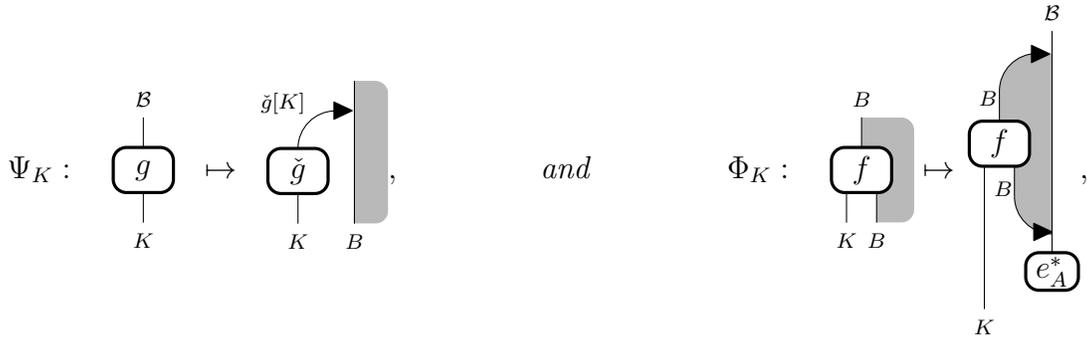
{
        \begin{scope}
            \clip[rounded corners=5pt] (-.7,-.7) rectangle (.7,.7);
        \end{scope}
        \draw (0,0) -- (0,.7) node[above]{$\scriptstyle \cB$};
        \draw (0 ,-.7) node[below]{$\scriptstyle K$} -- (0,0);
        \roundNbox{fill=white}{(0,0)}{.3}{.1}{.1}{$g$}
    }
    \mapsto
    \tikzmath{
        \begin{scope}
            \clip[rounded corners=5pt] (-.7,-.7) rectangle (.7,.7);
        \end{scope}
        \draw (0,0) -- (0, 0.3) arc (180:90:0.5);
        \draw (-0.2, 0.9) node{$\scriptstyle \check g[K]$};
        \draw (0.6, 0.8) node{$\blacktriangleright$};
        \draw (0, -0.7) node[below]{$\scriptstyle K$} -- (0,0);
        \begin{scope}
        \clip[rounded corners=5pt] (.5,-.7) rectangle (1.2, 1.2);
            \fill[\BColor] (0.75, -.7) -- (0.75, 1.2) -- (1.2, 1.2) -- (1.2, -0.7) -- (0.65,-.7);
        \end{scope}
        \draw (0.75 ,-.7) node[below]{$\scriptstyle B$} -- (0.75, 1.2);
        \roundNbox{fill=white}{(0,0)}{.3}{.1}{.1}{$\check g$}
    },
    \qquad\qquad \text{ and }\qquad \qquad
    \Phi_K:\ 
    \tikzmath{
        \begin{scope}
        \clip[rounded corners=5pt] (-.7,-.7) rectangle (.7,.7);
        \fill[\BColor] (.7,-.7) -- (.2,-.7) -- (.2,0) -- (0,0) -- (0,.7) -- (.7,.7);
        \end{scope}
        \draw (0,0) -- (0,.7) node[above]{$\scriptstyle B$};
        \draw (-.2,-.7) node[below]{$\scriptstyle K$} -- (-.2,0);
        \draw (.2,-.7) node[below]{$\scriptstyle B$} -- (.2,0);
        \roundNbox{fill=white}{(0,0)}{.3}{.1}{.1}{$f$}
    }
    &\mapsto
    \tikzmath{
        \fill[\BColor] (0,0) -- (0,.5) arc (180:90:0.5) -- (.7,1) -- (.7,0) -- (0,0);
        \fill[\BColor](.2,0) -- (.2,-.5) arc(180:270:.5) -- (.7,.5);
        \draw (0.7, -1.5) -- (0.7, 1.3) node[above]{$\scriptstyle\cB$};
        \draw (0,0) -- (0,.5) arc (180:90:0.5);
        \node at(-.1,.75){$\scriptstyle B$};
        \draw (0.55, 1) node{$\blacktriangleright$};
        \draw (0.58, -.98) node{$\blacktriangleright$};
        \draw (-0.2,-1.9) node[below]{$\scriptstyle K$} -- (-.2,0);
        \draw (0.2,-0.5) -- (.2,0);
        \node at(0.05,-.6){$\scriptstyle B$};
        \draw (0.2,-0.5) arc (180:270:0.5);
        \roundNbox{fill=white}{(0,0)}{.3}{.1}{.1}{$f$}
        \roundNbox{fill=white}{(0.7, -1.45)}{.25}{.1}{.1}{$e_A^*$}
    },
    \end{align*}
where the shadings denote whether we consider $A$-$A$ or $A$-$B$ bimodular maps.
\end{thm}
\begin{proof}
    We fist show that $\Psi_K$ is well-defined.
    Notice that $K\boxtimes_A B\in\rCorr(A\to B)$ is fgp as a right Hilbert $B$-module. 
    Since clearly $\Psi_K (g)$ is $A$-$B$ bimodular, and letting  $\{u_i\}_1^n\subset K$ be a Pimsner--Popa basis, then $\Psi_K (g)$ is $B$-compact, and even finite rank, by expressing 
    $\Psi_K (g)(k\boxtimes b)=\sum_1^n |\check{g}(u_i)\rangle\langle u_i\boxtimes 1|\ k\boxtimes b\rangle_B.$ 
    (This is also described in \cite[Proof of Theorem 8.1.27 (2)]{BlLM04}.) 
    Thus $\Psi_K(g)$ is an adjointable $A$-$B$ bimodular map. 
    
    That $\Psi_K$ and $\Phi_K$ are mutual inverses follows by direct computation, and the naturality of $\Psi$ and $\Phi$ is also straightforward.  
\end{proof}

Frobenius reciprocity lies at the heart of our C*-reconstruction technique (Section \ref{subsec:CAlgObjFun}) as it provides a direct link between certain C*-inclusions and abstract C*-algebra objects in $\cC$-graded vector spaces.
In the next proposition, we prove what might be the strongest feature of PQR inclusions. 
As an application of Frobenius reciprocity, in the presence of an irreducible inclusion, PQR implies the number of copies of a given fgp $A$-$A$ bimodule that $B\Omega\subset \cB$ can host is finite. 

\begin{thm}\label{thm:fgpInsideB}
    Let $A\overset{E}{\subset}B$ be a unital inclusion of C*-algebras. Then the diamond space linear functor 
    $$
        \fgpBim(A)^{\op}\to \Vec\ \ \text{ given by }\ \  L\mapsto \rCorr_{A-A}(L\to \cB)^\diamondsuit
    $$
    is in fact a C*-algebra object in $\Vec(\fgpBim(A)).$ 
    If $A\subset B$ is irreducible, then this object is connected and for each irreducible $K\in\fgpBim(A)$ we have 
    \[
            \dim_{\bbC} \left(\rCorr_{A-A}(K \to \cB)^\diamondsuit\right) \leq \mathsf{Ind}_W(K) < \infty.
    \]
    More generally, for pairwise non-isomorphic irreducible $K_1,\ldots, K_m\in \fgpBim(A)$ and integers $n_1,\ldots, n_m\in \bbN$ one has
        \[
            \dim_{\bbC} \left(\rCorr_{A-A}\left( \bigoplus_{j=1}^m K_j^{\oplus n_j} \to \cB\right)^\diamondsuit\right) \leq \sum_{j=1}^m n_j^2\mathsf{Ind}_W(K_j) <\infty.
        \]
    If one further assumes $A\subset B$ is a {\PQR}-inclusion, then for each $L\in\fgpBim(A)$ $$\rCorr_{A-A}(L \to \cB)^\diamondsuit = \rCorr_{A-A}(L \to \cB).$$
\end{thm}
\begin{proof}
By the C*-Frobenius Reciprocity Theorem \ref{thm:FR}, we see that $\rCorr_{A-A}(-\to \cB)^\diamondsuit \cong \rCorr_{A-B}(-\boxtimes_AB\to B)$ is manifestly a C*-algebra object in $\Vec(\fgpBim(A))$. (See Example \ref{ex:GenericCStarAO}.) Indeed, by \cite[\S 3]{JP17}, for any unitary tensor category $\cC$, there is an equivalence of categories between  $*$-algebra objects in $\Vec(\cC)$ and pointed cyclic $\cC$-module dagger categories $(\cM, m)$. Under this equivalence, with $\cC = \fgpBim(A),$ $\rCorr_{A-B}(-\underset{A}{\boxtimes}B\to B)$ corresponds to the connected C*-algebra object in $\Vec(\fgpBim(A))$ generated by basepoint $m = {}_AB_B\in \rCorr(A\to B) =: \cM,$ the ambient $\fgpBim(A)$-module C*-category of right  $A$-$B$ correspondences. 

If $A\subset B$ is irreducible, then $A'\cap B=\bbC$  translates to $\dim_\bbC(\rCorr_{A-A}(A\to\cB)^\diamondsuit) = 1$. Then, by \cite[Proposition 2.8]{MR3948170} the dimension of each fiber  $\dim_\bbC(\rCorr_{A-A}(K\to\cB)^\diamondsuit)$ for $K\in\fgpBim(A)$ is finite and bounded above by the Watatani index $\mathsf{Ind}_W(K)$. 

Now further suppose $A\subset B$ is a {\PQR}-inclusion. Observe that since $A'\cap A\subset A'\cap B= \bbC$, then $\rCorr_{A-A}(K \to \cB)$ is a Hilbert space with inner product $(g\mid  f):= g^*\circ f $ by Lemma~\ref{lem:Hilbert_space_structure}. 
 Since $\rCorr_{A-A}(K \to \cB)^\diamondsuit$ finite dimensional by the above, it is automatically closed and therefore complemented in this Hilbert space. Fix $g\in \rCorr_{A-A}(K\to \cB)$ and let 
    \[
        g = g' \oplus g'' \in \rCorr_{A-A}(K \to \cB)^\diamondsuit \oplus 
        \left(\rCorr_{A-A}(K \to \cB)^\diamondsuit\right)^\perp.
    \]
We must show $g''=0$ and this will be accomplished by establishing that $g''[K]$ is orthogonal to the dense subspace $\PQN(A\subset B)\Omega$. Let $b\in \PQN(A\subset B)$ with $b\Omega\in L\subset B\Omega$ for $L\in \fgpBim(A)$. Decompose $L \cong \bigoplus_{j=1}^m L_j$ into irreducible fgp $A$-$A$ bimodules, let $\iota_j \in \rCorr_{A-A}(L_j\to \cB)^\diamondsuit$ be the embedding $\iota_j(L_j)\subset L$ for each $j=1,\ldots,m$, and let $b\Omega = \sum_j \iota_j (\ell_j)$ for $\ell_j\in L_j$. Then for $k\in K$ we have
    \[
        \<g''(k)\mid \iota_j(\ell_j)\>_A =0
    \]
for each $j=1,\ldots, m$, since either: $L_j\cong K$ so that $\iota_j\in \rCorr_{A-A}(K \to \cB)^\diamondsuit$; or $L_j\not\cong K$ so that $\iota
_j^*\circ g''=0$ by Schur's Lemma. Consequently,
    \[
        \<g''(k) \mid b\Omega\>_A = \sum_{j=1}^m \<g''(k)\mid \iota_j(\ell_j)\>_A =0,
    \]
and hence $g'[K] \perp \PQN(A\subset B)\Omega$. This completes the case for irreducible $K\in \fgpBim(A)$. 

For general $L\in \fgpBim(A)$  let $L\cong\bigoplus_{j=1}^m K_j^{\oplus n_j}$ for pairwise non-isomorphic and irredicuble $K_1,\ldots, K_m\in \fgpBim(A)$ and integers $n_1,\ldots, n_m\in\bbN$. Then one has
        \[
            \rCorr_{A-A}(L\to \cB) \cong \bigoplus_{j=1}^m M_{n_j}\left( \rCorr_{A-A}(K_j\to \cB)\right), 
        \]
and similarly for the diamond space. Thus the result in this case follows by appealing to the irreducible case.
\end{proof}

\begin{remark}\label{cor:FinMultiplicities}
A close inspection of the proof of Theorem~\ref{thm:fgpInsideB} will reveal that the irreducibility assumption can be weakened in a few ways. For example, to obtain that the diamond spaces are finite dimensional one merely needs $A'\cap A=\bbC$ and a local finiteness assumption $\dim_{\bbC}(A'\cap B)<\infty$. To obtain the equality $\rCorr_{A-A}(L \to \cB)^\diamondsuit = \rCorr_{A-A}(L \to \cB)$, one needs only that the former is closed and $A'\cap A=\bbC$ We will however mostly adhere to the irreducible cases.  
\end{remark}

A remarkable characterization of PQR inclusions is that the $A$-$A$ correspondence $\cB$ decomposes as an orthogonal direct sum of fgp $A$-$A$ bimodules with finite multiplicities. We make this precise in Corollary \ref{cor:P-W}, which could be regarded as a \emph{generalized Peter-Weyl Theorem}. We view this as an initial  manifestation of discreteness in the C*-setting. (cf \cite[Proposition 7.13]{FalguieresThesis} and \cite[Corollary 3.24]{MR3948170}.)
On the classical side, Mishchenko described a family of inclusions of C*-algebras from representations of compact groups on Hilbert C*-Modules \cite{MR752175} obtaining orthogonality relations and a Peter-Weyl decomposition similar to ours (see also \cite[Thm 2.8.10]{MR2125398}).

\begin{cor}\label{cor:P-W}
Let $A\overset{E}{\subset} B$ be a unital irreducible inclusion of C*-algebras, let $\cB$ be the $B$-$A$ C*-correspondence obtained from $E$. Then, $A\subset B$ is PQR if and only if $\rCorr_{A-A}(K\to \cB)=\rCorr_{A-A}(K\to \cB)^\diamondsuit$ for all $K\in \Irr(\fgpBim(A))$ and $\cB$ admits an orthogonal  ``Peter-Weyl'' decomposition as an $A$-$A$ bimodule: 
        \[
            {}_A\cB_A \cong \overline{\bigoplus_{K\in \Irr(\fgpBim(A))}} K^{\oplus n_K},
        \]
where $n_K = \mathsf{dim}_{\bbC}\rCorr_{A-A}(K\to \cB)\in \bbN$ and the isomorphism is implemented by $\delta_\cB$.
 \end{cor}
\begin{proof}
$(\Rightarrow)$: Theorem~\ref{thm:fgpInsideB} implies $\rCorr_{A-A}(K\to \cB)=\rCorr_{A-A}(K\to \cB)^\diamondsuit$. Then assumed density of $\PQN(A\subset B)\Omega$ gives
    \[
        {}_A\cB_A \cong \overline{\bigoplus_{K\in \cK}} K\otimes \rCorr_{A-A}(K\to \cB) 
    \]
via the isometry from Proposition~\ref{prop:PQNAlgRealization}. Hence $K$ has multiplicity $n_K:=\mathsf{dim}_{\bbC}\rCorr_{A-A}(K\to \cB)$, and this is necessarily finite by Theorem~\ref{thm:fgpInsideB}. Thus $K\otimes \rCorr_{A-A}(K\to \cB) \cong K^{\oplus n_K}$ for each $K\in \cK$ yields the claimed decomposition for $cB$.

\noindent $(\Leftarrow)$: By Proposition~\ref{prop:PQNAlgRealization} we have
    \begin{align*}
        \overline{\PQN(A\subset B)\Omega} & = \delta_{\cB}\left[\overline{ \bigoplus_{K\in \cK} K\odot \rCorr_{A-A}(K\to \cB)^\diamondsuit }\right]\\
        &= \delta_{\cB}\left[\overline{\bigoplus_{K\in \cK}} K\otimes \rCorr_{A-A}(K\to \cB)\right] = \cB.
    \end{align*}
Thus the inclusion is PQR.
\end{proof}

The following proposition characterizes the PQR property for irreducible inclusions in terms of properties of its  C*-algebra of adjointable endomorphisms. 
   
\begin{prop}\label{prop:RelComm}
Let $A\overset{E}{\subset} B$ be a unital irreducible inclusion of C*-algebras.
Then $A\overset{E}{\subset} B$ is projective quasi-regular if and only if the C*-subalgebra  $[A\rhd]'\cap \End^\dag(\cB_A) = \End^\dag({}_A\cB_A)\subset \End^\dag(\cB_A)$ admits a family of central projections $\{p_i\colon i\in I\}$ satisfying
    \begin{enumerate}[label=(\arabic*)]
    \item $\sum_{i\in I} p_i =1$ with convergence in the strict topology;

    \item for each $i\in I$ there exists $n_i\in\bbN$ so that $p_i \End^\dag({}_A\cB_A) \cong M_{n_i}(\bbC)$;

    \item for every finite-rank projection $p\in\End^\dag({}_A\cB_A)$, the $A$-$A$ bimodule $p[\cB]\in\fgpBim(A)$.
    \end{enumerate}
In this case, $\End^\dag({}_A\cB_A)$ is $*$-isomorphic to the von Neumann algebra direct sum $\bigoplus_{i\in I} M_{n_i}(\bbC)$.
\end{prop}
\begin{proof}
$(\Rightarrow)$: Using the decomposition in Corollary~\ref{cor:P-W}, for $K\in \cK$ let $p_K\in \End^\dag({}_A\cB_A)$ be the projection on $K^{n_K}$. Then (1) is immediate. Irreducibility for $K,K'\in \cK$ gives
    \[
        p_K \End^\dag({}_A\cB_A) p_{K'} = \rCorr_{A-A}((K')^{\oplus n_{K'}}\to K^{\oplus n_K}) \cong \begin{cases} M_{n_K}(\bbC) & \text{if }K=K' \\ \{0\} & \text{otherwise} \end{cases}.
    \]
This gives (2) when $K=K'$. For  distinct $K,K'\in \cK$ and $T\in \End^\dag({}_A\cB_A)$, the above gives
    \[
        p_K T = p_K T \sum_{K'\in \cK} p_{K'} = p_K T p_K = \sum_{K'\in \cK} p_{K'} T p_K = T p_K.
    \]
Thus $p_K$ is a central projection. Finally, notice that each finite rank projection $p\in \End^\dag({_A}\cB_A)$ must be a finite sum of (subprojections of) the $p_i$'s and thus yields some $\fgpBim(A)\ni p[\cB]\subset B\Omega.$ Statement (3) is thus established.  

\noindent$(\Leftarrow)$: Let $\{p_i\in I\colon I\}$ be a family of central projections satisfying (1)-(3). Condition (2) and (3) imply that for each $i\in I$ there exists $K_i\in \fgpBim(A)$ irreducible so that $p_i[\cB] \cong K_i^{\oplus n_i}$. We claim
    \[
        n_i = \mathsf{dim}_{\bbC}\rCorr_{A-A}(K_i\to \cB).
    \]
Indeed, $p_i[\cB] \cong K_i^{\oplus n_i}$ implies we can find embeddings $\iota_1,\ldots, \iota_{n_i}\in \rCorr_{A-A}(K_i\to \cB)$ that are orthonormal with respect to the Hilbert space structure (see Lemma~\ref{lem:Hilbert_space_structure}) and satisfy $p_i[\cB]=\sum_{j=1}^{n_i} \iota_j[K_i]$. For
    \[
        f\in \text{span}\{\iota_1,\ldots, \iota_{n_i}\}^\perp,
    \]
we have $f[K_i] \perp p_i[\cB]$. In fact, irreducibility implies $f[K_i]\perp p_{i'}[\cB]$ for all $i'\in I$, and so
    \[
        f[K_i] = \sum_{i'\in I} p_{i'} f[K_i] = 0.
    \]
by condition (1). Thus $f=0$ and therefore
    \[
        \mathsf{dim}_{\bbC}\rCorr_{A-A}(K_i\to \cB) = \mathsf{dim}_{\bbC}\text{span}\{\iota_1,\ldots, \iota_{n_i}\}= n_i.
    \]
It then follows from Corollary~\ref{cor:P-W} that $A\subset B$ is PQR.
\end{proof}

\section{Crossed products by UTC-actions}\label{sec:ReducedRealization}
\subsection{Actions of UTCs on C*-algebras}

Recall that a functor is said to be \emph{faithful} if it is injective at the level of hom spaces, and is said to be \emph{full} if it is surjective at the level of hom spaces.
An action of a UTC $\cC$ on a unital C*/W*-algebra $A$ is a unitary tensor functor $F: \cC\to\fgpBim(A).$ 
We sometimes denote this by $\cC \overset{F}{\curvearrowright}A.$ 
By untarity, faithfulness of this functor is automatic (since the relevant endomorphism spaces are finite dimensional C*-algebras, and $F(1)=1$). 
However, the property of being full is not, and so we introduce the following definition:
\begin{defn}\label{defn:OuterAction}
    An \textbf{outer action} of a UTC $\cC$ on a unital C*/W*-algebra is a fully-faithful unitary tensor functor $F:\cC\to \fgpBim(A)$.
\end{defn}

In order to justify this nomenclature, we explain how to categorify an outer action by a discrete group on a unital C*-algebra:
\begin{ex}\label{ex:HilbgammOuterAction}
    For any countable discrete group $\Gamma$, we can consider the UTC of finite dimensional $\Gamma$-graded Hilbert spaces $\fdHilb(\Gamma),$ whose morphisms are uniformly bounded bounded linear transformations that respect the grading. 

    An action $\Gamma \overset{\alpha}{\curvearrowright} A$  over a simple unital C*-algebra $A$ consists of the same data as a  unitary tensor functor: 
    \begin{align*}
        \alpha:\fdHilb(\Gamma)&\to\fgpBim(A) \qquad \text{ with tensorator }  &\alpha^2_{g,h}:{}_{g}A\underset{A}{\boxtimes} {}_{h}A \to {}_{gh}A\\
        g&\mapsto {}_{g}A,  &{\xi\boxtimes\eta}\mapsto \alpha_{h^{-1}}(\xi)\cdot\eta.
    \end{align*}
    Here, for each $g\in\Gamma,$ ${}_{g}A\in\fgpBim(A)$ is $A$ as a vector space, with Hilbert C* $A$-$A$ bimodule structure given by $a\rhd\xi\lhd c=\alpha_{g^{-1}}(a)\cdot\xi\cdot c,$ and right and left inner products $\langle\xi\mid  \eta \rangle^{g}_A:=\xi^*\cdot \eta$ and ${}_A\!\langle\xi,\ \eta \rangle^{g}:=\alpha_{g^{1}}(\xi\cdot \eta^*)$ respectively.  
    
    In the usual language of group-actions, if $\{u_g\}_{g\in\Gamma}$ are the usual unitaries implementing the action, ie $\alpha_g(\cdot) = u_g(\cdot) u_g^*,$ outerness for $\alpha$ means that $g\neq h$ implies $u_g\neq u_h.$ 
    So, ${_g}A$ really is the same as $u_gA.$ 
    At the categorified level, this manifests as $g\neq h$ implies ${_g}A\not\cong {_h}A,$ so the categorified action maps distinct group elements to distinct irreducible bimodules.  
    And in greater generality this is precisely achieved by a fully-faithful functor, which will map non-isomorphic irreducible objects in $\cC$ into non-isomorphic simple objects in its target.  
\end{ex}
    
    We shall now discuss further symmetries of operator algebras which do not come from groups. 
    Popa established that every \emph{$\lambda$-lattice} arises as the standard invariant of some finite index subfactor \cite{MR1334479}. 
    Later on, together with Shlyakhtenko \cite{MR2051399} they showed this can be done using the $\rm{II}_1$-factor $\mathsf{L}\bbF_\infty$ associated to the free group in countably many generators, making it a universal receptacle for actions of UTCs \cite{BHP12}. 
    Subsequently, Guionnet, Jones, Shlyakhtenko and Walker provided a new diagrammatic proof of Popa's construction  via Jones' Planar algebras \cite{GJS10, JSW10}.     
    
\begin{ex}
    Using the \emph{GJS construction}, the first-named author and Hartglass show in \cite{MR4139893} that any (countably generated) UTC acts outerly on some unital separable simple exact C*-algebra $A$ with a unique trace $\tr$, whose $K$-theory depends only on $\cC$. 
    The GJS construction relies on diagrammatic and free-probabilistic techniques and yields an $A$-valued semicircular system C*-algebra, whose corners provide the bimodules used to represent $\cC.$ 
    See Section \ref{sec:GJS} for more details. 
\end{ex}

\subsection{Reduced crossed product C*-algebras by UTC-actions}\label{subsec:CrossProd}
In this section we describe how to take reduced crossed products by an action of a UTC $\cC$ on a unital simple C*-algebra $A$ to obtain a reduced C*-algebra or a C*-correspondence. 
The mindset is that $\cC$ acts on $A$ via fgp $A$-$A$ bimodules, and bimodules are generalized endomorphisms. So one should be able to construct extensions of $A$ via a generalized crossed product. 
In the subfactors literature, these are often found as \emph{realizations of algebra objects} (see \cite{MR3948170} and references therein), but these ideas have also appeared in the C*-algebra literature; cf \cite{MR1467459}, where in the case of a simple C*-algebra, they illustrate how to take a crossed product by an imprimitivity bimodule. 
    \begin{construction}[The realized reduced C*-algebra]\label{const:realizedredC*alg}
    Begin with a unital simple C*-algebra $A$, a UTC $\cC$ and a unitary tensor functor $F:\cC \to \Bim_{\sf{fgp}}(A).$
    Then, by forgetting the C*-correspondence structure via the forgetful tensor functor $\sf{For}$, we get ${\sf For}\circ F =: \bbF\in \Vec(\cC^{\sf{op}})$ is a C*-algebra object.
    Suppose that $\bbB\in \Vec(\cC)$ is a connected C*-algebra object, and consider the \emph{algebraic realization} \cite[Definition 4.1]{MR3948170} or \emph{(algebraic) generalized crossed product} from Proposition \ref{prop:PQNAlgRealization},  
    $$(A\rtimes_{F}\bbB)^{\diamondsuit}:=|\bbB|^{\diamondsuit}_{F} := \bigoplus_{c\in\Irr(\cC)} \bbF(c)\otimes \bbB(c),$$
    is an algebraic direct sum of vector spaces. (See also \cite[Remark 2.6]{MR3406647}.) 
    All of these notations and nomenclature are found across the literature, but here we shall elaborate mostly using the language of crossed products to better highlight the connection with discrete groups, their actions and usual crossed products. 
    
    We can endow this space with the structure of a unital associative $*$-algebra using the graphical calculus for $\Vec(\cC)$ afforded by the Yoneda embedding (cf.  Remark \ref{remark:GraphicalCalculus}). 
    (For a rigorous description of these constructions, we refer the reader to \cite[\S4.1]{JP17}.)   
    For $a\in\cC,$ a generic element in $(A\rtimes_{F}\bbB)^{\diamondsuit}$ is represented graphically as
\begin{align}
    \eta\otimes g =
    \tikzmath{
        \node at(0,1.4){$\scriptstyle \bbB$};
        \draw(0,.9) -- (0,1.2);
        \roundNbox{unshaded}{(0,.6)}{.3}{0}{0}{$g$};
        \draw[dotted] (-.6,0)--(.6,0);
        \node at(-.15,-.15){$\scriptstyle\mathbb{a}$};
        \draw (0,.3) -- (0,-.3);
        \node at(-.15,.15){$\scriptstyle\mathbb{a}$};
        \roundNbox{unshaded}{(0,-.6)}{.3}{0}{0}{$\eta$};
        \node at(0,-1.4){$\scriptstyle \bbF$};
        \draw(0,-.9) -- (0,-1.2);
    }.
\end{align}
    The horizontal dotted line stands for the \emph{permeable membrane} that allows to pass morphisms from $\cC$ back and forth (and reversing directions for composition). 
    We use this membrane to separate the two different but simultaneous conventions for the graphical calculi: above, for $\Vec(\cC)$ we read from the membrane to the top, and below the membrane, for $\Vec(\cC^{\op}),$ we read from the membrane to the bottom.

The (associative) multiplication for arbitrary $\sum_{c\in\Irr(\cC)}\xi_{(c)}\otimes f_{(c)}$ and $\sum_{d\in\Irr(\cC)}\eta_{(d)}\otimes g_{(d)}$ is given by:
\begin{align}\label{eqn:GradedMult}
	\sum_{c\in\Irr(\cC)}\xi_{(c)}&\otimes f_{(c)}\ \cdot \sum_{d\in\Irr(\cC)}\eta_{(d)}\otimes g_{(d)}
	=\sum_{c,d\in\Irr(\cC)}
	\tikzmath{
	    \node at (0.5,3){$\scriptstyle\bbB$};
	    \node at(1.15,1.4){$\scriptstyle \bbB$};	
	    \node at(-.15,1.4){$\scriptstyle \bbB$};	
	    \node at(-.15,.15){$\scriptstyle \mathbb{c}$};
	    \node at(-.15,-.15){$\scriptstyle \mathbb{c}$};
        \draw (.5,2.4)--(.5,2.8);
	    \draw (0,1.2) -- (0,1.4) arc(180:0:.5) -- (1,1.2);
	    \roundNbox{unshaded}{(.5,2)}{.4}{0}{0}{$\bbB^2_{c,d}$};
	    \roundNbox{unshaded}{(0,.8)}{.4}{0}{0}{$f_{(c)}$};
	    \draw(0,.4) --(0,-.4);
	    \roundNbox{unshaded}{(0,-.8)}{.4}{0}{0}{$\xi_{(c)}$};
	    \draw (0,-1.2) -- (0,-1.5);
	    \draw[dotted] (-.6,0) --(1.6,0);
	    \roundNbox{unshaded}{(1,.8)}{.4}{0}{0}{$g_{(d)}$};
	    \draw(1,.4) --(1,-.4);
	    \roundNbox{unshaded}{(1,-.8)}{.4}{0}{0}{$\eta_{(d)}$};
	    \draw (.5,-2.4)--(.5,-2.8);
	    \draw (0,-1.2) -- (0,-1.4) arc(180:360:.5) -- (1,-1.2);
	    \roundNbox{unshaded}{(.5,-2)}{.4}{0}{0}{$\bbF^2_{c,d}$};
	    \node at(1.15,.15){$\scriptstyle \mathbb{d}$};
	    \node at(1.15,-.15){$\scriptstyle \mathbb{d}$};
	    \node at(-.15,-1.4){$\scriptstyle \bbF$};
	    \node at(1.15,-1.4){$\scriptstyle \bbF$};
	    \node at (0.5,-3){$\scriptstyle \bbF$};
	}
	= \sum_{\substack{{c,d,e\in\Irr(\cC)}\\{\alpha\in\Isom(e, c\otimes d)}}}
	\tikzmath{
	    \node at(.5,3.5){$\scriptstyle\bbB$};
    	\draw (0,1.8) -- (0,2) arc(180:0:.5) -- (1,1.8);
	    \roundNbox{unshaded}{(.5,2.6)}{.4}{0}{0}{$\bbB^2_{c,d}$};
	    \roundNbox{unshaded}{(0,1.4)}{.4}{0}{0}{$f_{(c)}$};
	    \roundNbox{unshaded}{(0,-1.4)}{.4}{0}{0}{$\xi_{(c)}$};
	    \draw (.5,3)--(.5,3.3);
        \node at(-.15,.8){$\scriptstyle\mathbb{c}$};
        \node at(1.15,.8){$\scriptstyle\mathbb{d}$};
	    \draw(0,1) -- (0,.9) arc(180:360:.5)--(1,1);
	    \draw(.5,.4) -- (.5,-.4);
	    \filldraw[fill=white] (.5,.4) circle (3pt);
	    \node at(.35,-.15){$\scriptstyle\mathbb{e}$};
	    \node at(.35,.15){$\scriptstyle\mathbb{e}$};
        \draw[dotted] (-.6,0) --(1.6,0);
        \draw(0,-1) -- (0,-.9) arc(180:0:.5)--(1,-1);
        \filldraw[fill=white] (.5,-.4) circle (3pt);
        \node at(-.15,-.8){$\scriptstyle\mathbb{c}$};
        \node at(1.15,-.8){$\scriptstyle\mathbb{d}$};
        \roundNbox{unshaded}{(1,1.4)}{.4}{0}{0}{$g_{(d)}$};
	    \roundNbox{unshaded}{(1,-1.4)}{.4}{0}{0}{$\eta_{(d)}$};
	    \draw (.5,-3)--(.5,-3.3);
	    \draw (0,-1.8) -- (0,-2) arc(180:360:.5) -- (1,-1.8);
	    \roundNbox{unshaded}{(.5,-2.6)}{.4}{0}{0}{$\bbF^2_{c,d}$};
	    \node at(.5,-3.5){$\scriptstyle \bbF$};
	}\nonumber \\
	&=\sum_{\substack{{c,d,e\in\Irr(\cC)}\\{\alpha\in\Isom(e, c\otimes d)}}} \bbF(\alpha^*)\circ[\bbF^2_{c,d}(\xi_{(c)}\otimes\eta_{(d)})]\otimes\bbB(\alpha)\circ[\bbB^2_{c,d}(f_{(c)}\otimes g_{(c)})].
	\end{align}
	Here, for $c,d,e\in\Irr(\cC),$ $\Isom(e, c\otimes d)$ denotes a chosen complete set of isometries with orthogonal ranges, which therefore satisfy $\displaystyle\sum_{\substack{\alpha\in\Isom(e,c\otimes d)}}\alpha\circ\alpha^*=\id_{c\otimes d}.$ 
	The empty circles in the picture above stand for such $\alpha$ and $\alpha^*$.
	The unit of this algebra is $\bbF^1\otimes \bbB^1\in \bbF(1_\cC)\otimes \bbB(1_\cC).$ (Recall notation from discussion following Definition \ref{defn:AlgObjs}.)
	
	The $*$-structure is given by 
	\begin{align}
	    (\eta\otimes g)^*=
	    \sum_{d\in\Irr(\cC)}
	    \tikzmath{
        \node at(0,1.7){$\scriptstyle \bbB$};
        \draw(0,1.2) -- (0,1.5);
        \roundNbox{unshaded}{(0,.8)}{.4}{.5}{.2}{$j^{\bbB}_d(g_{(d)})$};
        \draw[dotted] (-.6,0)--(.6,0);
        \node at(-.15,-.15){$\scriptstyle\overline{\mathbb{d}}$};
        \draw (0,.4) -- (0,-.4);
        \node at(-.15,.15){$\scriptstyle\overline{\mathbb{d}}$};
        \roundNbox{unshaded}{(0,-.8)}{.4}{.5}{.2}{$j^{\bbF}_d(\eta_{(d)})$};
        \node at(0,-1.7){$\scriptstyle \bbF$};
        \draw(0,-1.2) -- (0,-1.5);
    }.
	\end{align}
\end{construction}

Our goal is to turn $(A\rtimes_{F}\bbB)^{\diamondsuit}$ into a C*-algebra containing a faithful copy of $A$ unitally. 
Consider the map 
\begin{align}\label{align:E}\begin{split}
	E^{\diamondsuit}:(A\rtimes_{F}\bbB)^{\diamondsuit} &\twoheadrightarrow A\\
	\sum_{c\in\Irr(\cC)} h_{(c)} \otimes b_{(c)} &\mapsto h^{ \bbOne_\cC}_{(1)}\otimes b^{\bbOne_\cC}_{(2)}=h_{ \bbOne_\cC}\otimes b_{\bbOne_\cC},
\end{split}\end{align}
where we dropped the legs of Sweedler's notation in the last equality  because $\bbB$ is connected. 
Consider $(A\rtimes_{F}\bbB)^{\diamondsuit}$ as a right  $A$-module with the (potentially degenerate) right $A$-valued  product $\langle x\ |\ y\rangle_A :=E^{\diamondsuit}(x^*y),$ and let $\cB$ be its completion after taking quotients by negligible vectors.  
After noticing we have the following inclusions $A\Omega\subset(A\rtimes_{F}\bbB)^{\diamondsuit}\Omega\subset \cB,$ we record the following Lemma:
\begin{lem}\label{lem:bddExpectation}
    The map $E^{\diamondsuit}$ is a unital completely positive, $A$-$A$ bimodular and faithful map. 
    Moreover, for any $y\in (A\rtimes_{F}\bbB)^{\diamondsuit}$ there exists $C_y\in[0,\infty),$ such that for every $x, y\in(A\rtimes_{F}\bbB)^{\diamondsuit}:$
    $$\langle y\rhd(x\Omega)\ |\ y\rhd(x\Omega) \rangle_A = E^{\diamondsuit}(x^*y^*yx)\leq C_yE^{\diamondsuit}(x^*x).$$
\end{lem}
\begin{proof}
    Since $A\Omega\subset \cB$ is fgp, by \cite[Lemma 2.3.7]{MR2125398} it is orthogonally complementable, and so the projection $e_A: \cB\to \cB$ with rank $A\Omega$ is adjointable. 
    By means of computation, we obtain that for any $x,y\in (A\rtimes_{F}\bbB)^{\diamondsuit},$ we have $E^{\diamondsuit}(x)e_A(y\Omega) = e_A x e_A (y\Omega), $ so  $E^{\diamondsuit}$ is implemented by $e_A.$ Thus $E^{\diamondsuit}$ is ucp.  
    Since $e_A$ commutes with the right $A$-action and left multiplication by elements in $A$, it follows that $E^{\diamondsuit}$ is $A$-$A$ bimodular. 
    We shall now verify the faithfulness of $E^{\diamondsuit}.$
    For arbitrary $x\in (A\rtimes_F\bbB)^{\diamondsuit}$, express $x=\sum_{c\in\Irr(\cC)}x_{(c)}\otimes f_{(c)}$ as a finite sum. 
    Since $\bbF$ is a C*-algebra object in $\Vec{(\cC^{\op})}$, by \cite[Proposition 3.3]{JP17}, each $\bbF(\overline {c} \otimes c)$ is a C*-algebra, and so 
    \begin{align*}
        E^{\diamondsuit}(x^*x) &= \sum_{c\in \Irr(\cC)} \left\langle \Omega \mid  \left[\bbF^2_{\overline c \otimes c}\left(j^{\bbF}_c \left(x_{(c)}\right)\boxtimes x_{(c)}\right)\right]\otimes \left[\bbB^2_{\overline c\otimes c}\left(j^{\bbB}_c\left(f_{(c)}\right)\boxtimes f_{(c)}\right)\right] \Omega\right\rangle_A\\
        &=\sum_{c\in\Irr(\cC)}
        \tikzmath{
	    \node at (0.5,3){$\scriptstyle\bbB$};
	    \node at(1.15,1.7){$\scriptstyle \bbB$};	
	    \node at(-.15,1.7){$\scriptstyle \bbB$};	
	    \node at(-.05,.35){$\scriptstyle \overline{\mathbb{c}}$};
	    \node at(-.05,-.35){$\scriptstyle \overline{\mathbb{c}}$};
        \draw (.5,2.4)--(.5,2.8);
	    \draw (0,1.2) -- (0,1.4) arc(180:0:.5) -- (1,1.2);
	   \draw(0.1,.66) arc(180:360:.4);
	                \draw[dotted] (-1.6,0) --(2,0);
	                \draw[dashed] (.5,.25) --(.5,-.25);
	   \draw(0.1,-.66) arc(180:0:.4);
	    \draw (.5,-2.4)--(.5,-2.8);
	    \draw (0,-1.2) -- (0,-1.4) arc(180:360:.5) -- (1,-1.2);
	    \node at(.95,.35){$\scriptstyle \mathbb{c}$};
	    \node at(.95,-.35){$\scriptstyle \mathbb{c}$};
	    \node at(-.15,-1.7){$\scriptstyle \bbF$};
	    \node at(1.15,-1.7){$\scriptstyle \bbF$};
	    \node at (0.5,-3){$\scriptstyle \bbF$};
	    \roundNbox{unshaded}{(1,-1)}{.4}{0}{0}{$x_{(c)}$};
	    \roundNbox{unshaded}{(1,1)}{.4}{0}{0}{$f_{(c)}$};
	    \roundNbox{unshaded}{(.5,-2)}{.4}{0}{0}{$\bbF^2_{c,c}$};
	    \roundNbox{unshaded}{(.5,2)}{.4}{0}{0}{$\bbB^2_{c,c}$};
	    \roundNbox{unshaded}{(-.2,1)}{.4}{.7}{0}{$j^\bbB_{c}\!\left(f_{(c)}\right)$};
	    \roundNbox{unshaded}{(0,-1)}{.4}{.7}{0}{$j^\bbF_c\!\left(x_{(c)}\right)$};
	}\ 
	\geq 0,
    \end{align*}
    since it is a sum of positive traced-out elements in $\bbF(\overline c\otimes c)$ multiplied by non-negative real numbers. 
    Indeed, the bottom and top diagrams are images of the categorical traces in the cyclic $\cC$-module C*-categories $\fM_\bbF$ and $\fM_\bbB$ determined by $\bbF$ and $\bbB$, respectively. And for each $c\in\Irr(\cC),$ $\Tr^\bbF_c: \bbF(\overline{c}\otimes c)\to A\cong\bbF(1_\cC)$ and $\Tr^\bbB_c: \bbB(\overline{c}\otimes c)\to \bbC\cong\bbB(1_\cC)$ are positive-definite linear maps.  
    Thus, the expression above is equal to zero if and only if for every $c\in\Irr(\cC)$ the $c$-components of $x$ are zero. Thus $E^{\diamondsuit}$ is faithful.
	  
    To prove the inequality, we make use of the graphical representation of the $\cC$ graded multiplication $x^*y^*yx$ from Equation \ref{eqn:GradedMult} (and replacing tensorators by bullets to save space):
	\begin{align}\label{align:BoundedMult}\begin{split}
    \|y \rhd x\Omega\|_A^2=E^{\diamondsuit}(x^*y^*yx)&=E^{\diamondsuit}\left(\underset{d,d'\in\Irr(\cC)}{\sum}\ \ 
    \tikzmath{
        \node at(-.7,4){$\scriptstyle\bbB$};
        \node at(-.7,3.35){$\bullet$};
        \draw(-.7,3.3)--(-.7,3.8);
        \draw (-2.4,2.8) .. controls (-2.2,3.5) and (.6,3.5) .. (.8,2.9);
        \roundNbox{unshaded}{(-2.2,2.6)}{.3}{.75}{0}{$\scriptstyle j^\bbB_{d'}\left(f_{(d')}\right)$};
        \roundNbox{unshaded}{(1.6,2.25)}{.3}{.1}{0}{$\scriptstyle f_{(d)}$};
            \draw[dashed, cyan, rounded corners] (-2.3,-1.9) rectangle(1.2,1.9); 
        \node at(.8,2.9) {$\bullet$};
        \node at(.1,1.7) {$\bullet$};
        \draw (.1,1.7) -- (.1,2.2) arc(180:30:.7);
        \draw (-.6,1) arc(180:0:.7);
        \node at(-1.6,0){$\scriptstyle \underset{c,c'\in\Irr(\cC)}{\sum}$};
        \draw (-.6,-1) arc(180:360:.7);
        \roundNbox{unshaded}{(-.3,.8)}{.3}{.85}{0}{$\scriptstyle j^{\bbB}_{c'}\left(g_{(c')}\right)$};
        \roundNbox{unshaded}{(.8,.8)}{.3}{.1}{0}{$\scriptstyle g_{(c)}$};
            \draw(-2.5,2.3) --(-2.5,-2.3);
            \draw (1.6,1.95) -- (1.6,-1.95);
            \node at(1.75,.25){${\scriptstyle\mathbb{d}}$};
            \node at(-2.7,.25){${\scriptstyle\overline{\mathbb{d'}}}$};
                \draw[dotted] (-3.4,0)--(2.3,0);
            \draw (-.6,.5) -- (-.6,-.5);
            \node at(.65,.25){${\scriptstyle\mathbb{c}}$};
        	\node at(-.4,.25){${\scriptstyle\overline{\mathbb{c'}}}$};
        	\node at(-.4,-.25){$\scriptstyle\overline{\mathbb{c'}}$};
        	\node at(.65,-.25){${\scriptstyle\mathbb{c}}$};
        	\node at(1.75,-.25){${\scriptstyle\mathbb{d}}$};
            \node at(-2.7,-.25){${\scriptstyle\overline{\mathbb{d'}}}$};
        	\draw (.8,.5) -- (.8,-.5);
        \roundNbox{unshaded}{(-.3,-.8)}{.3}{.85}{0}{$\scriptstyle j^{\bbF}_{c'}\left(y_{(c')}\right)$};
        \roundNbox{unshaded}{(.8,-.8)}{.3}{.1}{0}{$\scriptstyle y_{(c)}$};    
        \draw (.1,-1.7) -- (.1,-2.2) arc(180:330:.7);
        \node at(.1,-1.7) {$\bullet$};
        \node at(.8,-2.9) {$\bullet$};
        \roundNbox{unshaded}{(1.6,-2.25)}{.3}{.1}{0}{$\scriptstyle x_{(d)}$};
        \roundNbox{unshaded}{(-2.2,-2.6)}{.3}{.75}{0}{$\scriptstyle j^\bbF_{d'}\left(x_{(d')}\right)$};
        \draw (-2.4,-2.9) .. controls (-2.2,-3.5) and (.6,-3.5) .. (.8,-2.9);
        \draw(-.7,-3.3)--(-.7,-3.8);
        \node at(-.7,-4){$\scriptstyle\bbF$};
        \node at(-.7,-3.35){$\bullet$};
    }\right)\\
    &\leq \left\|\sum_{c\in\Irr(\cC)} y_{(c)}\right\|^2_{\fM_{\bbF}}\cdot \left\|\sum_{c\in\Irr(\cC)} g_{(c)}\right\|^2_{\fM_{\bbB}}\cdot E^{\diamondsuit}(x^*x).
	\end{split}\end{align}
	Here, we removed the positive diagram contained inside the cyan dashed box corresponding to $y^*y,$ and replaced it with the larger positive $\left\|\sum_c y_{(c)}\right\|^2_{\fM_{\bbF}}\cdot \left\|\sum_c g_{(c)}\right\|^2_{\fM_{\bbB}}\cdot\bbF^1\otimes\bbB^1,$ where the norms are computed inside the corresponding C*-module categories, $\fM_\bbF$ and $\fM_\bbB$. 
	Since $y$ is finitely supported, then these norms are finite, completing the proof. 
\end{proof}

The inequality in Lemma \ref{lem:bddExpectation} implies that left multiplication by elements in $(A\rtimes_{F}\bbB)^{\diamondsuit}$ defines a bounded operator on $\cB$. 
It is also easy to check that these left multiplication operators are adjointable, providing us with a faithful $*$-representation 
$-\rhd: (A\rtimes_{F}\bbB)^{\diamondsuit}\hookrightarrow \End^\dag_{\bbC-A}(\cB_A).$
We then define the \emph{reduced C*-realization}/\emph{generalized reduced C*-crossed product}: $$A\underset{F,r}{\rtimes}\bbB := \overline{\left[(A\rtimes_F \bbB)^{\diamondsuit}\right]\rhd}^{\|\cdot\|_{\op}},$$ 
the operator norm-completion of the image of $(A\rtimes_{F}\bbB)^{\diamondsuit}$ inside  $\End_{\bbC-A}^\dag(\cB_A).$ 
Therefore, the following norm comparison will hold on $(A\rtimes_{F}\bbB)^{\diamondsuit}$:
    $$\|b\Omega\|_A\leq \|b\|_{\op}.$$
We can further extend $E^{\diamondsuit}$ from (\ref{align:E}) to a conditional expectation $E:A\rtimes_{F,r}\bbB\to A.$  
This yields $\cB\in\rCorr(A\rtimes_{F,r}\bbB\to A).$ To show this conditional expectation is always faithful, we will need a fine-tuned version of the Pimsner-Popa inequality:

\begin{lem}\label{lem:PPIneq}
    Let $A$ be a unital C*-algebra with trivial center and $F:\cC\to\fgpBim(A)$ be an action. Let $\bbB$ be a connected $\cC$-graded C*-algebra in $\Vec(\cC)$ and $K\in\Irr(\fgpBim(A)$ with a finite right Pimsner-Popa basis $\{u_i\}_1^n$.   
    If $f\in\rCorr_{A-A}(K\to \cB)$ is nonzero and $f[K]\subset (A\rtimes_{F}\bbB)^\diamondsuit\Omega$, then $\{f(u_i)/\|f\|\}_1^n\subset f[K]\in\fgpBim(A)$ is a right Pimsner--Popa basis for the image of $f,$ with $\sum_1^n \check{f}(u_i)\check{f}(u_i)^*\in (A'\cap A\rtimes_{F,r}\bbB)^+\setminus\{0\}.$ 
    Furthermore, for each $k\in K$ we have a Pimsner--Popa type inequality:
    $$\dfrac{\|f\|}{\|\sum_1^n \check{f}(u_i)\check{f}(u_i)^*\|^{1/2}}\|\check{f}(k)\|\leq \|k\|_A\leq\dfrac{\|\check{f}(k)\|}{\|f\|}.$$
Consequently, for $L\in \fgpBim(A)$ satisfying $L\subset (A\rtimes_{F}\bbB)^\diamondsuit\Omega$ there exists a constant $\alpha_L\in(0, \infty)$ so that $\|x\Omega\|_A \leq \|x\| \leq \alpha_L \|x \Omega\|_A$ for all $x\Omega \in L$.
\end{lem}
\begin{proof}
    Denote $(A\rtimes_{F}\bbB)^\diamondsuit =: B^\diamondsuit$ and $A\rtimes_{F,r}\bbB=:B$. 
    Direct computation shows that $\sum_1^n \check{f}(u_i)\check{f}(u_i)^*\in (A'\cap B^\diamondsuit)^+.$ 
    The irreducibility of $K$ implies $f$ is injective, as $K$ being fgp gives $\mathsf{Ker}(f)$ as a complemented sub-module by Lemma~\ref{lem:fgp_properties}.\ref{part:fgp_closed_comp_kernel}, so in particular, the sum above is nonzero. 
    Clearly, $\left\{f(u_i)/{\|f\|}\right\}_1^n,$ the normalized image of a Pimsner--Popa basis under $f$ is a Pimsner--Popa basis for its image. 
    The rightmost inequality follows from
        \[
            \<f(k)\mid f(k)\>_A = (f\mid f)_{\bbC} \< k\mid k\>_A = \|f\|^2 \<k\mid k\>_A,
        \]  
    (see Lemma~\ref{lem:Hilbert_space_structure}) and the fact that $\|f(k)\|_A =\|E(\check{f}(k)^*\check{f}(k))\|^{1/2} \leq \|\check{f}(k)\|$. It remains to prove the first inequality.
    
    The proof is similar to \cite[Proposition 1.18]{MR1624182}, and we spell it out here to make the necessary adaptations. View $B^\diamondsuit\Omega$ as an algebraic $B^\diamondsuit$-bimodule under multiplication, with $B$-valued inner product  ${_B}\<x\Omega, y\Omega\> = xy^*$ for $x,y\in B^\diamondsuit$. Notice this expression is meaningful as $E^\diamondsuit$ is faithful on $B^\diamondsuit$ by Lemma \ref{lem:bddExpectation}, and also that $B^\diamondsuit$ is not necessarily closed with the induced norm.  
    Restricting this inner product to $f[K]$ does \emph{not} make it into a left Hilbert $B$-module (we have not even assumed $f[K]$ is invariant under the left/right $B^\diamondsuit$-action), but one does still have ${_B}\<f(k)\lhd b, f(k)\lhd b\> = {_B}\<\check{f}(k) b\Omega, \check{f}(k)b\Omega \> \leq {_B}\<f(k), f(k)\> \|b\|^2$ for all $k \in K$ and $b\in B^\diamondsuit$ since this holds on $B$. 
    Now, we first assume the right Pimsner--Popa basis for $K$ consists of a single element $u_1$. Then for all $k\in K$ we have
        \begin{align*}
            {_B}\<f(k), f(k)\>  &=\|f\|^{-4}\cdot {_B}\< f(u_1)\lhd \<f(u_1)\mid  f(k)\>_A, f(u_1)\lhd\<f(u_1)\mid  f(k)\>_A \> \\
            &\leq \|f\|^{-4}\cdot {_B}\<f(u_1), f(u_1)\> \|\<f(u_1)\mid  f(k)\>_A\|^2\\ 
            &=\|f\|^{-2}\cdot {_B}\<f(u_1), f(u_1)\> \|\<u_1\mid k\>_A\|^2\\ 
            &=\|f\|^{-2}\cdot {_B}\<f(u_1), f(u_1)\> \|k\|_A^2,
        \end{align*}
    where in the last equality we have used
        \[
            \langle k\mid k\rangle_A = \langle k \mid u_1\lhd\langle u_1\mid k\rangle_A\rangle_A=\langle k\mid u_1\rangle_A \langle u_1 \mid k\rangle_A.
        \]
    Thus
        \[
            \|f\|\cdot\|\check{f}(k)\| \leq \|\check{f}(u_1)\check{f}(u_1)^*\|^{1/2} \cdot  \|k\|_A.
        \]
    Note that this establishes the inequality when $n=1.$ 
    To prove the general case, when $n\in\bbN$ is arbitrary we consider $f[K]^{\oplus n} \subset (B^\diamondsuit\Omega)^{\oplus n}$, which we view as a right Hilbert $M_n(A)$-module where $\< (\xi_1,\ldots,\xi_n) \mid (\eta_1,\ldots, \eta_n)\>_{M_n(A)}\in M_n(A)$ for $\xi_1,\ldots, \xi_n, \eta_1,\ldots, \eta_n\in f[K]$ is the matrix with $(i,j)$-entry given by $\<\xi_i\mid \eta_j\>_A$. We also impart $f[K]^{\oplus n}$ with the $B$-valued inner product
        \[
            {_B}\< (\xi_1,\ldots,\xi_n), (\eta_1,\ldots, \eta_n)\> = \sum_{i=1}^n {_B}\< \xi_i, \eta_i \> \qquad\qquad \xi_1,\ldots, \xi_n, \eta_1,\ldots, \eta_n\in f[K].
        \]
    Set $\vec{u} := \|f\|^{-1}\cdot(f(u_1), f(u_2),\hdots, f(u_n))\in f[K]^{\oplus n},$ and observe that for \break
    $\vec{k} = (f(k_1),f(k_2),\hdots, f(k_n))\in f[K]^{\oplus n}$ one has
        \[
            \left[\vec{u}\lhd \<\vec{u}\mid \vec{k}\>_A\right]_i = \sum_{j=1}^n u_j\lhd \< u_j\mid k_i\>_A = k_i \qquad i=1,\ldots, n.
        \]
    Thus $\vec{u}\lhd \<\vec{u}\mid \vec{k}\>_A = \vec{k}$, and so the inequality for the $n=1$ case yields
        \begin{align*}
            \left\| \sum_{i=1}^n \check{f}(k_i)\check{f}(k_i)^* \right\|^{1/2} &= \|{_B}\<\vec{k}, \vec{k}\>\|^{1/2}\\
                &\leq \| {_B}\<\vec{u}, \vec{u}\>\|^{1/2}\cdot \|\vec{k}\|_{M_n(A)} = \left\| \sum_{i=1}^n \check{f}(u_i) \check{f}(u_i)^* \right\|^{1/2}  \dfrac{\|\vec{k}\|_{M_n(A)}}{\|f\|}.
        \end{align*}
    Applying this to $\vec{k}=(f(k), 0,\ldots, 0)$ for $k\in K$ then yields the claimed inequality.

    Finally, consider a (non-irreducible) $L\in \fgpBim(A)$ satisfying $L\subset (A\rtimes_{F}\bbB)^\diamondsuit\Omega$. Decompose $L=\bigoplus_{j=1}^m L_j$ into irreducible $L_j\in \fgpBim(A)$, and let $\iota_j\in \rCorr_{A-A}(L_j\to L)$ be the embedding. Fix $x\Omega\in L$ and let $\ell_j\in L_j$ be such that $x\Omega = \sum_{j=1}^m \iota_j(\ell_j)$. Then the above gives
        \[
            \|x\| \leq \sum_{j=1}^m \| \check{\iota_j}(\ell_j)\| \leq \sum_{j=1}^m \beta_j \|\iota_j(\ell_j)\|_A \leq \sum_{j=1}^m \beta_j \|x\Omega\|_A,
        \]
    where $\beta_j>0$ comes from a choice of Pimsner--Popa basis for $L_j$. So taking $\alpha_L:= \sum_{j=1}^m \beta_j$ completes the proof.
\end{proof}
In practice, Lemma~\ref{lem:PPIneq} tells us that the norm on $(A\rtimes_F\bbB)^\diamondsuit \subset A\rtimes_{F,r}\bbB$ can be controlled by the norm on $(A\rtimes_F\bbB)^\diamondsuit \Omega \subset \cB$ so long as one remains localized to a fixed $\fgpBim(A)\ni K\subset (A\rtimes_F\bbB)^\diamondsuit\Omega$.

\begin{prop}
The conditional expectation $E:A\rtimes_{F,r}\bbB\to A$ defined above is faithful.
\end{prop}
\begin{proof}
Suppose $x\in A\rtimes_{F,r}\bbB$ has $E(x^*x)=0$. Let $(x_n)_{n\in \bbN}\in (A\rtimes_{F,r}\bbB)^\diamondsuit$ be sequence converging in norm to $x$. For each $n\in \bbN$ we can write
    \[
        x_n = \sum_{c\in \Irr(\cC)} x_n^{(c)}
    \]
where the sum is finite and $x_n^{(c)}\in \bbF(c)\otimes \bbB(c)$. Now, $E(x^*x)=0$ is equivalent to $x\Omega =0$ and therefore
    \[
        x_n^{(c)} \Omega = P_c x_n\Omega \to P_c x\Omega =0.
    \]
However, since $\bbB$ is connected then all its fibers are finite dimensional \cite[Proposition 2.8]{MR3948170} and so $F(c)\otimes \bbB(c)\in \fgpBim(A)$. Applying the last part of Lemma~\ref{lem:PPIneq} to $F(c)\otimes \bbB(c)\subset (A\rtimes_{F}\bbB)^\diamondsuit\Omega$ we then obtain $\|x_n^{(c)}\|\to 0$.

Now, fix $d,e\in \Irr(\cC)$ and note that $\supp_{\cC}(e\otimes\bar{d}):=\{c\in \Irr(\cC)\colon \Hom(e\otimes \bar{d}\to c)\neq 0\}$ is finite. Also, for $y\in \bbF(d)\otimes \bbB(d)$ one has
    \[
        P_e x_n y\Omega = \sum_{c\in \supp_{\cC}(e\otimes\bar{d})} P_e x_n^{(c)} y\Omega.
    \]
Since the last sum is finite and each $x_n^{(c)}\to 0$ in norm, it follows that
    \[
        P_e x y\Omega = \lim_{n\to\infty} P_e x_n y\Omega = 0.
    \]
Using the $\sum_e P_e$ converges strictly to $1$ on $\cB$, we then have
    \[
        xy\Omega = \sum_{e\in \Irr(\cC)} P_e xy\Omega =0.
    \]
Thus $x=0$ since $\bigcup_d \bbF(d)\otimes \bbB(d)\Omega$ has dense span.
\end{proof}

\subsection{Discrete C*-inclusions and various examples}\label{sec:EgsCDisc}
Inspired by this construction, we shall now precisely define what it means for a unital inclusion of C*-algebras to be \emph{discrete}. Of course, outputs of our generalized reduced crossed product constitute the largest class of examples. Later on, in Section \ref{sec:Characterization} we will provide an abstract characterization of these, showing all irreducible  C*-discrete inclusions arise in this fashion.
\begin{defn}\label{defn:CDisc}
    We say a unital inclusion $A\overset{E}{\subset} B$ of C*-algebras is \textbf{C*-discrete} if $\overline{B^{\diamondsuit}}^{\|\cdot\|_B}= \overline{\PQR(A\subset B)}^{\|\cdot\|_B}=B.$
\end{defn}

We have thus proven a C*-analogue of \cite[Theorem 5.8, Proposition 5.10]{MR3948170}
\begin{thm}\label{thm:CrossProdPQR}
	Let A be a unital C*-algebra with trivial center, $\cC$ be a UTC, $F:\cC \to \Bim_{\sf{fgp}}(A)$ a unitary tensor (full) functor, and $\bbB$ a connected C*-algebra object in $\Vec(\cC).$
	Then $A\overset{E}{\subset} A\rtimes_{F,r}\bbB$ is an irreducible C*-discrete inclusion of C*-algebras with a faithful ucp $A$-$A$ bimodular conditional expectation $E: A\rtimes_{F,r}\bbB \twoheadrightarrow A$ projecting onto the $1_\cC$-component.
\end{thm}

Notice that by the usual norm comparison $\|b\Omega\|_A\leq\|b\|_B$ on $B,$ \textbf{C*-discreteness implies projective-quasi-regularity} from Definition~\ref{defn:PQR}. (We currently do not know if the converse holds.) Thus, all crossed products by UTC actions are also examples of PQR-inclusions. 

\begin{ex}\label{ex:DiscreteRealization}
    Recall from Example \ref{ex:HilbgammOuterAction} that an (outer)  action of a countable discrete group $\Gamma$ yields a (full) unitary tensor functor $\fdHilb(\Gamma)\to\fgpBim(A).$ 
    The usual reduced crossed product $A\rtimes_{r,\alpha}\Gamma$ corresponds to $A\rtimes_{r,\alpha}\bbC[\Gamma]$, and so these all constitute examples of irreducible C*-discrete inclusions, where $\bbC[\Gamma]\in \fdHilb(\Gamma)$ is the group C*-algebra object from Example \ref{ex:GroupAlg}.
    
    We can also capture the construction of \cite{MR1467459} by choosing $\Gamma=\bbZ,$ $\cC=\fdHilb(\bbZ),$ and the connected C*-algebra object $\bbB=\bbC[\bbZ],$ and recover the reduced crossed product by a single dualizable bimodule. 
\end{ex}

\begin{ex}\label{ex:algebrasSUq2}
    Recall Example \ref{ex:TKAlgObj} and the connected W*-algebra objects in  $\Vec(\Rep(\bbG))$, where $\bbG$ is a compact quantum group. By connectedness, these algebras are locally finite dimensional and so we can regard them as C*-algebra objects.  Then, given any outer action of $\bbG$ on a unital C*-algebra $A$ with trivial center (eg \cite{MR4139893}) and a C*-algebra $\bbB\in \Vec(\Rep(\bbG)$ we can use our crossed product construction to obtain a family of irreducible C*-discrete inclusions. 
    
    In case $\bbG=SU_q(2)$ for $q\in[-1,1]\setminus\{0\},$ then $\Rep(SU_q(2)) \cong \mathsf{TLJ}(\delta)$ is the Temperley-Lieb-Jones category of non-crossing partitions with loop parameter $\delta=-\textsf{sign}(q)(q+1/q)$. 
    As mentioned in Example \ref{ex:TKAlgObj} the connected C*-algebra objects are classified by \emph{$\delta$-fair and balanced graphs}. 
\end{ex}

\subsection{Realizations of \texorpdfstring{$\cC$}{cC}-graded vector spaces as  C*-correspondences}
Here, we deal with the problem on how to obtain a right $A$-$A$ Hilbert C*-correspondence $\cX$ from the data of a UTC-action $F:\cC\to\fgpBim(A)$ and a $\cC$-graded vector space $\bbX\in\Vec(\cC)$. 
Once we have constructed this correspondence from abstract data, in the case where $\bbX=\bbB$ is a connected C*-algebra object in $\Vec(\cC),$ we compare it with the canonical $A$-$A$ right C* correspondence $\cB$ one gets from the reduced crossed product C*-algebra equipped with its canonical conditional expectation $(A\rtimes_{F,r}\bbB, E),$ with the notation from Subsection \ref{subsec:CrossProd}.

\begin{construction}\label{construction:CorrespRealization}
    Given $\bbX\in\Vec(\cC)$ and an action $F:\cC\to\fgpBim(A)$ on a unital C*-algebra $A$, we define the realization/crossed product of $\bbX$ under $F$ by 
    $$\cX_F:=A\rtimes_F \bbX:=\overline{(A\rtimes_F \bbX)^{\diamondsuit}}^{\|\cdot\|_A}:= \overline{\bigoplus_{c\in\Irr(\cC)}}\bbF(c)\otimes \bbX(c)\in\rCorr(A\to A).$$
    (See Section~\ref{sec:modules_correspondences_bimodules} for details on orthogonal direct sums for C*-correspondences.)
\end{construction}

In particular, we will be interested in the following orthogonal direct sum right C*-correspondence:
\begin{lem}\label{lem:CorrRealization}
    Let $F:\cC\to\fgpBim(A)$ be an action of $\cC$ on a unital C*-algebra $A$ with trivial center,  and $\bbB\in\Vec(\cC)$ a C*-algebra object whose fibers are finite dimensional.
    There is a canonical unitary isomorphism of $A$-$A$ right C*-correspondences 
    $$\cB_F:= \overline{\bigoplus_{c\in \Irr(\cC)}} \bbF(c)\otimes \bbB(c) \cong \overline{\left(A\rtimes_{F,r}\bbB\right)\Omega}^{\|\cdot\|_A}=\cB.$$
\end{lem}
\begin{proof}
    For each $c\in\Irr(\cC),$ let $\{u(c)_i\}_{i\in I_c}$ be a finite right Pimsner--Popa basis for $F(c),$ and let $\{e(c)_j\}_{j=1}^{\dim(\bbB(c))}$ be a basis for $\bbB(c).$
    Then $\{u(c)_i\otimes e(c)_j\}_{i,j}$ is a right Pimsner--Popa basis for $\bbF(c)\otimes \bbB(c).$
    Now, clearly, the collection $\{u(c)_i\otimes e(c)_j\}_{c,i,j}$ is a (countable) generalized Pimsner--Popa basis in the sense of \cite{MR2085108} for $\overline{\bigoplus_{c\in \Irr(\cC)}} \bbF(c)\otimes \bbB(c).$
    Therefore, $\{u(c)_i\otimes e(c)_j\Omega\}_{c,i,j}$ is a generalized basis for $\overline{\left(A\rtimes_{F,r}\bbB\right)\Omega}^{\|\cdot\|_A}$ by Construction \ref{subsec:CrossProd}. 
    And the obvious mapping between these sets clearly extends to an $A$-$A$ unitary. 
    The details are routine, so we omit them. 
\end{proof}

\section{Characterization of C*-discrete inclusions}\label{sec:Characterization}
In this section we will prove the main theorem of this manuscript which abstractly characterizes those irreducible inclusions of C*-algebras that arise as crossed products by an outer action of a UTC over a  unital C*-algebra with trivial center. 
The theorem constitutes a C*-analogue of \cite[Theorem 5.35]{MR3948170}, where the authors characterize discrete subfactors in terms of $\cC$-graded W*-algebra objects. 
Our theorem and its proof concretely illustrate the feasibility and usefulness of adapting techniques from subfactors into the study of C*-algebras. 
For example, in Section \ref{sec:EgsApplications} using this result we will introduce later on a notion of a \emph{standard invariant} for irreducible C*-discrete inclusions.

\begin{defn}\label{defn:CStarDisc}
    Given a unital C*-algebra $A$ with trivial center, we define the category $\CDisc,$ whose objects are precisely unital irreducible C*-discrete inclusions $A\overset{E}{\subset} B,$ whose morphism spaces are defined as     
    $$\CDisc\!\left((A\overset{E}{\subset} B)\to (A\overset{F}{\subset} C)\right)
    :=\set{\psi:B\to C}{\ F\circ\psi = E,\ \psi\text{ is } A\text{-}A \text{ bimodular ucp map}}\!.$$    
\end{defn}

\begin{defn}\label{defn:ConCstaralg}
    Given a unitary tensor category $\cC$, let $\cCCAlgs$ be the category of connected $\cC$-graded C*-algebra objects consisting of  C*-algebra objects in $\Vec{(\cC)},$ whose morphisms are categorical ucp morphisms; \emph{i.e.} $*$-natural transformations $\theta:\bbA\Rightarrow\bbB$ such that for every $c\in\cC$ we have $\theta_{\overline{c}\otimes c}:\bbA(\overline{c}\otimes c)\to \bbB(\overline{c}\otimes c)$ is a positive operator with $\theta_{1_\cC}(\bbA^1)= \bbB^1.$
\end{defn}

Recall that we say a $\cC$-graded C*-algebra object $\bbA$ is connected if $\bbA(1_\cC)\cong \bbC$ (see Definition~\ref{defn:*strucConnected}). The main purpose of this section is to obtain the following theorem, which allows us to associate a \textbf{quantum symmetry}---an (outer) action $F\to\fgpBim(A)$ of a UTC on a simple unital C*-algebra $A,$ together with a chosen connected C*-algebra object in $\Vec(\cC)$---to a PQR-inclusion, which corresponds to the standard invariant of a subfactor:
\begin{thm}[{Theorem~\ref{thmalpha:Main}}]\label{thm:DiscCharacterization}
    Let $A$ be a unital C*-algebra with trivial center, and $F:\cC\to \fgpBim(A)$ be a fully-faithful unitary tensor functor (i.e. an outer action). We then have a categorical equivalence: 
    \begin{align*}
    &\left\{
        \begin{aligned}
            &\left(A\overset{E}{\subset}B\right)\in\CDisc\ \\
        \end{aligned}
         \middle|\
        \begin{aligned}
             & A'\cap B =\bbC1\\
              &\ \  \cC_{A\subset B}\subseteq F[\cC]
        \end{aligned}
    \right\}
     \cong\cCCAlgs.
    \end{align*}
\end{thm}
\noindent Here, $\cC_{A\subset B}$ is the UTC generated by the dualizable $A$-$A$ sub-bimodules of $\cB$ (see also Definition~\ref{defn:StdInv} below). To prove this statement, we will adapt the strategy followed by \cite{MR3948170} by extending the C*-algebra realization construction from Section \ref{sec:ReducedRealization} to a functor, and constructing its inverse, the \textit{underlying C*-algebra object} functor. 
Theorem \ref{thm:DiscCharacterization} should be compared with \cite[Cor 6.3]{MR1900138} and also with \cite[Lemma 3.8]{MR1622812}, where the fibers of $\bbB$ play the role of their \emph{systems of subspaces}.

\begin{remark}
    Notice that since here we only consider $\cC$-graded C*-algebra objects with finite dimensional fibers, being finite-depth for an inclusion $A\subset A\rtimes_{F,r}\bbB$ is  equivalent to the $\cC$-support of $\bbB$ (the subset of $c\in\Irr(\cC)$ with $\bbB(c)\neq 0$) being finite as well as finite (Watatani/Pimsner--Popa) index. In this case, $\bbB$ has the structure of a Frobenius algebra, and by dividing by the square root of the Watatani index can be normalized to a Q-system. (See \cite{MR4419534} and \cite[Theorem A]{MR4079745}.)
    For a general choice of connected $\bbB\in\rCorr(\cC)$, our crossed product construction will yield infinite-depth inclusions of infinite (Watatani/Pimsner--Popa) index.
\end{remark}

\subsection{The reduced generalized crossed product functor} 
In this section, given a fixed outer action action $\cC\overset{F}{\curvearrowright}A$ (i.e. a fully faithful unitary tensor functor $F:\cC\to\fgpBim(A)$), we show how to extend the reduced crossed product by a $\cC$-action construction from Section \ref{sec:ReducedRealization} to a functor $A\rtimes_{F,r}-:\cCCAlgs\to\CDisc.$ (We sometimes write just $A\rtimes-$ to save space if no ambiguity arises.) We mainly outline the exposition from \cite[\S 5.4]{MR3948170} for the reader's convenience, but also to make the necessary adjustments to adapt the construction to C*-algebras. Notice that here we make no reference to a specific choice of state on $A$ whatsoever. 
\begin{defn}\label{defn:RealizationFunctor}
We define the realization functor with respect to $F$ as the bounded extension of: 
\begin{align*}
    A\underset{F,r}{\rtimes}-:\cCCAlgs&\longrightarrow\CDisc\\
                \bbB&\mapsto \left(A \overset{E}{\subset}A\underset{F,r}{\rtimes}\bbB\right)\\
                \left(\bbB\overset{\theta}{\Rightarrow}\bbD\right)&\mapsto \left( A\underset{F,r}{\rtimes}\bbB\overset{A\rtimes\theta}{\to} A\underset{F,r}{\rtimes}\bbD\right)\\
                &\hspace{.5cm} \eta_{(c)}\otimes g_{(c)}\mapsto \eta_{(c)}\otimes \theta_c(g_{(c)}).
\end{align*}
\end{defn}
The expression above tells us how to evaluate $A\rtimes\theta$ on $(A\underset{F}{\rtimes}\bbB)^\diamondsuit$:
    \[
        (A\rtimes\theta)\left( \sum_{c\in\Irr(\cC)}\eta_{(c)}\otimes g_{(c)}\right) =    \sum_{c\in\Irr(\cC)} \eta_{(c)}\otimes\theta_c(g_{(c)}), 
    \]
and note that the latter element belongs to $(A\rtimes_F\bbD)^\diamondsuit$.

We now sketch how to make sense of $A\rtimes\theta(b)$ for arbitrary $b\in B$. Since $\bbD$ is connected, each $\bbD(c)$ for $c\in\Irr(\cC)$ comes canonically equipped with the inner product $\langle f_{(c)}|  f_{(c)'}\rangle_c= \bbD^2_{\bar{c}\otimes c}(j^\bbD_{c}(f_{(c)}\otimes f_{(c)}'))\in \bbC,$ turning each  fiber into a finite dimensional Hilbert space. The resulting \emph{$\cC$-graded Hilbert space} is denoted by $L^2\bbD\in\Hilb(\cC)$. The $\cC$-graded $\rm{W}^*$-algebra object $\bbL(L^2\bbD)$ of  \emph{$\cC$-graded bounded endomorphisms} of $\bbD$ is given fiber-wise by $\bbL(L^2\bbD)(c):=\Hilb(\cC)( \mathbb{c}\otimes L^2\bbD\to L^2\bbD),$ associated to the cyclic $\rm{W}^*$-subcategory of $\Hilb(\cC)$ generated by $L^2\bbD$. 
By a $\cC$-graded version of Stinespring Dilation Theorem \cite[Theorem 4.28]{JP17}, regarding $\theta$ as a ucp morphism $\theta:\bbB\Rightarrow \bbL(L^2\bbD),$ there is $\bbK\in\Hilb(\cC)$, a $*$-representation $\pi:\bbB\Rightarrow \bbL(\bbK)$ and an isometry $v: L^2\bbD\Rightarrow \bbK$ satisfying $\theta= \Ad(v)\circ \pi.$   

Similarly to \cite[Lemma 5.17]{MR3948170}, we can obtain an isometry $A\rtimes v\in\rCorr_{\bbC-A}(A\rtimes_F\bbD\to A\rtimes_F\bbK)$ extending the mapping $\sum_c{\xi_{(c)}\otimes f_{(c)}}\mapsto \sum_c{\xi_{(c)}\otimes (v\circ f_{(c)}})$ by noticing it preserves the $\|\cdot\|_A$-norms, commutes with the right $A$-action, and furthermore $(A\rtimes_Fv)^* = A\rtimes_F(v^*)$, which is similarly defined. Moreover, we can turn $\pi$ into an honest $\rm{C}^*$-algebra homomorphism $A\rtimes\pi: A\rtimes_{F,r}\bbB\to \End^\dag(A\rtimes_F\bbK_A),$ given by $(\pi\rtimes A)(\eta\otimes g)(\xi\otimes k):= [\bbF^2(\eta\boxtimes \xi)]\otimes [\pi(g)\circ(\id\otimes k)]$. Finally, at the algebraic level, on $(A\rtimes_F\bbB)^\diamondsuit$, one obtains $A\rtimes\theta = \Ad(A\rtimes v)\circ(A\rtimes\pi)$ (cf \cite[Proposition 5.19]{MR3948170}) which we now extend by continuity to all of $A\rtimes_{F,r}\bbB.$

By Theorem \ref{thm:CrossProdPQR} we know $A\underset{F,r}{\rtimes}-$ is valued in $\CDisc\subseteq \PQR$, and both $A\rtimes\id_\bbB = \id_{B}$ and $A\rtimes(\theta_1\circ\theta_2)=(A\rtimes\theta_1)\circ(A\rtimes\theta_2)$ hold trivially. 
That $A\rtimes\theta$ preserves the conditional expectations follows from a direct computation relying on the unitality assumption $\theta(\bbB^1) = \bbD^1.$

\subsection{The underlying C*-algebra object functor}\label{subsec:CAlgObjFun}
We shall now construct an inverse functor $\langle-\rangle_F$ to reduced generalized crossed product/realization by the outer action $\cC\overset{F}{\curvearrowright}A$. If no confusion may arise, we sometimes write simply $\langle-\rangle$ to save space. 
\begin{defn}\label{defn:UnderlyingAlg}
Given a UTC-action $F:\cC\to\fgpBim(A)$, 
we can consider the underlying C*-algebra functor defined by 
\begin{align*}
    \langle-\rangle_F: \CDisc&\longrightarrow \cCCAlgs\\
    \left(A\overset{E}{\subset}B\right)&\mapsto \langle B\rangle_F: \cC^{\op}\to\Vec,\ \text{ given by }\\
    &\hspace{2.1cm} c\mapsto \rCorr_{A-A}(F(c)\to \cB)\\
    &\hspace{0.6cm} \left(a\overset{\alpha}{\to} b\right)\mapsto \langle B\rangle (\alpha): \rCorr_{A-A}(F(b)\to \cB)\to \rCorr_{A-A}(F(a)\to \cB)\\
    &\hspace{8cm} f\mapsto f\circ F(\alpha),\\
    \left(B\overset{\psi}{\to}D\right)&\mapsto \langle\psi\rangle:\langle B\rangle_F\Rightarrow\langle D\rangle_F,\text{ is a natural transformation with components }\\ &\hspace{.6cm}\langle\psi\rangle_c:\langle B\rangle(c) \to \langle D\rangle(c)\\
    &\hspace{2.6cm} f\mapsto (\psi\circ \check{f})(\cdot)\Omega.
\end{align*}
(Recall that $\check{f}(\xi)\Omega = f(\xi)\in\cB,$ as described in Definition~\ref{defn:B-valuedmaps}.)
\end{defn}
Since all inclusions in $\CDisc$ are PQR (Definition \ref{defn:PQR}), in light of C*-Frobenius Reciprocity (Theorems \ref{thm:FR} and \ref{thm:fgpInsideB}), we can equivalently express:
\begin{equation}\label{eq:UnderlyingAlg}
        \rCorr_{A-A}(F(c)\to \cB) = \rCorr_{A-A}(F(c)\to \cB)^\diamondsuit
        \overset{\Psi_{F(c)}}{\cong}\rCorr_{A-B}\left(F(c)\underset{A}{\boxtimes} B \to B\right).
\end{equation}
Thus showing $\langle B\rangle_F$ defines a connected C*-algebra object corresponding to the cyclic left $\cC$-module C*-subcategory of $\rCorr(A\to B)$ generated by ${}_AB_B$. (Recall Example \ref{ex:GenericCStarAO}.)    
But also, for every $c\in\cC,$ $f[F(c)]\subset B\Omega,$ so for any $\xi\in F(c),$ as in Definition~\ref{defn:B-valuedmaps}, we denote by $\check{f}(\xi)$ the unique element in $B$ such that $f(\xi)=\check{f}(\xi)\Omega.$

We now explicity describe the $*$-algebra structure of $\langle B\rangle_F$ in the following lemmas: 
\begin{lem}
    Let $A$ be unital C*-algebra with trivial center and $F:\cC\to\fgpBim(A)$ be an outer action. 
    For a given irreducible unital inclusion $(A\subset B)\in\CDisc$ with $\cC_{A\subset B} \subset F[\cC],$ the connected $\cC$-graded vector space $\langle B\rangle_F\in\Vec(\cC)$ is an algebra object with tensorator 
\begin{align*}
    \langle B\rangle^2_{a,b}:\langle B\rangle_F(a)\otimes \langle B\rangle_F(b)&\to \langle B\rangle_F(a\otimes b)\\
    f\otimes g&\mapsto  \left(m_B\circ(\check f\otimes \check g)\circ (F^2_{a,b})^{-1}(\,\cdot\,)\right)\Omega.
\end{align*}
\end{lem}
\begin{proof}
    First notice that $\langle B\rangle_{a,b}^2$ is well-defined. Indeed, since $A\subset B$ is C*-discrete, it follows that it is projective quasi-regular, and by Theorem \ref{thm:fgpInsideB} it then follows the expression for $\langle B\rangle^2_{a,b}(f\otimes g)$ is well-formed. It is clear that $\langle B\rangle^2_{a,b}(f\otimes g)\in\rCorr_{A-A}(F(a\otimes b)\to \cB)^\diamondsuit,$ as $A$-$A$ bimodularity is clear, and since $F(a\otimes b)\in\fgpBim(A),$ then the map is $A$-compact and thus bounded.
    Naturality of $\langle B\rangle^2$ follows from Definition \ref{defn:UnderlyingAlg} and the naturality of $F^2$.
\end{proof}

\begin{lem}\label{lem:StarStruc}
The $\cC$-graded algebra object $\langle B\rangle_F\in\Vec(\cC)$ has a $*$-structure given by
\begin{align*}
    j^{\langle B\rangle}_c:\langle B\rangle(c)&\to \langle B\rangle(\overline{c})\\
    f&\mapsto \Phi_{F(\overline{c})}\left[(\ev_{F(c)}\boxtimes_A \id_B)\circ\left(\id_{F(\overline{c})}\boxtimes_A {(\Psi_{F(c)} f)}^\dag\right)\circ(\id_{F(\overline{c})}\boxtimes_A\id_B)\right].
\end{align*}
Here, we are using the natural isomorphisms $\Psi$ and $\Phi$ from Theorem \ref{thm:FR}, and $c\in\cC$ is arbitrary.
Graphically, (omitting the natural isomorphisms $F(\overline c) \cong \overline{F(c)}$,)  the $*$-structure is given by the composite:
\begin{align*}
    \langle B\rangle(c)\ni\tikzmath{
        \roundNbox{fill=white}{(0,0)}{.3}{.1}{.1}{$f$};
        \draw (0,-.30) -- (0,-.8);
        \draw (0,0.3) -- (0,.8);
        \node at(0,1){$\scriptstyle B\Omega$};
        \node at(0,-1){$\scriptstyle F(c)$};
    }\mapsto
    \tikzmath{
        \begin{scope}
            \clip[rounded corners=5pt] (-.7,-.8) rectangle (.6,.8);
            \fill[\BColor] (0,0) -- (0,.8) -- (.6,.8) --  (.6,0);
            \fill[\BColor] (.2,0) -- (.2,-.8) -- (.6,-.8) --  (.6,0);    
        \end{scope}
        \node at (0,1){$\scriptstyle B$};
        \draw (0,0.3) -- (0,.8);
        \roundNbox{fill=white}{(0,0)}{.3}{.1}{.1}{$\Psi f$};
        \draw (-.2,-.3) -- (-.2,-.8);
        \draw (.2,-.3) -- (.2, -.8);
        \node at (-.3,-1){$\scriptstyle F(c)$};
        \node at (.25,-1){$\scriptstyle B$};
    }\mapsto
    \tikzmath{
        \begin{scope}
            \clip[rounded corners=5pt] (-.7,-.8) rectangle (.8,.8);
            \fill[\BColor] (.2,0) -- (.2,.8) -- (.8,.8) --  (.8,0);
            \fill[\BColor] (0,.2) -- (0,-.8) -- (.8,-.8) --  (.8,.2);
        \end{scope}
        \node at (0,-1){$\scriptstyle B$};
        \draw (0,-0.3) -- (0,-.8);
        \roundNbox{fill=white}{(0,0)}{.3}{.2}{.2}{$\scriptstyle {\left(\Psi f\right)}^\dag$};
        \draw (-.2,.3) arc(0:180:.2) -- (-.6,-.8);
        \draw (.2,.3) -- (.2,.8);
        \node at (-.6,-1){$\scriptstyle F(\overline{c})$};
        \node at (.2,1){$\scriptstyle B$};
    }\overset{{\Phi_{F(\overline c)}}}{\mapsto}
    \tikzmath{
        \draw (0,-0.3) -- (0,-.5);
        \roundNbox{fill=white}{(0,0)}{.3}{.2}{.2}{$\scriptstyle {\left(\Psi f\right)}^\dag$};
        \draw (-.2,.3) arc(0:180:.2) -- (-.6,-.8);
        \draw (.2,.3) -- (.2,.8);
        \node at(0,-.6) {$\bullet$};
        \node at (-.6,-1){$\scriptstyle F(\overline{c})$};
        \node at (.2,1){$\scriptstyle B\Omega$};
    }=
    j_c^{\langle B\rangle}(f)\in\langle B \rangle (\overline c).
\end{align*}
For a map $f\in\rCorr_{A-A}(F(c) \to \cB)^\diamondsuit,$ we have $(\Psi f)^\dag\in\rCorr_{A-B}(B\to F(c)\boxtimes_A B)$ is given by 
$$(\Psi f)^\dag(b)= \left(\sum_{i=1}^n u(c)_i\boxtimes \check{f}(u(c)_i)^*\right)b,$$
where, as usual, $\{u(c)_i\}_1^n$ is a right Pimsner--Popa basis for $F(c)$.
\end{lem}
\begin{proof}
    This structure corresponds to the $\dag$-structure in the module C*-category corresponding to $\langle B\rangle_F$ as explained under Eqn. \ref{eq:UnderlyingAlg}, and Ex. \ref{ex:GenericCStarAO}. 
    The last assertion follows from direct computation. 
\end{proof}

To see that $\langle \psi\rangle$ is a categorical cp map, (since $\cC$ admits direct sums) it suffices to see that for each $c\in\cC,$ $\langle \psi\rangle_{\bar c\otimes c}:\langle B\rangle_F(\bar c\otimes c)\to \langle D\rangle_F(\bar c\otimes c)$ is a positive map connecting these C*-algebras. 
To do so, we compose the isomorphism from Equation \ref{eq:UnderlyingAlg} together with the usual Frobenius Reciprocity in $\fgpBim(A)$ to obtain an endomorphism $\hat f\in\End_{A-B}(F(c)\boxtimes_A B)$ from $f\in\rCorr_{A-A}(F(\overline{c}\otimes c)\to \cB)$. 
This allows us to work directly on the associated $\cC$-module cyclic C*-category corresponding to $\langle B\rangle_F$. 
Now in this picture, applying $\langle \psi\rangle_{\bar c \otimes c}$ to an element of the form $\hat{f}^\dag\circ \hat f$ yields (suppressing the natural isomorphisms $\overline{F(c)}\cong F(\bar c)$):
\begin{align*}
    \langle \psi\rangle_{\bar c\otimes c}:\ 
    \tikzmath{
        \node at (0.2,1.8){$\scriptstyle B$};
        \draw (0.2,0.9) -- (0.2,1.6);
        \draw (-1.1,-1.45) -- (-1.1,0.9) arc (180:0:.4);
        \node at(-1.5,0){$\scriptstyle F(\bar c)$};
        \roundNbox{fill=white}{(0,0.6)}{.35}{.1}{.1}{$\hat{f^\dag}$};
        \draw  (-.2,.25) -- (-.2,-.3);
        \draw (0.2,.25) --  (0.2,-.3);
        \node at(.35,0){$\scriptstyle B$};
        \draw (-.2,-1.5) -- (-.2,-.6);
        \node at (-.55,0) {$\scriptstyle F(c)$};
        \node at(-.55,-1.2){$\scriptstyle F(c)$};
        \draw (.2,-1.5)  -- (.2,-.6);
        \node at(.35,-1.2) {$\scriptstyle B$};
        \roundNbox{fill=white}{(0,-0.6)}{.35}{.1}{.1}{$\hat{f}$};
        \roundNbox{fill=white}{(-0.8,-1.8)}{.35}{.3}{.4}{${F^2_{\bar c,c}}^\dag$};
        \node at (0.2,-1.5) {$\bullet$};
        \draw (-.85,-2.15) -- (-.85,-2.6);
        \draw (-.75,-2.15) -- (-.75,-2.6);
        \node at(-.8,-2.8){$\scriptstyle F(\bar c\otimes c)$};
    }\Omega
    &\hspace{.3cm}\mapsto \hspace{.3cm}
    \tikzmath{
        \node at (0.2,2.4){$\scriptstyle D$};
        \draw (0.2,1.7) -- (0.2,2.2);
        \roundNbox{fill=white}{(0.2,1.55)}{.25}{.1}{.1}{$\psi$};
        \draw (0.2,0.9) -- (0.2,1.3);
        \draw (-1.1,-1.45) -- (-1.1,0.9) arc (180:0:.4);
        \node at(-1.5,0){$\scriptstyle F(\bar c)$};
        \roundNbox{fill=white}{(0,0.6)}{.35}{.1}{.1}{$\hat{f^\dag}$};
        \draw  (-.2,.25) -- (-.2,-.3);
        \draw (0.2,.25) --  (0.2,-.3);
        \node at (.35,0){$\scriptstyle B$};
        \draw (-.2,-1.5) -- (-.2,-.6);
        \node at (-.55,0) {$\scriptstyle F(c)$};
        \node at(-.55,-1.2){$\scriptstyle F(c)$};
        \draw (.2,-1.5)  -- (.2,-.6);
        \node at(.35,-1.2) {$\scriptstyle B$};
        \roundNbox{fill=white}{(0,-0.6)}{.35}{.1}{.1}{$\hat{f}$};
        \roundNbox{fill=white}{(-0.8,-1.8)}{.35}{.3}{.4}{${F^2_{\bar c,c}}^\dag$};
        \node at (0.2,-1.5) {$\bullet$};
        \draw (-.85,-2.15) -- (-.85,-2.6);
        \draw (-.75,-2.15) -- (-.75,-2.6);
        \node at(-.8,-2.8){$\scriptstyle F(\bar c\otimes c)$};
    }\Omega\in\rCorr_{A-A}(F(\bar c\otimes c)\to \cD)^\diamondsuit_+,
\end{align*}
which is positive and obviously $A$-$A$ bimodular since $\psi$ is completely positive.  
With the assumption $\psi(1_A)=1_A$ we now conclude $\langle \psi\rangle$ is a categorical ucp map.
Functoriality of $\langle-\rangle_F$ is obvious, as  $\langle\psi_2\circ\psi_1\rangle=\langle \psi_2\rangle\circ\langle \psi_2\rangle$ and $\langle\id_B\rangle = \id_{\langle B\rangle}$ hold trivially by definition. 
Alternatively, splitting the $F(\bar c\otimes c)$ string using $F^2$ and \emph{raising strings} via the Frobenius reciprocity in $\fgpBim(A)$, we see that $\langle \psi\rangle$ acts on $\End_{A-B}(F(c)\boxtimes_A B)$ by $(\id_{F(c)}\boxtimes\psi)\circ-.$

We now record the following lemma:
\begin{lem}
    Given $\psi\in\CDisc\left(\left(A\overset{E}{\subset} B\right)\to \left(A\overset{F}{\subset} D\right)\right)\!,$ it is possible to extend $\psi$ to an $A$-$A$ bimodular bounded (contractive) map $\psi:\cB\to \cD.$
    Moreover, $\langle\psi\rangle$ is a categorical ucp map, 
    and $\langle-\rangle_F$ is functorial. 
\end{lem}
\begin{proof}
    That $\psi$ extends to $\cB$ is a direct consequence of $F\circ\psi = E$ and Kadison's inequality: $\psi(b)^*\psi(b)\leq\psi(b^*b)$ for all $b\in B$. 
    The rest was proven above. 
\end{proof}
Thus, the underlying algebra object $\langle-\rangle_F$ gives a well-defined functor. 

\begin{remark}
    Observe that we could have considered a category consisting of PQR inclusions (Definition \ref{defn:PQR}) in order to define the underlying C*-algebra object functor. However, we will make use of the (potentially stronger) C*-discrete condition to show it is inverse to reduced generalized crossed product/realization. 
\end{remark}

\subsection{Proof Theorem \ref{thmalpha:Main}}\label{subsec:characterization}
We shall prove Theorem~\ref{thm:DiscCharacterization} in this section. 
To do so, we show the reduced generalized crossed product/realization functor from Definition \ref{defn:RealizationFunctor} and the underlying C*-algebra object functor from Definition \ref{defn:UnderlyingAlg} are mutual inverses.

\begin{defn}\label{defn:delta}
Let $A$ be a unital C*-algebra with trivial center. 
Consider the following family of densely-defined maps indexed by $(A\subset B)\in\CDisc$  
\begin{align*}
    \delta^{\diamondsuit}: A\underset{F,r}{\rtimes}\langle-\rangle_F &\Rightarrow \id_{\CDisc} \\
    \bigoplus_{c\in\Irr(\cC)} F(c)\otimes_{\bbC}\langle B \rangle(c)\owns \xi\otimes f&\overset{\delta^{\diamondsuit}_B}{\mapsto} \check{f}(\xi).
\end{align*}
Recall that for $\xi\in F(c),$ $\check{f}(\xi)\Omega = f(\xi)$ (Definition~\ref{defn:B-valuedmaps}).
\end{defn}

The component maps $\delta^{\diamondsuit}_B$ come from the \emph{algebraic realization} from Proposition~\ref{prop:PQNAlgRealization} and have co-domains $\PQR(A\subset B)$. In the following we will extend these to the corresponding reduced C*-completions, and convince ourselves that our notion of C*-discrete inclusions yields the right co-domain for the resulting natural isomorphism $\delta.$

\begin{lem}\label{lem:deltaproperties}
    Let $A$ be a unital C*-algebra with trivial center.
    For each $A\overset{E}{\subset}B\in\CDisc$, the map $\delta^{\diamondsuit}_B$ is an $A$-$A$ bimodular, unital $*$-algebra homomorphism preserving the conditional expectation. In particular, $\PQN(A\subset B)$ is unital $*$-subalgebra of $B$.
\end{lem}
\begin{proof}
    Unitality follows from $\delta^{\diamondsuit}_B(\Omega\otimes \langle B\rangle^1) = \check{\langle B\rangle^1}(\Omega) =1_A,$ and $A$-$A$ bimodularity is automatic from definition. 
    
    To see $\delta^{\diamondsuit}_B$ is a $*$-algebra homomorphism, by linearity, it suffices to check on elementary tensors. 
    Let $\xi_1\otimes f_1\in F(a)\otimes \langle B\rangle(a)$ and $\xi_2\otimes f_2\in F(b)\otimes \langle B\rangle(b)$ be arbitrary.
    We have 
    \begin{align*}
        \delta^{\diamondsuit}_B&\left((\xi_1\otimes f_1)\cdot(\xi_2 \otimes f_2) \right)\\
        &= \delta^{\diamondsuit}_B\left(\sum_{\substack{c\in\Irr(\cC)\\ \alpha\in \Isom(c,a\otimes b)}}F(\alpha^*)\circ(F^2_{a,b}(\xi_1\boxtimes \xi_2))\otimes \left(m_B\circ(f_1\otimes f_2\circ (F^2_{a,b})^{-1})\circ F(\alpha)\right)\right)\\
        &= \sum_{\substack{c\in\Irr(\cC)\\ \alpha\in \Isom(c,a\otimes b)}}\left[m_B\circ (\check{f_1}\otimes \check{f_2})\circ (F^2_{a,b})^{-1} \right] \circ(F(\alpha)\circ F(\alpha^*))\circ \left[F^2_{a,b}(\xi_1\boxtimes \xi_2) \right]\\
        &= m_B\circ (\check{f_1}\otimes \check{f_2})(\xi_1\boxtimes \xi_2)\\
        &= \check{f_1}(\xi_1)\cdot \check{f_2}(\xi_2).
    \end{align*}
For an arbitrary $\xi\otimes f\in F(c)\otimes \langle B\rangle(c)$ we compute: 
    $$\delta^{\diamondsuit}_B((\xi\otimes f)^*)\!=\!\left(\Phi\left(\Psi f\right)^\dag\right)\!\check{\,}\,(\bar \xi)\!=\!(\ev_{F(c)}\boxtimes \id_B)(\bar \xi\boxtimes (\Psi(f))^\dag(1))= \langle \xi\mid  f^\dag(\Omega)\rangle_A\!=\!\langle\check{f}(\xi)\mid  1 \rangle_B\!=\!(\check{f}(\xi))^*.$$
    Thus $\delta^{\diamondsuit}_B$ is a unital $*$-algebra homomorphism.
    
    Finally, if we denote the canonical conditional expectation $E':A\rtimes_{r,F}\langle B\rangle\twoheadrightarrow A$, we obtain that 
    \begin{align*}
        (E\circ \delta^{\diamondsuit}_B)\left(\sum_{\substack{c\in\Irr(\cC)}} \xi_{(c)}\otimes f_{(c)}\right) &= E\left(\sum_{\substack{c\in\Irr(\cC)}}\check{f}{(c)}(\xi_{(c)}) \right)\\
        &=\check{f}_{1_\cC}(\xi_{1_\cC}) = \xi_{1_\cC}\otimes f_{1_\cC} = E'\left(\sum_{\substack{c\in\Irr(\cC)}} \xi_{(c)}\otimes f_{(c)}\right).
    \end{align*}
    Where the equality $\check{f}_{1_\cC}(\xi_{1_\cC}) = \xi_{1_\cC}\otimes f_{1_\cC}$ follows from connectedness $\langle B\rangle(1_\cC) \cong \bbC\cong A'\cap B.$
\end{proof}

In the next lemma, the C*-discrete property for inclusions becomes a crucial tool for our correspondence, as opposed to the potentially weaker projective quasi-regularity:

\begin{lem}\label{lem:ExactReconstruction}
    The map $\delta^{\diamondsuit}_B$ extends to an $A$-$A$ bimodular, expectation preserving C*-algebra isomorphism $\delta_B: A\rtimes_{r,F}\langle B\rangle\to B.$
\end{lem}
\begin{proof}
    Since $\delta^{\diamondsuit}_B$ preserves the expectations by Lemma \ref{lem:deltaproperties}, it extends to an isometric map of right $A$-$A$ C*-correspondences (recall Proposition \ref{prop:PQNAlgRealization})
    $$\delta_\cB: \overline{A\rtimes_{r,F}\langle B\rangle \Omega}^{\|\cdot\|_A}\to \cB.$$
    The map $\delta_\cB$ has dense range by Lemma \ref{lem:CorrRealization} and since $(A\subset B)\in\CDisc$ implies projective quasi-regularity. 
    This is, by Proposition \ref{prop:PQNAlgRealization}, we know that $\bigoplus_{\substack{c\in\Irr(\cC)}}F(c)\otimes_\bbC \rCorr_{A-A}(F(c)\to \cB)^\diamondsuit \overset{\delta^{\diamondsuit}_B}{\cong} B^{\diamondsuit}\Omega:=\PQN(A\subset B)\Omega,$ which after deleting $\Omega$ becomes a dense $*$-subalgebra of $B$ since $(A\subset B)\in \CDisc$ by assumption.
    To show it is actually a unitary, it then suffices to check it is adjointable.
    Thus, on the $\|\cdot\|_A$-dense subset $B^{\diamondsuit}\Omega$, the inverse map  $(\delta_\cB)^{-1}$ satisfies: 
    $$B^{\diamondsuit}\Omega\ni b\Omega\overset{(\delta_\cB)^{-1}}{\longmapsto}\sum_{d\in\supp(b)}\eta_{(d)}\otimes g_{(d)}\in \overline{(A\rtimes_{F,r}\langle B\rangle_F) \Omega}^{\|\cdot\|_A}.$$
    By a direct computation we see that $(\delta_\cB)^{-1} = (\delta_\cB)^\dag$ densely, and thus $\delta_\cB$ is unitary. 
    
    To extend $\delta^{\diamondsuit}_B$ to the reduced C*-completion $A\rtimes_{r,F}\langle B\rangle,$ notice that 
    $$\Ad(\delta_\cB):\End^\dag_{\bbC-A}\left(\overline{A\rtimes_{r,F}\langle B\rangle \Omega}^{\|\cdot\|_A}\right) \to \End^\dag({_\bbC}\cB_A)$$ is a C*-algebra isomorphism taking $[A\rtimes_{r,F}\langle B\rangle]\rhd \subset \End^\dag_{\bbC-A}\left(\overline{A\rtimes_{r,F}\langle B\rangle \Omega}^{\|\cdot\|_A}\right)$ onto $[B]\rhd\subset \End^\dag_{\bbC-A}(\cB_A).$ 
    Indeed, since $\delta_B^{\diamondsuit}$ preserves products, $\Ad(\delta_\cB)((\xi\otimes f)\rhd)$ acts on $B^{\diamondsuit}\Omega$ as multiplication: $\check{f}(\xi)\cdot(-)\Omega.$ Here of course $\check{f}(\xi)\in B^{\diamondsuit}.$ 
    Thus, the map $\delta_B:=\Ad(\delta_\cB)|_{A\rtimes_{F,r}\langle B\rangle}$ implements the desired C*-algebra isomorphism extending $\delta^{\diamondsuit}_B.$ 
\end{proof}

\begin{remark}
    The existence of the expectation preserving $A$-$A$ bimodular C*-algebra isomorphism $\delta_B: A\rtimes_{r,F}\langle B\rangle\to B$ from Lemma \ref{lem:ExactReconstruction} allowed us to reconstruct the original C*-algebra $B$ exactly from the left $A\rtimes_{r,F}\langle B\rangle$-action on $\overline{A\rtimes_{r,F}\langle B\rangle \Omega}^{\|\cdot\|_A}$, up to unitary conjugation. And so, remembering the data of $E$---or the representation of $B$ acting on the left of $\cB$---squarely recovers the C*-algebra $B$ and not just a dense $*$-subalgebra. 
\end{remark}

We now verify naturality of $\delta.$ If $\psi\in\CDisc((A\subset B)\to(A\subset D)),$ then 
\begin{align*}
    A\underset{F,r}{\rtimes}\langle\psi \rangle: A\underset{F,r}{\rtimes}\langle B \rangle &\to A\underset{F,r}{\rtimes}\langle D \rangle\\
    \xi_{(c)}\otimes f_{(c)} &\mapsto \xi_{(c)}\otimes \psi(f_{(c)}).
\end{align*} 
Then for each $c\in\cC$, a direct computation gives us that  $(\psi\circ\delta_B)(\xi_{(c)}\otimes f_{(c)}) = \psi\left(\check{f}_{(c)}\left(\xi_{(c)}\right)\right) =(\delta_D \circ A\underset{F,r}{\rtimes}\langle\psi\rangle) (\xi_{(c)}\otimes f_{(c)}),$ for any $\xi_{(c)}\otimes f_{(c)}\in A\underset{F,r}{\rtimes}\langle B\rangle.$
Since these inclusions are C*-discrete, must hold on all $A\rtimes_{F,r}\langle B\rangle,$ and together with linearity, we have thus proven that $\delta$ defines a natural isomorphism. 

We shall now construct a natural isomorphism going the other way. 
This transformation was spelled out in \cite[\S 5.6]{MR3948170} in the context of subfactors and Hilbert spaces, and we outline it here in the context of C*-algebras and C*-correspondences: 
\begin{defn}\label{defn:kappa}
Consider the family of morphisms whose components are indexed by $\bbD\in\cCCAlgs$
\begin{align*}
    \kappa: \id_{\cCCAlgs}&\Rightarrow\langle A\underset{F,r}{\rtimes}-\rangle_F\\
    \kappa_\bbD:\bbD&\Rightarrow \langle A\underset{F,r}{\rtimes}\bbD\rangle_F.
\end{align*}
Were each $\kappa_\bbD,$ a morphism in $\cCCAlgs$ indexed by $c\in\cC,$ is given by 
$$
\tikzmath{
    \node at(0,.8){$\scriptstyle \bbD$};
    \draw (0,.3) -- (0,.6);
    \roundNbox{fill=white}{(0,0)}{.3}{.1}{.1}{$f$}
    \node at(0,-.8){$\scriptstyle\mathbb{c}$};
    \draw (0,-.3) -- (0,-.6);
}\hspace{.3 cm}
\overset{\kappa_{\bbD, c}}{\mapsto}\hspace{.3 cm}
\left(
    \bbD(c)\ni\xi\mapsto
    \tikzmath{
    \node at(0,1.7){$\scriptstyle\bbD$};
    \draw (0,1) -- (0,1.5);
    \roundNbox{fill=white}{(0,.7)}{.3}{.1}{.1}{$f$}
    \draw[dotted](-.6,0) --(.6,0);
    \draw (0,-.4)--(0,.4);
    \node at(0,-1.7){$\scriptstyle \bbF$};
    \draw (0,-1) -- (0,-1.5);
    \roundNbox{fill=white}{(0,-.7)}{.3}{.1}{.1}{$\xi$}
    \node at(-.15,.15){$\scriptstyle \mathbb{c}$};
    \node at(-.15,-.15){$\scriptstyle \mathbb{c}$};
}
\right)
\in A\underset{F,r}{\rtimes}\bbD.
$$
\end{defn}

\begin{lem}
    The family $\kappa = \{\kappa_\bbD\}_{\bbD\in\cCCAlgs}$ defines a natural isomorphism.
\end{lem}
\begin{proof}
    Naturality follows from a simple graphical computation using the permeability of the dotted membrane. (See \cite[\S 4.1]{JP17}). 
    
    To show each $\kappa_{\bbD,a}$ is an isomorphism we first notice that $\bbD(a)$ and $\langle A\underset{F,r}{\rtimes}\bbD \rangle(a)$ are isomorphic finite dimensional $\bbC$-vector spaces as shown by the following composite of natural isomorphisms (wrt $a\in\cC$): 
    \begin{align*}
        \bbD(a)&\cong \bigoplus_{c\in\Irr(\cC)}\cC(a\to c)\otimes\bbD(c)\overset{\heartsuit}{\cong} \bigoplus_{c\in\Irr(\cC)} \rCorr_{A-A}(F(a)\to F(c))\otimes \bbD(c)\\
        & \cong\rCorr_{A-A}\left(F(a)\to \bigoplus_{c\in\Irr(\cC) \cap\supp(a)} F(c)\otimes \bbD(c)\right)\\
        &\overset{\clubsuit}{\cong} \rCorr_{A-A}\left(F(a)\to \overline{\left(A\underset{F,r}{\rtimes}\bbD\right)\!\Omega}^{\|\cdot\|_A} \right)\\
        &\overset{\spadesuit}{=} \rCorr_{A-A}\left(F(a)\to \overline{\left(A\underset{F,r}{\rtimes}\bbD\right)\Omega}^{\|\cdot\|_A} \right)^\diamondsuit\cong \langle A\underset{F,r}{\rtimes}\bbD \rangle(a).
    \end{align*}
    Here, the first isomorphism follows from the decomposition of $a\in\cC$ into simples in $\cC,$ $a\cong\bigoplus_{c\in\Irr(\cC)}c^{\oplus m_c}$, where $m_c:=\dim_\bbC(\cC(a\to c)).$ Isomorphism $\heartsuit$ follows since $F$ is fully-faithful by assumption, and $\clubsuit$ is an application of Lemma \ref{lem:CorrRealization}. The last one, $\spadesuit,$ is Theorem \ref{thm:fgpInsideB} since $A\subset A\rtimes_{F,r}\bbD\in\CDisc.$
    
    Now observe that $\kappa_{\bbD,a}$ is injective, as it maps the $\Irr(\cC)$-graded components of vectors $F(a)\ni \xi=\sum_{c\in\Irr(\cC)}\xi_{(c)}$ into the direct summands of the amplification $\bbF(c)\otimes \bbD(c).$
    Then, each $\kappa_{\bbD,a}$ is automatically an isomorphism, and this completes the proof.
\end{proof}
In fact, by carrying out the above composite of isomorphisms graphically, one can see they yield exactly $\kappa_{\bbD,a}.$ 
We now show it has the right properties:

\begin{lem}\label{lem:KappaHom}
    The natural isomorphism $\kappa_\bbD$ is a unital homomorphism of $*$-algebra objects.
\end{lem}
\begin{proof}
    We first check multiplicativity. 
    On the one hand, given $f\otimes g\in \bbD(a)\otimes \bbD(b),$ the map $\kappa_{\bbD, a\otimes b}(\bbD^2_{a,b}(f\otimes g))$ is given by 
    \begin{align*}
        \bbF(a\otimes b)\ni F^2_{a,b}(\eta\boxtimes\xi)\mapsto
        \tikzmath{ 
            \node at(0,2.5){$\scriptstyle\bbD$};
            \draw (0,2) -- (0,2.3);
            \draw (-.6,.9)arc(180:0:.6);
            \roundNbox{fill=white}{(0,1.7)}{.35}{.15}{.15}{$\bbD^2_{a,b}$}
            \roundNbox{fill=white}{(-.6,.65)}{.25}{.1}{.1}{$f$}
            \roundNbox{fill=white}{(.6,.65)}{.25}{.1}{.1}{$g$}
            \node at(-.75,.15) {$\scriptstyle\mathbb{a}$};
            \node at(-.75,-.15) {$\scriptstyle\mathbb{a}$};
            \draw(-.6,-.35) -- (-.6,.35);
            \draw(.6,-.35) -- (.6,.35);
                \draw[dotted](-1,0) --(1,0);
            \roundNbox{fill=white}{(-.6,-.65)}{.25}{.1}{.1}{$\eta$}
            \roundNbox{fill=white}{(.6,-.65)}{.25}{.1}{.1}{$\xi$}
            \node at(.75,.15) {$\scriptstyle\mathbb{b}$};
            \node at(.75,-.15) {$\scriptstyle\mathbb{b}$};
            \draw (0,-1.5) -- (0,-1.8);
            \draw (-.6,-.9)arc(180:360:.6);
            \node at(0,-2.5){$\scriptstyle\bbF$};
            \draw (0,-2) -- (0,-2.3);
            \draw (-.6,-.9)arc(180:360:.6);
            \roundNbox{fill=white}{(0,-1.7)}{.35}{.15}{.15}{$\bbF^2_{a,b}$}
        }\in A\underset{F,r}{\rtimes}\bbD,
    \end{align*}
    which is equal to the map 
    $\langle A\underset{F,r}{\rtimes}\bbD\rangle^2_{a,b}\circ(\kappa_{\bbD,a}\otimes\kappa_{\bbD,b})$ on the nose.
    
    We now check that $\kappa_{\bbD,\overline{a}}\circ j^\bbD_a = j^{\scriptstyle\langle A\rtimes_F\bbD\rangle}_a\circ \kappa_{\bbD,a}.$
    For any arbitrary $g\in\bbD(a),$ and arbitrary $\overline{\eta}\in\bbF(\overline (a)),$ we have that, on the one hand, 
    \begin{align*}
        \left(\kappa_{\bbD,\overline{a}}\circ j^\bbD_a\right)(g)(\overline{\eta})
        =\tikzmath{
            \draw (0,1) -- (0,1.3);
            \node at (0,1.5){$\scriptstyle \bbD$};
            \roundNbox{fill=white}{(0,.7)}{.3}{.3}{.3}{$j^{\bbD}_a(g)$}
            \draw (0,.4) --(0,-.4);
            \node at(-.15,.15){$\scriptstyle \overline{\mathbb{a}}$};
            \draw[dotted](-.8,0) -- (.8,0);
            \node at(-.15,-.15){$\scriptstyle \overline{\mathbb{a}}$};
            \roundNbox{fill=white}{(0,-.7)}{.3}{.3}{.3}{$j^{\bbF}_a(\eta)$}
            \draw (0,-1) -- (0,-1.3);
            \node at(0,-1.5){$\scriptstyle \bbF$};
        }\Omega
        =\left(j^{\bbF}_a(\eta)\otimes j^{\bbD}_a(g)\right)\Omega
        = (\eta\otimes g)^*\Omega.
    \end{align*}
    On the other, by Lemma \ref{lem:StarStruc} we have
    \begin{align*}
        \left(j^{\scriptstyle\langle A\rtimes_{F,r}\bbD\rangle}_a\circ \kappa_{\bbD,a}\right)(g)(\overline\eta) &=\bigg(\bigg[(\ev_{F(a)}\boxtimes\id_{\cD_F})\circ\left(\id_{F(\overline c)}\boxtimes \big(\Psi(\kappa_{\bbD,a}(g))\big)^\dag\right)\bigg](\overline{\eta}\boxtimes 1)\bigg)\Omega\\
        &=(\ev_{F(a)}\boxtimes\id_{\cD_F})\left(\overline{\eta}\boxtimes\sum_{\substack{i=1}}^n u_i\boxtimes\left(\widetilde{\kappa_{\bbD,a}(g)}(u_i)\right)^* \right)\Omega\\
        &= \left[\sum_{\substack{i=1}}^n \langle\eta\mid  u_i \rangle_A\left(\widetilde{\kappa_{\bbD,a}(g)}(u_i)\right)^*\right]\Omega\\
        &= \left[\sum_{\substack{i=1}}^n \langle\eta\mid  u_i \rangle_A\left(j^{\bbF}(u_i)\otimes j^{\bbD}_a(g) \right)\right]\Omega\\
        &=\left[j^{\bbF}_a\left(\sum_{\substack{i=1}}^n u_i\lhd\langle u_i\mid  \eta\rangle_A\right)\otimes j^{\bbD}_a(g)\right]\Omega\\
        &= \left[j^{\bbF}(\eta)\otimes j^{\bbD}(g)\right]\Omega\\
        &= (\eta\otimes g)^*\Omega.
    \end{align*}
    Thus, $\kappa_{\bbD}$ is a $*$-algebra object map.
    
    Finally, that $\kappa_{\bbD, 1_\cC}(\bbD^1)= \langle A\rtimes_{F,r}\bbD\rangle^1$ follows from direct examination.
\end{proof}

Thus $\kappa$ is the claimed natural isomorphism, completing the proof of Theorem \ref{thm:DiscCharacterization}.

\section{Examples and some applications}\label{sec:EgsApplications}
In Section \ref{sec:EgsCDisc} we described C*-discrete inclusions arising from actions of discrete (quantum) groups and their reduced crossed products. We shall now describe abundant examples of C*-discrete infinite index inclusions arising from the the Guionnet-Jones-Shlyakhtenko construction.

\subsection{Discrete extensions of GJS C*-algebras}\label{sec:GJS}
We show how to produce families of irreducible C*-discrete inclusions from UTCs. All the ingredients one needs are an action of a unitary tensor category and its connected C*-algebra objects.

First, let $\cC$ be an arbitrary countably generated UTC. 
In \cite{MR4139893}, the authors construct a unital, simple, separable, exact, with stable rank 1, and monotracial GJS C*-algebra $A:=A_{\scriptstyle \cC},$ and a bi-involutive fully-faithful unitary tensor functor $$H:\cC \hookrightarrow \fgpBimtr(A_\cC).$$ 
Here, $\tr$ denotes the unique trace on $A$, and $\fgpBimtr(A)$ stands for the tensor subcategory of fgp $A$-$A$ bimodules $K$ which are \emph{compatible with the trace} \cite[Definition 5.7]{MR1624182}; i.e. $\tr\left( \langle \eta\mid  \xi\rangle_{A}\right)= \tr\left( {}_{A}\langle \xi, \eta \rangle\right)$ for all $\eta, \xi\in K.$ 
Let $N:=N_\cC:=A''\subset \cB(L^2(A, \sf{tr})),$ which in finite depth is known to be an interpolated free-group $\rm{II_1}$-factor \cite{MR3110503} and $L\bbF_\infty$ otherwise. 
Moreover, in \cite[\S 5]{MR4139893} Hartglass and the first-named author built a \emph{Hilbertification} fully-faithful unitary tensor functor 
$${}_{N}(-\boxtimes_{A}L^2(A))_{N}: \fgpBimtr(A)\hookrightarrow\spbfBim(N),$$ 
which recovers the fully-faithful unitary tensor functor 
$H'': \cC\hookrightarrow\spbfBim(N)$ constructed in \cite{BHP12}. 
(See \cite[8.5.36]{BlLM04} for a related construction.)

The second ingredient consisting of C*-algebra objects in UTCs can be obtained from discrete subfactors. 
By \cite{MR3948170}, we know that irreducible spherical discrete subfactors $L\subset M$ correspond to actions of its underlying supporting UTC of fgp bimodules $\cC_{L\subset M},$ together with its underlying connected W*-algebra object $\bbM\in \mathsf{W^*Alg}(\cC_{L\subset M})$.
But by connectedness, $\mathsf{W^*Alg}_{\mathsf{con}}(\cC_{L\subset M})=\rCorr_{\mathsf{con}}(\cC_{L\subset M}).$ 
So one can obtain C*-algebra objects in $\cC=\cC_{L\subset M}$ from subfactors.  
From this fixed UTC-action $\cC\overset{H}{\curvearrowright}A$ we can get a family of irreducible C*-discrete inclusions indexed by connected objects in $\rCorr(\cC)$ simply by taking reduced crossed products (Construction  \ref{const:realizedredC*alg}.)

\subsection{Standard invariants in \texorpdfstring{$\CDisc$}{CDisc}}
One of the strongest motivations to study the category of irreducible C*-discrete inclusions lies in the fact that it provides a sizable class of  inclusions which can be characterized by the \emph{standard invariant} from subfactor theory. 
In fact, each irreducible extremal discrete subfactor yields such a C*-discrete inclusion. 
    
Our approach to defining the standard invariant is a generalization of  M{\"u}ger's approach for finite depth and finite index $\rm{II}_1$-factors via unitary tensor categories and their Q-systems \cite{MR1966524}.  
As was spelled out in \cite{MR3948170}, one needs to replace Q-systems by more general $*$-algebra objects to include infinite depth and infinite index discrete subfactors. 

Taking our inspiration from \cite[\S 7.1]{MR3948170}, in this section we shall establish that inclusions $\CDisc$ admit a well-behaved standard invariant in Corollary \ref{cor:reconstruction}. 
Conversely, in Definition \ref{defn:AbsStdInv} we shall introduce the notion of an  \emph{abstract standard invariant}, which formally corresponds to an action of a unitary tensor category on a unital C*-algebra as well as a C*-algebra object, and show there is always a connected C*-discrete inclusion realizing the abstract dynamical data in Corollary \ref{cor:RedReconstruction}. 
These consequences of our main result Theorem \ref{thm:DiscCharacterization}, should be compared with Popa's analogous result for subfactors stating every abstract standard invariant can be realized on some discrete subfactor \cite{MR1055708, MR1334479,MR1339767}. 

\begin{defn}\label{defn:StdInv}
    Let $\left( A\overset{E}{\subset} B\right)\in\CDisc$ be connected. The \textbf{standard invariant} for this inclusion consists of:
    \begin{itemize}
        \item The unitary tensor category $\cC_{A\subset B}\subset \fgpBim(A)$ consisting of all  fgp (equivalently dualizable) Hilbert $A$-$A$ bimodules appearing as orthogonal summands in  ${}_A\cB_A,$ and
        \item The underlying connected $\cC$-graded C*-algebra object $\langle B\rangle\in\Vec(\cC)$ given by $\langle B\rangle(K):=\rCorr_{A-A}(K\to \cB)\cong\rCorr_{A-B}(K\boxtimes_A B\to B).$
    \end{itemize}
\end{defn}

As a direct consequence of Theorem \ref{thm:DiscCharacterization} and its proof, we conclude that connected C*-discrete inclusions are completely characterized by this data. This is, we can \emph{reconstruct} C*-discrete inclusions from its standard invariant:
\begin{cor}[{Corollary~\ref{coralpha:StdInvs}}]\label{cor:reconstruction}
    Let $\left( A\overset{E}{\subset} B\right)\in\CDisc$ be unital and irreducible. Then 
    $$\left(A\overset{E}{\subset} B\right) \cong \left( A\overset{E'}{\subset} A\underset{F,r}{\rtimes} \langle B\rangle\right).$$
    Here, $E'$ is the canonical faithful conditional expectation mapping elements in the crossed product onto the coefficient algebra $A$ and $F:\cC_{A\subset B}\hookrightarrow \fgpBim(A)$ is the inclusion. 
\end{cor}
\begin{proof}
    Directly from the characterization from Subsection \ref{subsec:characterization}, as  $\delta:A\underset{F,r}{\rtimes}\langle-\rangle_F \cong \id_{\CDisc}$.
\end{proof}

\begin{ex}\label{ex:StdInvGp}
    Consider a unital C*-algebra $A$ with trivial center and an outer action $\alpha$ of a discrete group $\Gamma.$ Let $A\subset B$ be the C*-discrete inclusion obtained by taking reduced crossed products. As described in Example \ref{ex:HilbgammOuterAction}, it is straightforward to promote $\alpha$ to a unitary tensor functor $F:\fdHilb(\Gamma)\to\fgpBim(A),$ and $B\cong A\rtimes_{F,r}\bbC[\Gamma]$ (Example \ref{ex:DiscreteRealization}) for the $\Gamma$ graded group C*-algebra object $\bbC[\Gamma]\in\Vec(\fdHilb(\Gamma))\cong \Vec(\Gamma)$ from Example \ref{ex:GroupAlg}.
    It is then clear that the standard invariant remembers $\Gamma$ and  $\alpha,$ and that $\langle B\rangle_F = \bbC[\Gamma].$
\end{ex}

Establishing a meaningful converse to this corollary requires some interpretation; this is, one could ask if it is possible to consider standard invariants abstractly. By extrapolating from the discrete group case, where one considers a C*-dynamical system as an (outer) action on a (unital simple) C*-algebra $A,$ we think of an abstract standard invariant as the data that encodes a \emph{C* quantum symmetry}. 
This is formally an action $\cC\overset{F}{\curvearrowright}A$ of a unitary tensor category---a (fully-)faithful unitary tensor functor $F:\cC\to\fgpBim(A)$---together with the algebraic data of a connected C*-algebra object $\bbB\in\rCorr(\cC)$ governing how the bimodules in the image of $F$ behave under multiplication.  
Recall we have already made sense of reduced crossed products in this context in Construction \ref{const:realizedredC*alg} and  will justify this nomenclature in Corollary \ref{cor:RedReconstruction}. 
Formally, we have the following Definition: 

\begin{defn}\label{defn:AbsStdInv}
    An \textbf{abstract standard invariant} consists of 
    \begin{itemize}
        \item a unitary tensor category $\cC,$ and
        \item a connected $\cC$-graded C*-algebra object $\bbB\in\rCorr(\cC)$ generating $\cC$ in the sense that for each $c\in\cC$ there is some $a\in\cC$ such that $c$ is a subobject of $a$ and $\bbB(a)\neq\{0\}.$
    \end{itemize}
\end{defn}

We now establish that abstract standard invariants always arise from C*-discrete inclusions, thus obtaining a \textbf{Reconstruction Theorem in the C*-discrete setting}: 
\begin{cor}[{Corollary~\ref{coralpha:RedReconstruction}}]\label{cor:RedReconstruction}
    Let $(\cC,\bbB)$ be an abstract standard invariant. Then, there exists $(A\subset B)\in\CDisc$ whose standard invariant is isomorphic to $(\cC,\bbB).$ Moreover, $A$ can be chosen to be separable, unital, simple, stable rank 1, monotracial, and exact.
\end{cor}
\begin{proof}
    In \cite{MR4139893}, the first author and Hartglass showed that for any UTC $\cC$ there is such  C*-algebra $A=A(\cC)$ and a fully-faithful bi-involutive strong monoidal functor $F:\cC\to\fgpBim(A);$ i.e. an outer action $\cC\overset{F}{\curvearrowright}A.$ 
    Since $\kappa: \id_{\cCCAlgs}\cong \langle A\rtimes_{F,r}-\rangle_F$ as shown in Subsection \ref{subsec:characterization}, it then follows that  $(\cC,\bbB)\cong(F[\cC],\bbB) \cong(\cC_{A\subset A\rtimes\bbB}, \langle A\rtimes_{F,r}\bbB\rangle_F).$
\end{proof}

\bibliographystyle{amsalpha}
\bibliography{bibliography}
\end{document}